\documentclass[a4paper]{article}

\usepackage[T1]{fontenc}
\usepackage{lmodern}
\usepackage{microtype}
\usepackage[a4paper,margin=24mm]{geometry}
\usepackage{amsmath,amssymb,amsfonts,amsthm,mathtools}
\usepackage{mathrsfs}
\newcommand{\mathds}[1]{\mathbb{#1}}
\usepackage{graphicx}
\usepackage{booktabs}
\usepackage{longtable}
\usepackage{tabularx}
\usepackage{array}
\usepackage{enumerate}
\usepackage{xcolor}
\usepackage{hyperref}

\renewcommand{\arraystretch}{1.2} 
\definecolor{primary}{RGB}{20,50,115}
\definecolor{accent}{RGB}{140,35,35}

\hypersetup{
    colorlinks=true,
    linkcolor=primary,
    citecolor=accent,
    urlcolor=primary,
    pdftitle={Projective Models of Genus 21 K3 Surfaces},
    pdfauthor={Junyu Meng}
}

\theoremstyle{plain}
\newtheorem{thm}{Theorem}[section]
\newtheorem{lemma}[thm]{Lemma}
\newtheorem{Corollary}[thm]{Corollary}
\newtheorem{proposition}[thm]{Proposition}

\theoremstyle{definition}
\newtheorem{definition}[thm]{Definition}

\theoremstyle{remark}
\newtheorem{remark}[thm]{Remark}

\newcommand{\she}{\mathscr{E}}
\newcommand{\shE}{\mathcal{E}}

\newcommand{\shF}{\mathcal{F}}

\newcommand{\shh}{\mathscr{H}}

\newcommand{\sho}{\mathcal{O}}

\newcommand{\shv}{\mathscr{V}}

\newcommand{\shU}{\mathcal{U}}
\newcommand{\shV}{\mathcal{V}}

\newcommand{\WW}{\mathcal{W}}
\newcommand{\smallD}{\mathcal{D}}

\newcommand{\bz}{\mathbb{Z}}

\newcommand{\bc}{\mathbb{C}}

\newcommand{\bp}{\mathbb{P}}

\newcommand{\slf}{\mathfrak{sl}}

\DeclareMathOperator{\ev}{ev}

\DeclareMathOperator{\Pic}{Pic}
\DeclareMathOperator{\Aut}{Aut}

\DeclareMathOperator{\image}{im}

\DeclareMathOperator{\wt}{wt}

\DeclareMathOperator{\rk}{rank}
\DeclareMathOperator{\Sym}{Sym}
\DeclareMathOperator{\Ext}{Ext}
\DeclareMathOperator{\End}{\mathcal End}
\DeclareMathOperator{\Hom}{Hom}

\DeclareMathOperator{\ch}{ch}
\DeclareMathOperator{\td}{td}

\DeclareMathOperator{\Ker}{Ker}
\DeclareMathOperator{\coker}{coker}
\DeclareMathOperator{\length}{length}
\DeclareMathOperator{\Supp}{Supp}

\newcommand{\DD}{\mathbf{D}}

\newcommand{\RG}{R\Gamma}
\newcommand{\RHom}{R\!\Hom}
\newcommand{\sheafhom}{\mathscr{H}\kern -.5pt om}
\newcommand{\sheafext}{\mathscr{E}\kern -.5pt xt}
\newcommand{\sheafend}{\mathscr{E}\kern -.5pt nd}
\newcommand{\HHom}{\shh\!om}
\newcommand{\EExt}{\she\!xt}

\title{\textbf{Projective Models of Genus 21 K3 Surfaces}}
\author{Junyu Meng}
\date{}

\begin{document}
\maketitle

\begin{abstract}
We construct projective models for general primitively polarized K3 surfaces $(S,h)$ of genus $21$ (degree $h^2=40$) inside the 12-dimensional fine quiver moduli space $X = N(3;3,4)$ parameterizing stable representations of the $3$-Kronecker quiver with dimension vector $(3,4)$. 
As a corollary, we show that the corresponding moduli space $\mathscr{M}_{21}$ is unirational and prove O'Grady's generalized Franchetta conjecture for genus $21$.
Then we investigate the geometry of $(S,h)$ through its Fourier--Mukai partner $S' = M_{S}(4,h,5)$ and construct
 an explicit isomorphism $f : S^{[2]} \xrightarrow{\sim} S'^{[2]}$ together with a common closed embedding into the projectivization $\bp^{17}$ of the space of global sections of the normalized universal bundle, using the polarization $H_2 - 4\delta_2$. Finally, we identify the Brill--Noether locus $\mathfrak{B} = \{Z \in S^{[3]} : h^0(F \otimes I_Z) \ge 2\}$ scheme-theoretically with the $\bp^2$-bundle $Y' = \bp_{S'}((F')^*)$ and exhibit an isomorphism $S^{[3]} \xrightarrow{\sim} M_{S'}(3,h',6)$.
\end{abstract}


\section{Introduction and main results}\label{sec:intro}

\subsection{Historical background}

A central problem in the projective geometry of K3 surfaces is to describe explicit projective models for general polarized K3 surfaces $(S,h)$ of genus $g$, where $h^2 = 2g-2$, and to determine the birational geometry of the nineteen-dimensional coarse moduli space $\mathscr{M}_g$ of primitively polarized K3 surfaces of genus $g$. For $2 \le g \le 5$, the projective models of general polarized K3 surfaces are classical: they are either double covers of $\bp^2$ branched along sextic curves ($g=2$) or complete intersections in $\bp^g$ ($g=3,4,5$) \cite{SaintDonat}. For genus $g>5$, there are no longer enough parameters for complete intersections in projective spaces to sweep out the nineteen-dimensional moduli space $\mathscr{M}_g$.

In a foundational series of papers \cite{Mukai88,Mukai92,Mukai96,Mukai06,Mukai10}, Mukai introduced the \emph{vector-bundle method}, where he replaced the ambient projective space by a rational homogeneous variety $G/P$ (such as a Grassmannian, an orthogonal Grassmannian, a symplectic Grassmannian, or a $G_2$-homogeneous variety) and cut it by general global sections of some homogeneous vector bundle $\WW$ satisfying the Calabi--Yau adjunction condition $\det\WW \simeq K_{G/P}^{-1}$ and $\rk(\WW) = \dim(G/P) - 2$. Using this method, Mukai constructed projective models for general polarized K3 surfaces of genera
\[
g \in \{6, 7, 8, 9, 10, 12, 13, 18, 20\},
\]
thereby proving the unirationality of $\mathscr{M}_g$ in each of these cases (with $g=11$ treated via curves and Brill--Noether theory in \cite{Mukai96}).

However, as the genus $g$ increases, homogeneous spaces $G/P$ and the completely reducible homogeneous bundles become too rigid to accommodate the numerical constraints. To overcome this barrier, Mukai \cite{Muk16} turned to replacing homogeneous varieties with \emph{fine moduli spaces of quiver representations} in genus $g=16$ ($h^2=30$).  For the $m$-Kronecker quiver---having two vertices and $m$ arrows from the first vertex to the second---and a coprime dimension vector $(a,b)$, King \cite{King} used geometric invariant theory to produce a smooth projective Fano variety $N(m;a,b)$ (also denoted $K^m_{(a,b)}$) of dimension $mab - (a^2+b^2-1)$, Picard number $1$ and Fano index $m$ \cite{Fanoquivermoduli}. Crucially, $N(m;a,b)$ carries two normalized universal vector bundles of ranks $a$ and $b$ 
fitting into the universal quiver representation on this fine moduli space.

In \cite{Muk16}, Mukai realized a general polarized K3 surface of genus $16$ inside the twelve-dimensional quiver moduli space $N(4;2,3)$ as the zero locus of a general section of the vector bundle $\shE_2^{\oplus 2} \oplus \shF_2^{\oplus 2}$, where $\shE_2,\shF_2$ are duals of universal bundles of ranks 3 and 2 respectively. This is a further advancement of the vector-bundle method: over a general polarized K3 surface $(S,h)$ of genus $16$, there exists a unique slope-stable vector bundle $F_2$ of Mukai vector $(2,h,8)$, such that for any slope-stable bundle $[E_2]\in M_S(3,h,5)$, the space of bundle morphisms $\Hom(F_2^*,E_2^*)$ is four-dimensional, and the evaluation morphism $F_2^*\otimes \Hom(F_2^*,E_2^*)\rightarrow E_2^*$ is a stable family of quiver representations of dimension vector $(2,3)$, inducing a morphism from $S$ to the fine quiver moduli space $N(4;2,3)$ such that the pullbacks of $\shF_2,\shE_2$ are equal to $F_2,E_2$ respectively.
The rich internal geometry of the genus $16$ model---including the projective geometry of the $\bp^1$-bundle $\bp_S(F_2^*)$ and the identification of $S^{[2]}$ with some Debarre--Voisin hyper-K\"ahler fourfolds via the unique non-trivial Fourier--Mukai partner $S' = M_{S}(3,h,5)$---was recently developed in detail by Han \cite{Han25} and the author \cite{Meng16}.

Introducing the roles of quiver moduli spaces is consistent with the past constructions of projective models. Recall, for example, that on a general genus 8 K3 surface, there exists a unique stable vector bundle $\sho_S$ of Mukai vector $(1,0,1)$ and a unique stable vector bundle $F_3$ of Mukai vector $(4,h,2)$, such that the evaluation map $\sho_S\otimes\Hom(\sho_S,F_3)\rightarrow F_3$ is a stable family of quiver representations of the 6-Kronecker quiver of dimension vector $(1,4)$ (which is equivalent to the fact that $F_3$ is globally generated). By the universal property of the quiver moduli space $N(6;1,4)$, we get an induced morphism from $S$ to $N(6;1,4)\cong \mathrm{Gr}(2,6)$. Based on the result of Mukai \cite{Mukai88}, the map $S\rightarrow \mathrm{Gr}(2,6)$ is an embedding in general, allowing us to get lots of information of a general genus $8$ polarized K3 surface.

As a continuation of this philosophy, the present article deals with general polarized K3 surfaces of genus 21. More precisely, the purpose of the present article is twofold:
\begin{enumerate}[(1)]
\item First, we establish a new projective model for general primitively polarized K3 surfaces of genus $21$ ($h^2=40$) inside the twelve-dimensional Fano quiver moduli space
\[
X = N(3;3,4)
\]
of representations of the $3$-Kronecker quiver with dimension vector $(3,4)$, and deduce the unirationality of $\mathscr{M}_{21}$. This is related to the evaluation morphism $F^*\otimes \Hom(F^*,E^*)\rightarrow E^*$ over $S$, for the unique stable bundle $F$ with $v(F)=(3,h,7)$, and a general stable bundle $E$ with $v(E)=(4,h,5)$.
\item Second, we link the geometry of $S$  with its Fourier--Mukai partner $S'=M_S(4,h,5)$, which has naturally appeared in the construction of projective models. We construct an explicit isomorphism $S^{[2]} \xrightarrow{\sim} S'^{[2]}$ between the Hilbert squares together with a common projective embedding in $\bp^{17}$. We also prove that the Hilbert cube $S^{[3]}$ is isomorphic to $M_{S'}(3,h',6)$, whose non-locally-free locus is equal to the $\bp^2$-bundle $Y' = \bp_{S'}((F')^*)$.
\end{enumerate}

\subsection{Setup and statement of main results}

Let $X = N(3;3,4)$ be the fine moduli space of stable representations of the $3$-Kronecker quiver with dimension vector $(3,4)$ with respect to the canonical stability parameter $\theta=(4,-3)$. As we review in Section~\ref{sec:ambient}, $X$ is a smooth projective Fano twelvefold with $\Pic(X) = \bz. L$ and $K_X = L^{-3}$. There exist two normalized universal bundles $\mathcal U_3$ and $\mathcal U_4$ of ranks $3$ and $4$ over $X$ corresponding to the two vertices of the quiver. We denote their duals by $\shF = \mathcal U_3^*$ and $\shE = \mathcal U_4^*$, so that
\[
\rk(\shE)=4,\qquad \rk(\shF)=3,\qquad \det\shE = \det\shF = L.
\]
The universal representation can be expressed as $\shF^*\otimes W^*\rightarrow\shE^*$, where $W^*$ is the 3-dimensional space freely generated by the three arrows in the quiver. We define a rank-ten vector bundle on $X$ by
\[
\WW := \shF^{\oplus 2} \oplus \shE.
\]
Note that $\rk(\WW) = 2\cdot 3 + 4 = 10 = \dim X - 2$ and $\det\WW = L^3 = K_X^{-1}$. For a smooth 2-dimensional zero locus $S = Z(\sigma)$ of a global section $\sigma \in H^0(X,\WW)$, we write
\begin{equation}\label{eq:restriction-notation}
E = \shE|_S,\qquad F = \shF|_S,\qquad h = c_1(L|_S).
\end{equation}
Throughout the paper, calligraphic letters $\shE,\shF$ denote bundles on the ambient variety $X$, whereas ordinary letters $E,F$ denote their restrictions to $S$.

Our first main theorem establishes that this construction yields projective models for general primitively polarized K3 surfaces of genus $21$.

\begin{thm}\label{thm:main}
The zero locus of a general section $\sigma \in H^0(X,\WW)$ is a smooth K3 surface $S = Z(\sigma)$. The ample polarization $h = c_1(L|_S)$ and restricted bundles $E,F$ satisfy the following properties:
\begin{enumerate}[(i)]
\item $h^2 = 40$, and $h$ is primitive in $H^2(S,\bz)$.
\item The bundles $E$ and $F$ are $h$-slope-stable on $S$, with Mukai vectors
\[
v(E) = (4,h,5),\qquad v(F) = (3,h,7).
\]
\item The family of smooth surfaces $S=Z(\sigma)$ dominates the nineteen-dimensional moduli space $\mathscr{M}_{21}$ of primitively polarized K3 surfaces of genus $21$.
\end{enumerate}
In particular, a general primitively polarized K3 surface of genus $21$ occurs in this construction, the coarse moduli space $\mathscr{M}_{21}$ is unirational (Corollary~\ref{cor:unirationality}) and O'Grady's generalized Franchetta
conjecture holds for genus 21 (Corollary~\ref{franchetta}).
\end{thm}
There are two technical subtleties in proving Theorem~\ref{thm:main}:
\begin{itemize}
\item First, an ample class $h$ of self-intersection $h^2 = 40$ on a K3 surface could a priori be imprimitive; namely, we need to exclude the case $h = 2A$, where $A$ is a primitive polarization of degree $A^2 = 10$ (genus $6$).
\item Second, stability of the universal quiver representation on $X$ is a pointwise linear-algebraic condition and does not automatically imply slope-stability of the restricted bundles $E$ and $F$ on a codimension-ten zero locus $S \subset X$.
\end{itemize}
To overcome these obstacles without circularity, we first compute the cohomology of $E, F$ and the bundles of homomorphisms on every smooth 2-dimensional zero locus $S = Z(\sigma)$ by combining the Koszul resolution complex of $\sho_S$ with the technique of Teleman quantization. The latter utilizes the Harder--Narasimhan stratification of the unstable locus and the full exceptional collection on the auxiliary sixfold $Y_{2,3} = N(3;2,3)$ considered in \cite{Derquiver}.

\begin{proposition}\label{prop:cohom-summary}
For every smooth 2-dimensional zero locus $S = Z(\sigma)$, the restriction maps $H^0(X,\shE) \to H^0(S,E)$ and $H^0(X,\shF) \to H^0(S,F)$ are surjective with kernels of dimensions $1$ and $5$ respectively, and the cohomology groups of restricted bundles are as follows:
\begin{equation}\label{eq:cohom-summary}
\begin{array}{c|ccc}
\toprule
\text{Bundle on } S & h^0 & h^1 & h^2 \\
\midrule
E & 9 & 0 & 0 \\
F & 10 & 0 & 0 \\
\End(E) & 1 & 2 & 1 \\
\End(F) & 1 & 0 & 1 \\
E^*\otimes F & 3 & 0 & 0 \\
\bottomrule
\end{array}
\end{equation}
In particular, $E$ and $F$ are simple on every such smooth surface $S$, $E$ has isotropic Mukai vector $v(E)^2 = 0$ and $F$ is a spherical bundle with $v(F)^2 = -2$.
\end{proposition}

Using Proposition~\ref{prop:cohom-summary}, we prove in Proposition~\ref{prop:KS} that the Kodaira--Spencer map of the universal family of smooth 2-dimensional zero loci has maximal rank $19$ at every point, without assuming primitivity of $h$. By the local Torelli theorem, the family therefore contains members $S$ with rank one Picard groups $\Pic(S) = \bz. A$ on which $E$ and $F$ are globally generated. In Proposition~\ref{prop:no-imprimitive}, we use Kuleshov's existence theorem for spherical bundles \cite{Kuleshov} and an analysis of the Jordan--H\"older filtration to prove that a K3 surface of degree $A^2 = 10$ (with $\Pic(S)=\bz.A$) cannot carry a simple globally generated rank-four bundle with Mukai vector $(4,2A,5)$ and $h^0=9$. This rules out the undesired case $h = 2A$, establishes the primitivity of $h$ across the entire connected family and simultaneously yields $h$-slope-stability of $E$ and $F$ on a general member.

In the second half of the paper (Sections~\ref{sec:fm}--\ref{sec:hilb-schemes}), we study the intrinsic geometry of a general genus $21$ K3 surface $(S,h)$ via its Fourier--Mukai partner
\[
S' = M_{S}(4,h,5).
\]
Whereas the spherical rank-three bundle $F$ is rigid and uniquely determined by $(S,h)$, the isotropic rank-four bundle $E$ varies in the two-dimensional moduli space $S'$, which is also a K3 surface. Mukai's duality theorem \cite{Mukai2} equips $S'$ with a primitive polarization $h'$ of degree $(h')^2 = 40$ and a normalized universal bundle $\she$ on $S \times S'$ with $\det\she \simeq \sho_S(h) \boxtimes \sho_{S'}(h')$. More importantly, the sheaf $\she$ is universal in both directions, i.e. it also induces an identification $S \simeq M_{S'}(4,h',5)$ so that the roles of $S$ and $S'$ are completely interchangeable.

The varying family of spaces $\Hom(\she|_{S\times\{s'\}},F)$ ($s' \in S'$) glues together and forms the unique slope-stable spherical bundle $(F')^* = q_*(\she^*\otimes p^*F)$ on $S'$ with Mukai vector $v((F')^*) = (3,-h',7)$. The varying family of global sections $H^0(S,\she|_{S\times\{s'\}})$ ($s' \in S'$) glues together and forms the unique slope-stable spherical bundle $\mathscr{V}' = q_*\she$ on $S'$ of Mukai vector $v(\mathscr{V}') = (9,2h',9)$.
The universal bundle $\she$ and its rank-nine pushforward bundles $\mathscr{V} = p_*\she$ and $\mathscr{V}' = q_*\she$ directly govern an explicit isomorphism between $S^{[2]}$ and $S'^{[2]}$ as follows.

\begin{thm}[{= Proposition~\ref{isoHilbertSquare} and Theorem~\ref{closed-embedding}}]\label{thm:intro-hilb-square}
Let $(S,h)$ be a general primitively polarized K3 surface of genus $21$ and $(S',h')$ be its Fourier--Mukai partner $M_{S}(4,h,5)$. Put $\mathscr{V} = p_*\she$ on $S$, $\mathscr{V}' = q_*\she$ on $S'$, $\she_Z = q_*(\she|_{Z\times S'})$ for $Z \in S^{[2]}$ and let $e_Z : \mathscr{V}' \to \she_Z$ be the restriction morphism.
\begin{enumerate}[(i)]
\item For every $Z \in S^{[2]}$, the degeneracy locus $D_7(e_Z) \subset S'$ of $e_Z : \mathscr{V}' \to \she_Z$ defined by the $8\times 8$ maximal minors has length $2$. The assignment $f(Z) = D_7(e_Z)$ defines an isomorphism between hyper-K\"ahler fourfolds
\[
f : S^{[2]} \xrightarrow{\ \sim\ } S'^{[2]},
\]
 whose inverse is given by the same construction with the factors $S$ and $S'$ exchanged.
\item For each $Z \in S^{[2]}$, there is a unique nonzero section (up to scalar) $[\sigma_Z] \in \bp(H^0(S\times S',\she)) \simeq \bp^{17}$ vanishing on $Z\times S'$ (and simultaneously on $S\times f(Z)$). The assignment $Z \mapsto [\sigma_Z]$ defines a closed embedding $j : S^{[2]} \hookrightarrow \bp(H^0(S\times S',\she))$ such that 
\[
j^*\sho_{\bp^{17}}(1) \simeq \sho_{S^{[2]}}(H_2 - 4\delta_2).
\]
The embedding $j' : S'^{[2]} \hookrightarrow \bp(H^0(S\times S',\she))$ that is defined in the same way satisfies $j = j'\circ f$.
\end{enumerate}
\end{thm}

Whereas the geometry of $S^{[2]}$ is governed by the rank-four universal bundle $\she$ alone, incorporating the spherical rank-three bundle $F$ and the three-dimensional arrow spaces $\Hom(\she|_{S\times\{s'\}},F)$ naturally reveals interesting geometry of the six-dimensional Hilbert cube $S^{[3]}$.

\begin{thm}[{= Theorems~\ref{hilb:bn} and \ref{hilb:global-fm}}]\label{thm:intro-hilb}
Let $(S,h)$ be a general primitively polarized K3 surface of genus $21$ and $(S',h')$ be its Fourier-Mukai partner $M_{S}(4,h,5)$.
\begin{enumerate}[(i)]
\item Every nonzero bundle morphism $u: \she|_{S\times\{s'\}} \to F$ ($s' \in S'$) over $S$ is generically surjective and its degeneracy locus is a length 3 subscheme $Z\subset S$. This induces a closed embedding of the $\bp^2$-bundle $Y' = \bp_{S'}((F')^*)$ into the Hilbert cube $S^{[3]}$, identifying $Y'$ scheme-theoretically with the four-dimensional Brill--Noether locus
\[
\mathfrak{B} = \{Z \in S^{[3]} : h^0(S, F \otimes I_Z) \ge 2\}.
\]
\item For every $Z \in S^{[3]}$, the shifted Fourier--Mukai transform $K_Z := \Phi_{\she}(I_Z)[1]$ of $I_Z$ is a slope-stable torsion-free sheaf on $S'$ of Mukai vector $(3,h',6)$, and $\Phi_{\she}(-)[1]$ induces an isomorphism between hyper-K\"ahler sixfolds
\[
S^{[3]} \xrightarrow{\ \sim\ } M_{S'}(3,h',6).
\]
Under this isomorphism, the non-locally-free locus in $M_{S'}(3,h',6)$ is precisely $\mathfrak{B}\simeq Y'=\bp_{S'}((F')^*)$, whose points correspond to elementary modifications $0 \to K_Z \to F' \to k(s') \to 0$.
\end{enumerate}
\end{thm}

\subsection{Organization of the paper}

In Section~\ref{sec:preliminaries}, we review King's GIT construction of quiver moduli spaces, the Harder--Narasimhan and Hesselink stratifications, and Teleman's quantization theorem via the Cousin spectral sequence. We also supply a formula for restrictions of the normal bundles of unstable strata for future use. In Section~\ref{sec:ambient}, we describe the ambient twelvefold $X = N(3;3,4)$ and study the spaces of global sections and global generation properties of its universal bundles. We also establish the vanishing results needed on the auxiliary sixfold $Y_{2,3} = N(3;2,3)$. 
Section~\ref{sec:unstable-strata} classifies the eighteen unstable Harder--Narasimhan types of $X$ and performs the weight analysis on the three unstable strata corresponding to Harder-Narasimhan types $((1,1),(2,3))$, $((2,2),(1,2))$ and $((3,3),(0,1))$. Subsection~\ref{sec:restriction} applies the Koszul resolution to compute the cohomology of the restricted bundles on $S$. In subsections~\ref{degreePolarization}--\ref{sec:deformations}, we determine the Mukai vectors of bundles on $S$, prove Proposition~\ref{prop:cohom-summary} and show that the Kodaira--Spencer map has rank $19$. Subsection~\ref{sec:primitive} rules out an imprimitive polarization and proves Theorem~\ref{thm:main}, Corollary~\ref{cor:unirationality} and Corollary~\ref{franchetta}. Section~\ref{sec:fm} develops the Fourier--Mukai duality between $S$ and $S' = M_{S}(4,h,5)$, and establishes the slope-stability and global generation properties of several bundles constructed via Fourier-Mukai transform. Section \ref{sec:hilb-schemes} deals with the hyper-K\"ahler geometry of Hilbert squares and cubes of genus 21 K3 surfaces. Subsections~\ref{sec:length-two}--\ref{sec:projective} establish the isomorphism $f : S^{[2]} \xrightarrow{\sim} S'^{[2]}$ between the Hilbert squares and their common projective embedding in $\bp^{17}$ (Theorem~\ref{thm:intro-hilb-square}). Subsections~\ref{sec:p2bundle}--\ref{sec:isomorphicCube} study the Hilbert cube $S^{[3]}\cong M_{S'}(3,h',6)$ and prove Theorem~\ref{thm:intro-hilb}. Finally, Appendix~\ref{app:hn} records the weights of the actions of the one-parameter subgroups corresponding to the eighteen unstable Harder-Narasimhan types.
\paragraph{Acknowledgement}
The author would like to thank Laurent Manivel, Jieao Song, Grzegorz Kapustka, Micha\l{}  
Kapustka and Jaros\l{}aw Buczy\'nski for helpful discussions.
\paragraph{AI disclosure}
The question of investigating the locally complete intersection in $N(3;3,4)$ was proposed to the author by his advisor Laurent Manivel. An affirmative answer that these locally complete intersections are K3 surfaces of genus 21 was obtained by the author using exactly Teleman quantization, which only involves the strict inequality case in Lemma \ref{Teleman}. But the question of whether the induced polarization is primitive remained unsolved until ChatGPT 6 Astra came to help. ChatGPT 6 Astra also pointed out how to deal with the case $m_\lambda(\mathcal{T}) =\eta_\lambda$ in Teleman quantization (Lemma \ref{Teleman}), which reveals more cohomological information unknown to the author previously. Moreover, large language models corrected the mistake in the original proof about smoothness of $\mathfrak B$ in Theorem \ref{hilb:bn}. AI also assisted with auditing the arguments for correctness and polishing the exposition. The author takes full responsibility for the correctness of the contents in this article.

\section{Preliminaries on quiver moduli spaces and Teleman quantization}\label{sec:preliminaries}

\subsection{Classical GIT construction and universal bundles}\label{subsection: classical quiver moduli}

We begin by recalling several foundational facts concerning moduli spaces of quiver representations, including King's GIT construction and the universal representations; we refer to \cite{King,Reincke} for further details.

Throughout this article, $Q=(Q_0,Q_1)$ denotes an acyclic quiver (that is, a finite directed graph without oriented cycles). We write $i \in Q_0$ for a vertex and $[a:i\rightarrow j]\in Q_1$ for an arrow. Fixing a dimension vector $\underline{d} = (d_i)_{i\in Q_0} \in \mathds{N}^{Q_0}$, we define the affine representation space and the base change group by
\[
R_{\underline{d}} = \prod_{a:i\rightarrow j}\Hom(\mathds{C}^{d_i},\mathds{C}^{d_j}),\qquad
GL(\underline{d}) = \prod_{i\in Q_0}GL(\mathds{C}^{d_i}).
\]
The reductive group $GL(\underline{d})$ acts on $R_{\underline{d}}$ by change of bases. Since the diagonal scalar subtorus $\Delta = \{(\alpha\cdot\mathrm{Id}_{d_i})_{i\in Q_0} : \alpha\in \mathds{C}^*\} \subset GL(\underline{d})$ acts trivially on $R_{\underline{d}}$, the action factors through the quotient group $G_{\underline{d}} = GL(\underline{d})/\Delta$.

To construct the moduli space via Geometric Invariant Theory, we fix a stability parameter $\theta = (\theta_i)_{i\in Q_0} \in \mathds{Z}^{Q_0}$ satisfying $\sum_{i\in Q_0} d_i\theta_i = 0$. We associate to $\theta$ the character $\chi_{\theta} : GL(\underline{d}) \to \mathds{C}^*$ defined by
\[
\chi_{\theta}\bigl((g_i)_{i\in Q_0}\bigr) = \prod_{i\in Q_0}\det(g_i)^{-\theta_i},
\]
which vanishes on $\Delta$ and therefore descends to be a character of  $G_{\underline{d}}$. Viewing $\chi_{\theta}$ as a $G_{\underline{d}}$-linearization of the trivial line bundle over $R_{\underline{d}}$, we can form the GIT quotient $M = R_{\underline{d}} /\!\!/_{\chi_\theta} G_{\underline{d}}$. By King's criterion \cite{King}, GIT (semi)stability with respect to the linearization $\chi_\theta$ has the following representation-theoretic characterization.

\begin{definition}\label{def:king-stability}
A representation $V$ of $Q$ of dimension vector $\underline{d}$ is called \emph{$\theta$-semistable} if for every subrepresentation $W\subset V$, one has $\sum_{i\in Q_0}\dim(W_i)\theta_i \leq 0$, and \emph{$\theta$-stable} if strict inequality holds for every proper nonzero subrepresentation $0 \subsetneq W \subsetneq V$.
\end{definition}

We assume henceforth that $\underline{d}$ is \emph{$\theta$-coprime}, meaning that $\sum_{i\in Q_0}e_i\theta_i \neq 0$ for every subdimension vector $0 < \underline{e} < \underline{d}$. Under this assumption, $\underline{d}$ is indivisible, $\theta$-semistability coincides with $\theta$-stability and the stable locus $R_{\underline{d}}^{\theta\text{-st}} \subset R_{\underline{d}}$ is a principal $G_{\underline{d}}$-bundle over the GIT quotient $M$. Consequently, $M$ is a smooth projective variety.

Furthermore, $M$ is a fine moduli space admitting a universal quiver representation constructed as follows. Choose some integral vector $\underline{a} = (a_i)_{i\in Q_0} \in \mathds{Z}^{Q_0}$ such that $\sum_{i\in Q_0}a_i d_i = 1$, and let $L(\underline{a})$ denote the trivial line bundle on $R_{\underline{d}}$ equipped with the $GL(\underline{d})$-linearization given by $\chi_{\underline{a}}$. For each vertex $i\in Q_0$, let the trivial bundle $\mathds{C}^{d_i}$ carry the standard representation of the $i$-th factor $GL(\mathds{C}^{d_i})$. On the tensor product $\mathds{C}^{d_i} \otimes L(\underline{a})$, the diagonal scalar subtorus $\Delta$ acts with weight $1 - \sum_{j\in Q_0} a_j d_j = 0$, so $\mathds{C}^{d_i} \otimes L(\underline{a})$ descends to a vector bundle $\mathcal{U}_i$ of rank $d_i$ on $M$. The tautological $GL(\underline{d})$-equivariant maps $\bc^{d_i} \otimes L(\underline{a})\rightarrow\bc^{d_j} \otimes L(\underline{a})$ corresponding to arrows $a:i\to j$ likewise descend to bundle morphisms $\mathcal{U}_i \to \mathcal{U}_j$ on $M$, which assemble together to give a universal family of quiver representations. Note that the individual bundles $\mathcal{U}_i$ depend on the choice of normalization $\underline{a}$ up to a common line-bundle twist, so the universal family is also determined only up to a twist.

\paragraph{Kronecker quiver moduli.}
Let $a,b$ be coprime positive integers and $m \ge 3$ be an integer. The \emph{$m$-Kronecker quiver} has two vertices $1,2$ and $m$ arrows from $1$ to $2$. Its moduli space of stable representations of dimension vector $\underline{d}=(a,b)$ with respect to the canonical stability parameter $\theta=(b,-a)$ is denoted by $N(m;a,b)$ or $K^m_{(a,b)}$. Because $\gcd(a,b)=1$ and $\theta$ is coprime to $\underline{d}$, whenever $N(m;a,b)$ is nonempty, it is a smooth projective Fano variety of dimension $mab - (a^2+b^2-1)$, with Picard number $1$ and Fano index $m$ by \cite[Corollary~5.2]{Fanoquivermoduli}.

\subsection{Harder--Narasimhan and Hesselink stratifications}\label{subsec:HN-Hesselink}

To analyze the cohomology of descended vector bundles on $M$ via local cohomology on the representation space $R = R_{\underline{d}}$, we recall the Harder--Narasimhan stratification of the unstable locus $R^{\mathrm{us}} = R - R^{\mathrm{ss}}$ and its identification with the Hesselink stratification \cite{HNquiver,Hoskins,Kirwan}.

Fix the stability parameter $\theta$, for any nonzero dimension vector $\underline{e} \in \mathds{N}^{Q_0}$, we define its \emph{slope} by
\[
\mu(\underline{e}) = \frac{\sum_{i\in Q_0}e_i\theta_i}{\sum_{i\in Q_0}e_i}.
\]
For any nonzero quiver representation $W$, we use $\mu(W)$ to denote $\mu\bigl(  \underline{\dim} (W)  \bigr)$.
By assumption $\mu(\underline{d})=0$, so a representation $V$ of dimension vector $\underline{d}$ is $\theta$-semistable (resp.\ $\theta$-stable) if and only if $\mu(W) \le 0$ (resp.\ $\mu(W) < 0$) for every proper nonzero subrepresentation $W \subset V$. Note that even when $\mu(\underline{d})\neq0$, we can still define the notion of semistability by requiring $\mu(W)\leq\mu(\underline{d})$ for all subrepresentations $W\subset V$.

\begin{definition}\label{def:HN-filtration}
A filtration $0 = W^0 \subset W^1 \subset \cdots \subset W^l = V$ of a quiver representation $V$ is called a \emph{Harder--Narasimhan (HN) filtration} if each subquotient $W^s/W^{s-1}$ is semistable ($1 \le s \le l$) and their slopes are strictly decreasing:
\[
\mu(W^1/W^0) > \mu(W^2/W^1) > \cdots > \mu(W^l/W^{l-1}).
\]
\end{definition}

Every quiver representation $V$ admits a unique Harder--Narasimhan filtration \cite{Reincke,HNquiver}. Writing $\underline{d}^s = \underline{\dim}(W^s/W^{s-1})$, the tuple $\tau = (\underline{d}^1,\dots,\underline{d}^l)$ is called the \emph{Harder--Narasimhan (HN) type} of $V$; it satisfies $\sum_{s=1}^l \underline{d}^s = \underline{d}$, $\mu(\underline{d}^1) > \cdots > \mu(\underline{d}^l)$, and each $\underline{d}^s$ admits a semistable representation.

For each HN type $\tau$, the stratum
\[
S_{\tau} = \{V \in R_{\underline{d}} : \text{the HN type of } V \text{ is } \tau\}
\]
is a smooth, locally closed, $G_{\underline{d}}$-invariant subvariety of $R_{\underline{d}}$ \cite{HNquiver}. Fix a standard reference flag $0 = V_i^0 \subset V_i^1 \subset \cdots \subset V_i^l = \mathds{C}^{d_i}$ with $\dim(V_i^s/V_i^{s-1}) = \underline{d}_i^s$ for each $i \in Q_0$ and a compatible decomposition
\[
\bc^{d_i}=\bigoplus_{s=1}^l U_i^s,
\qquad
V_i^s=\bigoplus_{r=1}^s U_i^r.
\]
Let $P_\tau \subset G_{\underline{d}}$ be the parabolic subgroup preserving this flag, $L_\tau $ be its Levi subgroup, and $\mathfrak{u}_\tau$ be the Lie algebra of its unipotent radical. We define the \emph{attracting blade} $Y_\tau \subset S_\tau$ and the \emph{semistable center} $Z_\tau^\circ \subset Y_\tau$ by
\begin{align*}
Y_{\tau} &= \{V \in S_{\tau} : \text{the $s$-th step of the HN filtration of $V$ has underlying spaces } (V_i^s)_{i\in Q_0}\},\\
Z_{\tau}^\circ &= \{V \in Y_{\tau} : V \text{ splits as direct sum } \bigoplus_{s=1}^l W^s/W^{s-1} \text{ of HN subquotients with underlying spaces } U_i^s
\}.
\end{align*}
Then $S_\tau \simeq G_{\underline{d}} \times_{P_\tau} Y_\tau$, and the quotient stack $[S_\tau / G_{\underline{d}}]$ is isomorphic to $[Y_\tau / P_\tau]$. Moreover the natural map $Y_\tau \to Z_\tau^\circ$, sending a representation to the direct sum of all its HN subquotients, is an affine bundle morphism, whose fiber at a point of $Z_\tau^\circ$ consists of strictly upper-triangular block matrices.

On the GIT side, fixing the standard maximal torus $T \subset G_{\underline{d}}$ and the standard Weyl-invariant inner product $\|\cdot\|$ on the cocharacter lattice $\mathrm{Y}(T)$, the Hesselink--Kempf--Ness stratification \cite{Kirwan,Hoskins} partitions $R^{\mathrm{us}}$ into strata $S_{[\lambda]}$ indexed by conjugacy classes of primitive one-parameter subgroups $\lambda : \mathds{C}^* \to G_{\underline{d}}$ that minimize the normalized Hilbert--Mumford weight $\bar{\mu}(x,\lambda) = \langle \lambda, \chi_\theta \rangle / \|\lambda\|$ among all $\eta$ for which $\lim_{t\to 0}\eta(t)\cdot x$ exists. By \cite{Hoskins}, the Hesselink stratification of $R_{\underline{d}}$ coincides with the Harder--Narasimhan stratification: to an HN type $\tau = (\underline{d}^1,\dots,\underline{d}^l)$ with slopes $\mu_s = \mu(\underline{d}^s)$, we associate the primitive one-parameter subgroup $\lambda = \lambda_\tau$ whose $i$-th block is proportional to
\[
\lambda_i(t) = \mathrm{diag}\bigl(\underbrace{t^{\mu_1},\dots,t^{\mu_1}}_{\underline{d}_i^1},\,\underbrace{t^{\mu_2},\dots,t^{\mu_2}}_{\underline{d}_i^2},\,\dots,\,\underbrace{t^{\mu_l},\dots,t^{\mu_l}}_{\underline{d}_i^l}\bigr),\qquad t\in\mathds{C}^*, \quad i\in Q_0.
\]
More precisely, we need to clear the denominators of $\mu_s$ whenever necessary to obtain the corresponding primitive integral one-parameter subgroup $\lambda$. We will use the notation $S_{\lambda}=S_\tau $, $Y_\lambda=Y_\tau $, and $Z_\lambda^\circ = Z_\tau^\circ$. Now the affine bundle map $Y_\lambda\rightarrow Z_{\lambda}^\circ$ can be interpreted as sending $x\in Y_\lambda$ to $\lim_{t\rightarrow 0}\lambda(t).x$.

In our applications to $Y_{2,3} = N(3;2,3)$ and $X = N(3;3,4)$, we will care about the case where the HN filtration consists of two parts with $\tau = (\underline{d}^1,\underline{d}^2)$, where $\underline{d}^1 = (a_+,b_+)$, $\underline{d}^2 = (a_-,b_-)$ and $\mu(\underline{d}^1) > \mu(\underline{d}^2)$. Let us record here the $L_\lambda$-equivariant presentation of the restricted normal bundle $\mathcal{N}_\lambda := N_{S_\lambda/R}|_{Z_\lambda^\circ}$ for the $3$-Kronecker quiver.

For any point of $Z_\lambda^\circ$, the corresponding representation splits as $M_+ \oplus M_-$, where $M_+$ has source and target spaces $(C_+,D_+)$ of dimensions $(a_+,b_+)$ and $M_-$ has source and target spaces $(C_-,D_-)$ of dimensions $(a_-,b_-)$. Because $Y_\lambda$ consists of upper-triangular block matrices, the restriction of the normal bundle of $Y_\lambda$ inside $R$ to $Z_\lambda^\circ$ is canonically equal to the lower-left block $\Hom(C_+,D_-)^{\oplus 3}$. Since $S_\lambda = G_{\underline{d}} . Y_\lambda$ and $P_\lambda$ leaves $Y_\lambda$ invariant, quotienting by the infinitesimal action of the off-diagonal Lie algebra block $\Hom(C_+,C_-) \oplus \Hom(D_+,D_-)$ yields the $L_\lambda$-equivariant exact sequence
\begin{equation}\label{eq:normal-two-block-complex}
0 \longrightarrow \Hom(C_+,C_-)\oplus\Hom(D_+,D_-) \xrightarrow{\ \mathrm{tr}\ } \Hom(C_+,D_-)^{\oplus 3} \longrightarrow \mathcal{N}_\lambda \longrightarrow 0,
\end{equation}
where $\mathrm{tr}$ is given by composition with the arrow maps appearing in $M_+$ and $M_-$. The above sequence is actually exact: the kernel of $\mathrm{tr}$ is equal to the spaces of homomorphisms between quiver representations $\Hom_Q(M_+,M_-)$, which vanishes identically on $Z_\lambda^\circ$ because $M_+$ and $M_-$ are semistable representations with $\mu(M_+) > \mu(M_-)$. Consequently, $\mathcal{N}_\lambda$ has rank
\[
c_\lambda = \mathrm{codim}_R(S_\lambda) = 3a_+ b_- - a_+ a_- - b_+ b_-,
\]
and taking determinants in \eqref{eq:normal-two-block-complex} yields the $L_\lambda$-equivariant formula:
\begin{equation}\label{eq:normal-two-block-det}
\det\mathcal{N}_\lambda = (\det C_+)^{a_- - 3b_-}\otimes (\det D_-)^{3a_+ - b_+}\otimes (\det C_-)^{-a_+}\otimes (\det D_+)^{b_-}.
\end{equation}

\subsection{Teleman quantization}\label{subsection: Teleman quantization}

Let $V$ be a finite-dimensional representation of $G = G_{\underline{d}}$, let $\widetilde{\mathcal{T}} = \sho_R \otimes V$ be the corresponding $G$-equivariant vector bundle on $R$, and let $\mathcal{T}$ be the descended vector bundle on $M = R^{\mathrm{ss}}//G \simeq [R^{\mathrm{ss}}/G]$. To compute $H^q(M,\mathcal{T})$, we use the local cohomology exact triangle for the open immersion $[R^{\mathrm{ss}}/G] \hookrightarrow [R/G]$ with closed complement $[R^{\mathrm{us}}/G]$:
\begin{equation}\label{eq:unstable-localization}
\RG_{[R^{\mathrm{us}}/G]}([R/G],\widetilde{\mathcal{T}})
\longrightarrow \RG([R/G],\widetilde{\mathcal{T}}) 
\longrightarrow \RG(M,\mathcal{T})
\longrightarrow \RG_{[R^{\mathrm{us}}/G]}([R/G],\widetilde{\mathcal{T}})[1].
\end{equation}
Since $R$ is affine and $G$ is reductive, we have $\RG([R/G],\widetilde{\mathcal{T}})=\bigl(\bc[R]\otimes V\bigr)^G[0]$. As a result, for every $q \ge 1$, we have a canonical isomorphism
\begin{equation}\label{eq:positive-cohomology-supports}
H^q(M,\mathcal{T}) \simeq H^{q+1}_{[R^{\mathrm{us}}/G]}([R/G],\widetilde{\mathcal{T}}).
\end{equation}

We can order the finitely many unstable HN strata to obtain a chain of closed substacks $\emptyset = Z_{n+1} \subset Z_n \subset \cdots \subset Z_1 = [R^{\mathrm{us}}/G]$, whose locally closed differences are $Z_{p+1}- Z_{p+2} = [S_{\lambda}/G]$. Then the right hand side of the above isomorphism can be calculated using Cousin spectral sequence:
\begin{equation}\label{fundmentalspectral}
E_1^{p\geq0,q}=H^{p+q}_{Z_{p+1}/Z_{p+2}}([R/G],\widetilde{\mathcal T})=H^{p+q}_{Z_{p+1}-Z_{p+2}}([R/G]-Z_{p+2}, \ \widetilde{\mathcal T})\implies H^{p+q}
_{[R^{\rm us}/G]}([R/G], \ \widetilde{\mathcal T}).
\end{equation}
For convenience, we will replace $[R/G]-Z_{p+2}$ by $[R'/G]$, where $R'$ is any $G$-invariant open 
subset containing $S_\lambda$ as a closed subvariety. 

On a stratum $S_\lambda$ of codimension $c_\lambda := \mathrm{codim}_R(S_\lambda) > 0$, let $\mathcal{N}_\lambda = N_{S_\lambda/R'}|_{Z_\lambda^\circ}$ be the restricted normal bundle, and let $\mathcal{A}_\lambda = N_{Z_\lambda^\circ/Y_\lambda}$ be the normal bundle of $Z_\lambda^\circ$ inside $Y_\lambda$. Since $\lim_{t\to 0}\lambda(t)\cdot x \in Z_\lambda^\circ$ for $x\in Y_\lambda$, the $\lambda$-weights on $\mathcal{A}_\lambda$ and $\mathfrak{u}_\lambda$ are strictly positive, while the $\lambda$-weights on $\mathcal{N}_\lambda$ are strictly negative. We define the positive integer
\begin{equation}\label{eq:eta-definition}
\eta_\lambda = \wt_\lambda(\det\mathcal{N}_\lambda^*) > 0.
\end{equation}
Because $\lambda(\mathds{C}^*)$ is central in $L_\lambda$ and acts trivially on $Z_\lambda^\circ$, the restriction $\widetilde{\mathcal{T}}|_{Z_\lambda^\circ}$ decomposes into $L_\lambda$-equivariant subbundles $(\widetilde{\mathcal{T}}|_{Z_\lambda^\circ})_\beta$ according to the $\lambda$-weights $\beta$. We denote by $m_\lambda(\mathcal{T})$ the maximal $\lambda$-weight occurring in $\widetilde{\mathcal{T}}|_{Z_\lambda^\circ}$. We define the \emph{critical weight subbundle} on $Z_\lambda^\circ$ by
\begin{equation}\label{eq:critical-coefficient}
\mathcal{C}_\lambda(\mathcal{T}) = (\widetilde{\mathcal{T}}|_{Z_\lambda^\circ})_{\eta_\lambda} \otimes \det\mathcal{N}_\lambda.
\end{equation}
Note that the torus $\lambda(\mathds{C}^*)$ acts trivially on this bundle.

\begin{lemma}[{Teleman quantization}]\label{Teleman}
Suppose that $m_\lambda(\mathcal{T}) \le \eta_\lambda$ on the stratum $S_\lambda$, 
then $$H^i_{[S_\lambda/G]}([R'/G],\widetilde{\mathcal{T}}) = 0$$ for all $i < c_\lambda$, and for every $j \ge 0$ there is a canonical isomorphism
\begin{equation*}\label{eq:boundary-supported-isomorphism}
H^{c_\lambda+j}_{[S_\lambda/G]}([R'/G],\widetilde{\mathcal{T}}) \simeq H^j\bigl([Z_\lambda^\circ/L_\lambda],\,\mathcal{C}_\lambda(\mathcal{T})\bigr).
\end{equation*}
In particular, if the strict inequality $m_\lambda(\mathcal{T}) < \eta_\lambda$ holds, then $C_\lambda(\mathcal T)=0$ and $H^*_{[S_\lambda/G]}([R'/G],\widetilde{\mathcal{T}}) = 0$ in all degrees.
\end{lemma}

\begin{proof}
Let $\mathcal{J} \subset \sho_{R'}$ be the ideal sheaf of the smooth subvariety $S_\lambda \subset R'$. Local cohomology supported on $[S_\lambda/G]$ is computed by $\RG_{[S_\lambda/G]}(\widetilde{\mathcal{T}}) = \varinjlim_n \RHom(\sho_{R'}/\mathcal{J}^n, \widetilde{\mathcal{T}})$. There is a canonical filtration of $\sho_{R'}/\mathcal J^n$ by $\mathcal J^a/\mathcal J^{a+1}$, thus we would like to consider
$R\Hom(\mathcal J^a/\mathcal J^{a+1}, \widetilde{\mathcal T})$,
which is isomorphic to $$R\Hom\bigl(\Sym^a  N_{S_\lambda/R'}^*,\widetilde{\mathcal T}|_{[S_\lambda/G]}\otimes\det N_{S_\lambda/R'}[-c_\lambda]\bigr)=\RG\bigl([S_\lambda/G],\,\widetilde{\mathcal{T}}|_{S_\lambda}\otimes \det N_{S_\lambda/R'} \otimes \Sym^a N_{S_\lambda/R'}[-c_\lambda]\bigr)$$ by Grothendieck-Verdier duality. 
Since $S_\lambda \simeq G \times_{P_\lambda} Y_\lambda$, we have $[S_\lambda/G] \simeq [Y_\lambda/P_\lambda]$. Moreover, $Y_\lambda \to Z_\lambda^\circ$ is an affine bundle morphism with relative tangent bundle $\mathcal{A}_\lambda$, so pushing forward along $Y_\lambda \to Z_\lambda^\circ$ and filtering by polynomial degrees along the fibers produce subquotients on $[Z_\lambda^\circ/P_\lambda]$ of the form
\[
\widetilde{\mathcal{T}}|_{Z_\lambda^\circ} \otimes \det\mathcal{N}_\lambda \otimes \Sym^a\mathcal{N}_\lambda \otimes \Sym^b\mathcal{A}_\lambda^*[-c_\lambda],\qquad a,b \ge 0.
\]
To calculate the cohomology of the above vector bundle over $[Z_\lambda^\circ/P_\lambda]$, we use the spectral sequence $E^{r,q}=H^r_{P_\lambda}(H^q(Z_\lambda^\circ,-))\implies H^{r+q}([Z_\lambda^\circ/P_\lambda],-)$. Recall that $H^r_{P_\lambda}(-)=H^r(u_\lambda,-)^{L_\lambda}$, where $u_\lambda$ is the Lie algebra of the unipotent part of $P_\lambda$, and the $L_\lambda$-invariant part of the Lie algebra cohomology is calculated by the invariant parts of Lie algebra cochains:
\begin{equation}\label{eq:complete-weight-term}
\Bigl(H^{q-c_\lambda}\bigl(Z_\lambda^\circ,\,\widetilde{\mathcal{T}}|_{Z_\lambda^\circ}\otimes\det\mathcal{N}_\lambda\otimes\Sym^a\mathcal{N}_\lambda\otimes\Sym^b\mathcal{A}_\lambda^*\bigr)\otimes \Lambda^r\mathfrak{u}_\lambda^*\Bigr)^{L_\lambda},\qquad a,b,r \ge 0.
\end{equation}
Because $\lambda(\mathds{C}^*) \subset L_\lambda$ acts trivially on $Z_\lambda^\circ$, any $L_\lambda$-invariant piece in \eqref{eq:complete-weight-term} must have total $\lambda$-weight zero. By assumption, $\widetilde{\mathcal{T}}|_{Z_\lambda^\circ}\otimes\det\mathcal{N}_\lambda$ has $\lambda$-weights at most $m_\lambda(\mathcal{T}) - \eta_\lambda \le 0$, whereas $\mathcal{N}_\lambda$, $\mathcal{A}_\lambda^*$, and $\mathfrak{u}_\lambda^*$ all have strictly negative $\lambda$-weights. Consequently, whenever $a+b+r > 0$, the entire group \eqref{eq:complete-weight-term} has strictly negative $\lambda$-weight and its $L_\lambda$-invariant part vanishes. As a result, only the $(a,b,r)=(0,0,0)$ term survives, and within it only the $\lambda$-weight-zero summand $\mathcal{C}_\lambda(\mathcal{T}) = (\widetilde{\mathcal{T}}|_{Z_\lambda^\circ})_{\eta_\lambda}\otimes\det\mathcal{N}_\lambda$ contributes. Since $L_\lambda$ is reductive, the $L_\lambda$-invariant part in the above formula for $a=b=r=0$ is equal to $H^{q-c_\lambda}([Z_\lambda^\circ/L_\lambda],(\widetilde{\mathcal T}|_{Z_\lambda^\circ})_{\eta_\lambda}
 \otimes\det\mathcal N_\lambda)$.
 Now we can conclude that $H^q_{[S_\lambda/G]}([R'/G],\widetilde{\mathcal T})=H^{q-c_\lambda}([Z_\lambda^\circ/L_\lambda],(\widetilde{\mathcal T}|_{Z_\lambda^\circ})_{\eta_\lambda}
 \otimes\det\mathcal N_\lambda).$
\end{proof}
The above lemma has essentially appeared in \cite{Teleman,HalpernL}. An application of such Teleman quantization to quiver moduli spaces in the case where the strict inequality $m_\lambda(\mathcal{T}) < \eta_\lambda$ holds has also appeared in \cite{belmans2023rigidity}. We include the proof here for keeping explicit track of the case $m_\lambda(\mathcal{T}) = \eta_\lambda$.

\begin{remark}
    Combining \eqref{eq:positive-cohomology-supports}, the Cousin spectral sequence \eqref{fundmentalspectral}, and Lemma~\ref{Teleman}, we see that, in the case $m_\tau(\mathcal{T}) \le \eta_\tau$ holds for all unstable HN types, in order to prove the vanishing of a specific cohomology group $H^p(M,\mathcal{T}), p\geq1$, it suffices to verify on each stratum $S_\lambda$ where $m_\lambda(\mathcal{T}) = \eta_\lambda$ that
\begin{equation}\label{eq:boundary-degree-target}
H^{p+1-c_\lambda}\bigl([Z_\lambda^\circ/L_\lambda],\,\mathcal{C}_\lambda(\mathcal{T})\bigr) = 0.
\end{equation}

\end{remark}

\begin{remark}\label{gerbe}
    The calculation of cohomology of critical weight subbundles over the quotient stack $[Z_\tau^\circ/L_\tau]$ can be made more explicit in the case where the HN type consists of two blocks, i.e. $\tau = (\underline{d}^1,\underline{d}^2)$ with $\underline{d}^1+\underline{d}^2=\underline{d}$. Under such assumption, $Z_\tau^\circ=R^{\mathrm{ss}}_{\underline{d}^1}\times R^{\mathrm{ss}}_{\underline{d}^2}$ and the quotient of $L_{\tau}$ by $\lambda$ is equal to $G_{\underline{d}^1}\times G_{\underline{d}^2}$, where $\lambda$ still denotes the primitive 1-ps associated to $\tau$. Since the action of $\lambda$ on the critical weight subbundle $\mathcal{C}_\lambda(\mathcal T)$ is trivial, the bundle $\mathcal{C}_\lambda(\mathcal T)$ descends to a vector bundle over $[R^{\mathrm{ss}}_{\underline{d}^1}/G_{\underline{d}^1}]\times[ R^{\mathrm{ss}}_{\underline{d}^2}/G_{\underline{d}^2}]$ along the gerbe map
    $$[Z_\tau^\circ/L_\tau]\longrightarrow [R^{\mathrm{ss}}_{\underline{d}^1}\times R^{\mathrm{ss}}_{\underline{d}^2}/ G_{\underline{d}^1}\times G_{\underline{d}^2}].$$
    Consequently, we only need to calculate the cohomology of $\mathcal C_\lambda(\mathcal T)$ on the product of smaller quotient stacks. Note here the semistable loci $R^{\mathrm{ss}}_{\underline{d}^1}, R^{\mathrm{ss}}_{\underline{d}^2}$ are defined using slope-stability as in the paragraph before Definition \ref{def:HN-filtration}.
\end{remark}

\section{\texorpdfstring{The ambient moduli space $X = N(3;3,4)$ and the auxiliary space $Y_{2,3} = N(3;2,3)$}{The ambient moduli space X = N(3;3,4) and the auxiliary space Y(2,3) = N(3;2,3)}}\label{sec:ambient}

\subsection{\texorpdfstring{Geometry and universal bundles of $X = N(3;3,4)$}{Geometry and universal bundles over X = N(3;3,4)}}

Let $V_3, V_4$ and $W$ be complex vector spaces of dimensions $3, 4$ and $3$, respectively. The representation space of the $3$-Kronecker quiver with dimension vector $(3,4)$ is
\[
R = \Hom(V_3,V_4)\otimes W.
\]
Choosing a basis of the arrow space $W$ identifies a point of $R$ with a triple of linear maps $(M_1,M_2,M_3)$, where $M_i : V_3 \to V_4$. The reductive group $GL(V_3)\times GL(V_4)$ acts on $R$ by change of basis, with the diagonal scalar subgroup $\bc^*$ acting trivially; we set
\[
G = \bigl(GL(V_3)\times GL(V_4)\bigr)/\bc^*.
\]
With the canonical stability parameter $\theta=(4,-3)$, a representation in $R$ is stable if and only if the dimension vector $(a,b)$ of any proper nonzero subrepresentation satisfies
\begin{equation}\label{eq:king-inequality}
4a - 3b < 0.
\end{equation}
Since $\gcd(3,4)=1$, semistability coincides with stability, the $G$-action on $R^{\mathrm{st}}$ is free, and King's construction yields the smooth projective Fano twelvefold
\[
X = N(3;3,4) = R^{\mathrm{st}}/G,\qquad \dim X = 36 - (9 + 16 - 1) = 12.
\]

To construct the universal bundles on $X$, we consider the character
\[
\chi = \det(V_4)\otimes\det(V_3)^{-1}
\]
of $GL(V_3)\times GL(V_4)$. The diagonal scalar subgroup acts on the twisted representations
\begin{equation*}\label{eq:normalization}
V_3\otimes\chi^{-1},\qquad V_4\otimes\chi^{-1}
\end{equation*}
trivially and these two representations descend to the universal bundles $\mathcal U_3, \mathcal U_4$ on $X$. We define their duals by
\[
\shE = \mathcal U_4^*,\qquad \shF = \mathcal U_3^*.
\]
The universal family of quiver representations is now given by $\shF^*\otimes W^* \to \shE^*$.
The determinants of $\shE,\shF$ coincide with the generator $L$ of $\Pic(X)$ because the codimension of the unstable locus is at least two and we have an identification of $\Pic(X)$ with the character group $\chi(G)=\bz.[(\det V_4)^3\otimes(\det V_3)^{-4}]$ (see \cite{Fanoquivermoduli} for details).
\begin{equation*}\label{eq:L-character}
\det\shE = \det\shF = L,\qquad L \longleftrightarrow (\det V_4)^3\otimes(\det V_3)^{-4}.
\end{equation*}
The generator $L$ is actually ample by the argument below and the fact that $X$ is Fano (\cite{Fanoquivermoduli}).

 Since the (semi)stable locus $R^{\mathrm{ss}}$ is a principal $G$-bundle over $X$, we have the following Euler-type tangent sequence on $X$:
\begin{equation}\label{eq:tangent-X}
0 \longrightarrow \sho_X \longrightarrow \End(\shE)\oplus\End(\shF) \longrightarrow (\shE^*\otimes\shF)\otimes W \longrightarrow T_X \longrightarrow 0.
\end{equation}
Taking determinants and using $c_1(\shE^*\otimes\shF) = 4c_1(\shF) - 3c_1(\shE) = c_1(L)$ yields
$K_X = L^{-3}$.
Note that if we define $\WW:=\shF^{\oplus 2}\oplus\shE$, then 
\begin{equation}\label{eq:canonical}
\det\WW = \det(\shF^{\oplus 2}\oplus\shE) = L^3 = K_X^{-1}.
\end{equation}
\begin{lemma}\label{lem:ambient-sections}
On $X = N(3;3,4)$, the spaces of global sections of the universal bundles and their homomorphism bundles admit the following $GL(W)$-equivariant identifications:
\[
H^0(X,\shE) \simeq \Sym^3 W^*\ (\dim = 10),\qquad H^0(X,\shF) \simeq \mathbb{S}_{(3,1)}W^*\ (\dim = 15),
\]
\[
H^0(X,\End(\shE)) \simeq H^0(X,\End(\shF)) \simeq \bc,\qquad H^0(X,\shE^*\otimes\shF) \simeq W^*\ (\dim = 3),
\]
and $H^0(X,\shF^*\otimes\shE) = 0$. In particular, $H^0(X,K_1(\shE)) \simeq \bc$ and $H^0(X,K_1(\shF)) \simeq \bc^5$, where $K_1(\mathcal{T}) := \mathcal{T}\otimes\WW^*$.
\end{lemma}

\begin{proof}
Because $\mathrm{codim}_R(R - R^{\mathrm{st}}) \ge 2$, Hartogs's theorem identifies $H^0(X,\mathcal{T})$ with the space of $G$-invariants $(\bc[R]\otimes V)^G$ for any representation $V$ of $G$. Consider the one-parameter subgroup $\gamma(t) = (\mathrm{diag}(t,t,t),\mathrm{Id}_4)$ of $GL(V_3)\times GL(V_4)$. On the representations $V_4^*\otimes\chi$ and $V_3^*\otimes\chi$ that descend to $\shE$ and $\shF$, $\gamma(t)$ acts with weights $-3$ and $-4$ respectively. Meanwhile for the coordinate ring $\bc[R] = \bigoplus_{d\ge 0}\Sym^d(V_3\otimes V_4^*\otimes W^*)$, $\gamma$ acts on its degree-$d$ summand with weight $d$.

Consequently, the space of global sections of $\shF$ must lie in the degree-$4$ polynomial component $\Sym^4(V_3\otimes V_4^*\otimes W^*)\otimes V_3^*\otimes\det(V_4)\otimes\det(V_3)^{-1}$ by considering the action of $\gamma$. By representation theory, taking $GL(V_4)$-invariants must single out the direct summand $$[\Lambda^4(V_3\otimes W^*)\otimes \Lambda^4V_4^*]\otimes V_3^*\otimes \det V_4\otimes \det V_3^*=\Lambda^4(V_3\otimes W^*)\otimes V_3^*\otimes \det V_3^*.$$
Furthermore, since $V_3^*\otimes \det V_3^*=\mathbb S_{(2,1,1)}V_3^*$ and the transpose of the Young diagram $(2,1,1)$ is $(3,1)$, taking $GL(V_3)$-invariants yields $H^0(X,\shF) \simeq \mathbb{S}_{(3,1)}W^*$, which has dimension $15$.

By the same reasoning, taking $\gamma$-invariants and $GL(V_3)$-invariants first shows that $H^0(X,\shE)$ comes from $\Lambda^3(V_4^*\otimes W^*)\otimes V_4^*\otimes\det(V_4)$. Then taking $GL(V_4)$-invariants shows that $ H^0(X,\shE)\simeq \Sym^3 W^*$. 

Geometrically, at a representation $M : V_3 \to V_4\otimes W$, which is viewed as a $4\times 3$ matrix with entries in $W$, the map $\Sym^3 W^* \to \shE_{[M]}$ evaluates the four signed $3\times 3$ maximal minors of $M$, while $\mathbb{S}_{(3,1)}W^* \to \shF_{[M]}$ evaluates the linear syzygies among the maximal minors arising from Laplace expansions (recall that $\mathbb{S}_{(3,1)}W^*=\ker(S^3W^*\otimes W^*\rightarrow S^4W^*)$).

Similarly, on $\End(\shE)$, $\End(\shF)$, $\shE^*\otimes\shF$, and $\shF^*\otimes\shE$, the $\gamma$-weights on the fibers are $0, 0, -1$, and $+1$, respectively. Since $\bc[R]$ has only nonnegative $\gamma$-weights, $H^0(X,\shF^*\otimes\shE)=0$, $H^0(X,\End(\shE)) = (V_4^*\otimes V_4)^G = \bc$, $H^0(X,\End(\shF)) = (V_3^*\otimes V_3)^G = \bc$, and $H^0(X,\shE^*\otimes\shF) = (V_3^*\otimes V_4\otimes W^*\otimes V_4^*\otimes V_3)^G = W^*$. 
\end{proof}

\begin{proposition}\label{cor:gg-on-S}
On the twelvefold $X = N(3;3,4)$, the bundles $\shF$ and $\shE$ satisfy the following global generation properties:
\begin{enumerate}[(i)]
\item The rank-three bundle $\shF$ is globally generated everywhere on $X$.
\item The rank-four bundle $\shE$ is globally generated outside a closed subvariety $Z_{\mathrm{deg}} \subset X$ of dimension $5$. Along $Z_{\mathrm{deg}}$, the evaluation map $\ev_{\shE} : H^0(X,\shE)\otimes\sho_X \to \shE$ has rank $3$.
\item For a general global section $\sigma \in H^0(X,\WW)$ of $\WW = \shF^{\oplus 2}\oplus\shE$, the zero locus $S = Z(\sigma)$ is a smooth surface disjoint from $Z_{\mathrm{deg}}$, and the restricted bundles $E = \shE|_S$ and $F = \shF|_S$ are globally generated on $S$.
\end{enumerate}
\end{proposition}

\begin{proof}
Let $[M] \in X$ be a stable representation. Choosing bases of $V_3 \simeq \bc^3$, $V_4 \simeq \bc^4$ and an isomorphism $W=\bc\langle x,y,z\rangle$, we can view $M$ as a $4\times 3$ matrix with entries in $W$. Its maximal minors can be viewed as elements in $\Sym^3W$. Recall from the proof of Lemma~\ref{lem:ambient-sections} that the evaluation map at the point $[M]$ is as follows:
$$
\ev_{\shE,[M]} : H^0(X,\shE) \simeq \Sym^3 W^* \longrightarrow \shE_{[M]}.
$$
This evaluation map is not surjective if and only if the four maximal minors of the $4\times 3$ matrix are not linearly independent in $\Sym^3W$.

Suppose $[M] \in X$ lies in the degeneracy locus $Z_{\mathrm{deg}} \subset X$ where $\rk(\ev_{\shE,[M]}) \le 3$, so that the maximal minors are not linearly independent. By choosing suitable bases for $V_3$ and $V_4$, we can suppose that the determinant of the top $3\times 3$ submatrice $\overline{M}$ is zero. By the classification of $3\times 3$ matrices of linear forms with vanishing determinant in \cite[Theorem 1.1]{EisenbudHarris}, $\overline{M}$ is either:
\begin{enumerate}[(a)]
\item \emph{compressible}, meaning that for some $a \in \{1,2,3\}$ there is an $a$-dimensional subspace $C_a$ mapped by $\overline{M}(x,y,z)$ into an $(a-1)$-dimensional subspace for all $(x,y,z)$; or
\item in the $GL(\bc^3)\times GL(\bc^3)$-orbit of the $3\times 3$ skew-symmetric matrix representing $\bc^3=\Lambda^2 W \to W\otimes W=\bc^3\otimes W$.
\end{enumerate}
Case~(a) is impossible for a stable representation $[M] \in X$ because otherwise $M$ will admit a subrepresentation of dimension vector $(a,a)$, violating the stability inequality $4a - 3a = a > 0$.

Thus only case~(b) can occur and we can choose suitable bases for $V_3$ and $V_4$ such that $\overline{M}(x,y,z)$ is the standard skew-symmetric matrix:
\begin{equation}\label{b-form-matrix}
    M(x,y,z) = \begin{pmatrix}
0 & z & -y \\
-z & 0 & x \\
y & -x & 0 \\
\ell_1 & \ell_2 & \ell_3
\end{pmatrix},\qquad \ell_i \in W.
\end{equation}

A direct calculation shows that the four $3\times 3$ maximal minors of $M(x,y,z)$ are
\[
\bigl(x Q,\, y Q,\, z Q,\, 0\bigr),
\]
where $Q := \ell_1x + \ell_2y + \ell_3z \in \Sym^2 W$. Notice that $Q$ cannot vanish identically: if $\ell_1 x + \ell_2 y + \ell_3 z \equiv 0$, exactness of the Koszul complex associated to the regular sequence $(x,y,z)$ implies that the fourth row $(\ell_1,\ell_2,\ell_3)$ is a $\bc$-linear combination of the first three rows of $M(x,y,z)$, which means $\image(M) \subset V_4$ is contained in a $3$-dimensional subspace $D_3 \subset V_4$, giving a destabilizing subrepresentation of dimension vector $(3,3)$. Since $Q \neq 0$ in $\Sym^2 W$, the three cubic polynomials $xQ, yQ, zQ \in \Sym^3 W$ are linearly independent. Consequently, $\rk(\ev_{\shE,[M]}) = 3$ at every point $[M]\in Z_{\mathrm{deg}}$.

The above quadratic form $Q$ is determined up to scalar by the $GL(V_3)\times GL(V_4)$-orbit of $M$ (with $[M]\in Z_{\mathrm{deg}}$): it is the unique quadratic element which divides all maximal minors of the $4\times 3$ matrix $M$. Conversely, such $Q$ actually determines the orbit of $M$: if $M,M'$ are two elements in the form of (\ref{b-form-matrix}) with the same $Q$, then $(\ell_1-\ell_1')x+(\ell_2-\ell_2')y+(\ell_3-\ell_3')z=0$. By the Koszul complex, the row vector $(\ell_1-\ell_1', \ \ell_2-\ell_2', \ \ell_3-\ell_3')$ must be a linear combination of the first three rows in (\ref{b-form-matrix}). Consequently, $M$ can be obtained by applying elementary row operations to $M'$, which shows that $M$ and $M'$ lie in the same $GL(V_3)\times GL(V_4)$-orbit.
As a result, we see that $\dim Z_{\mathrm{deg}}=\dim(\bp(\Sym^2W))=5$.

To show that $\mathcal F$ is globally generated, we only need to show that the composition 
$$S^3W^*\otimes W^*=H^0(X,\shE)\otimes W^*\rightarrow \mathcal{E}\otimes W^*\rightarrow\mathcal F$$ is surjective.
By the stability of quiver representations, the universal map $\mathcal{E}\otimes W^*\rightarrow\mathcal F$ is a surjective bundle morphism. This fact together with the above analysis for $\mathcal E$ shows that the bundle $\mathcal F$ is globally generated at all points $[M] \notin Z_{\mathrm{deg}}$. 
 For $[M] \in Z_{\mathrm{deg}}$, we only need to show that the transpose of the above composition map $\mathcal{F}^*\rightarrow \mathcal E^*\otimes W\rightarrow S^3W\otimes W$ is injective, which can be proved using the explicit form (\ref{b-form-matrix}).

By~(i) and~(ii), the bundle $\WW = \shF^{\oplus 2}\oplus\shE$ is globally generated on $X - Z_{\mathrm{deg}}$, while along $Z_{\mathrm{deg}}$ the rank of the evaluation map $H^0(X,\WW)\otimes\sho_X \to \WW$ is equal to $2\cdot 3 + 3 = 9$. Consequently, for a general section $\sigma \in H^0(X,\WW)$, its zero locus doesn't intersect $Z_{\mathrm{deg}}$ (since $\dim Z_{\mathrm{deg}}=5$) and is a smooth surface $S$.
The bundles $E = \shE|_S$ and $F = \shF|_S$ are globally generated on $S$ by~(i) and~(ii).
\end{proof}

\subsection{\texorpdfstring{Cohomology on the auxiliary sixfold $Y_{2,3} = N(3;2,3)$}{Cohomology on the auxiliary sixfold Y(2,3) = N(3;2,3)}}\label{subsec:auxiliary-Y}
For future use, we consider
$
Y_{2,3} = N(3;2,3),
$
which is the six-dimensional fine moduli space of stable representations of the $3$-Kronecker quiver with dimension vector $(2,3)$ and stability parameter $(3,-2)$. This is actually an interesting variety which is equal to a blow-down of $\mathrm{Hilb}^3(\bp^2)$, with many other properties studied in \cite{IlievManivel,meng2024quantumcohomologyfanoquiver,Derquiver}. We collect here some cohomological information of vector bundles on $Y_{2,3}$ required for the calculations in Section~\ref{sec:unstable-strata}.

The sixfold $Y_{2,3}$ is smooth and Fano of Picard number $1$, with ample generator $\sho_{Y_{2,3}}(1)$ and canonical bundle $K_{Y_{2,3}} = \sho_{Y_{2,3}}(-3)$. Choosing the normalization vector $\underline{a}=(-1,1)$ yields universal bundles $\shU$ (of rank $2$) and $\shV$ (of rank $3$) on $Y_{2,3}$ such that $\shU\otimes W^* \to \shV$ is the universal representation and
\[
\det\shU = \det\shV = \sho_{Y_{2,3}}(-1).
\]
In particular, we have the canonical identifications
\begin{equation}\label{eq:small-rank-identities}
\shU^* = \shU(1),\qquad \Lambda^2\shV = \shV^*(-1),\qquad \Sym^2\shU = \slf(\shU)(-1),
\end{equation}
where $\slf(\shU)$ denotes the trace-free endomorphism bundle of the rank-two bundle $\shU$.

By \cite[Theorem~1.1]{Derquiver}, the bounded derived category $D^b(Y_{2,3})$ admits the full strong exceptional collection
\begin{equation}\label{eq:exceptional-collection}
\begin{aligned}
D^b(Y_{2,3}) = \bigl\langle
&\sho_{Y_{2,3}},\,\shV^*,\,\shU^*,\,\shV(1),\,\slf(\shU)(1),\ \,
\sho_{Y_{2,3}}(1),\,\shV^*(1),\,\shU^*(1),\,\shV(2),\\
&\sho_{Y_{2,3}}(2),\,\shV^*(2),\,\shU^*(2),\,\shV(3)
\bigr\rangle.
\end{aligned}
\end{equation}
Setting $\smallD := \shU^{\oplus 2}\oplus\shV$, we have the following computations of cohomology.

\begin{lemma}\label{lem:small-cohom}
On $Y_{2,3} = N(3;2,3)$, the vector bundles $\shU$, $\Lambda^2\shV$, $\shU\otimes\smallD$, $\Lambda^2\shV\otimes\smallD$, $\sho_{Y_{2,3}}(-1)$, and $\smallD(-1)$ are acyclic ($\RG(Y_{2,3},-) = 0$). Furthermore, if we set
\[
\mathcal{J}_3 = \shU\otimes\Lambda^2\smallD,\qquad \mathcal{J}_4 = \Lambda^2\shV\otimes\Lambda^2\smallD,
\]
then
\begin{align}
\RG(Y_{2,3},\mathcal{J}_3) &= 0,\label{eq:Kprime}\\
\RG(Y_{2,3},\mathcal{J}_4) &= \bc^2[-3].\label{eq:K}
\end{align}
\end{lemma}

\begin{proof}
In each row of Table~\ref{tab:semiorthogonal-Y}, the bundle in the first column is isomorphic to $\mathcal{A}^*\otimes\mathcal{B}$, where $\mathcal{A}$ (second column) appears to the right of $\mathcal{B}$ (third column) in the exceptional collection \eqref{eq:exceptional-collection}. By semiorthogonality, $\RG(Y_{2,3},\mathcal{A}^*\otimes\mathcal{B}) = \RHom(\mathcal{A},\mathcal{B}) = 0$, so every bundle in the first column is acyclic.

\begin{table}[htbp]
\centering
\small
\begin{tabular}{lll}
\toprule
\textbf{Bundle on $Y_{2,3}$} & \textbf{Later object $\mathcal{A}$} & \textbf{Earlier object $\mathcal{B}$} \\
\midrule
$\shU$ & $\shU^*$ & $\sho_{Y_{2,3}}$ \\
$\Lambda^2\shV \simeq \shV^*(-1)$ & $\sho_{Y_{2,3}}(1)$ & $\shV^*$ \\
$\sho_{Y_{2,3}}(-1)$ & $\sho_{Y_{2,3}}(1)$ & $\sho_{Y_{2,3}}$ \\
$\Sym^2\shU \simeq \slf(\shU)(-1)$ & $\slf(\shU)(1)$ & $\sho_{Y_{2,3}}$ \\
$\shU\otimes\shV$ & $\shU^*(1)$ & $\shV(1)$ \\
$\shU\otimes\shV^*(-1)$ & $\shU^*(1)$ & $\shV^*$ \\
$\End(\shV)(-1)$ & $\shV^*(1)$ & $\shV^*$ \\
$\shU(-1)$ & $\shU^*(1)$ & $\sho_{Y_{2,3}}$ \\
$\shV(-1)$ & $\shV^*(1)$ & $\sho_{Y_{2,3}}$ \\
$\slf(\shU)\otimes\shV(-1)$ & $\shV^*(2)$ & $\slf(\shU)(1)$ \\
\bottomrule
\end{tabular}
\caption{Acyclicity of elementary tensor bundles on $Y_{2,3}$ via semiorthogonality of \eqref{eq:exceptional-collection}.}\label{tab:semiorthogonal-Y}
\end{table}

Since $\shU\otimes\shU \simeq \slf(\shU)(-1)\oplus\sho_{Y_{2,3}}(-1)$ and $\Lambda^2\shV\otimes\smallD \simeq (\shU\otimes\shV^*(-1))^{\oplus 2}\oplus\End(\shV)(-1)$, the first set of acyclicity assertions follows immediately from Table~\ref{tab:semiorthogonal-Y}.

To compute $\mathcal{J}_3$, we expand the exterior square of $\smallD = \shU^{\oplus 2}\oplus\shV$:
\begin{equation}\label{eq:small-exterior-square}
\Lambda^2\smallD \simeq \sho_{Y_{2,3}}(-1)^{\oplus 3}\oplus\Sym^2\shU\oplus[(\shU\otimes\shV)^{\oplus 2}]\oplus\Lambda^2\shV.
\end{equation}
Using $\shU\otimes\Sym^2\shU \simeq \Sym^3\shU\oplus\shU(-1)$ and $\Sym^2\shU \simeq \slf(\shU)(-1)$, we obtain
\begin{equation*}\label{eq:J3-decomposition}
\mathcal{J}_3 \simeq \Sym^3\shU \oplus \shU(-1)^{\oplus 4} \oplus [\bigl(\slf(\shU)\otimes\shV(-1)\bigr)^{\oplus 2} \oplus \shV(-1)^{\oplus 2}] \oplus \bigl(\shU\otimes\shV^*(-1)\bigr).
\end{equation*}
Every direct summand in the above formula except $\Sym^3\shU$ appears in Table~\ref{tab:semiorthogonal-Y}. To verify that $\Sym^3\shU$ is also acyclic, we apply Lemma~\ref{Teleman} on $Y_{2,3} = N(3;2,3)$ with stability parameter $(3,-2)$. The seven unstable HN strata of $(2,3)$ have the following codimensions and weights:
\begin{center}
\small
\begin{tabular}{lrrr}
\toprule
\textbf{HN type on $Y_{2,3}$} & $c_\lambda$ & $\eta_\lambda$ & $m_\lambda(\Sym^3\shU)$ \\
\midrule
$((1,0),(1,1),(0,2))$ & $12$ & $90$ & $60$ \\
$((1,0),(1,2),(0,1))$ & $9$ & $100$ & $75$ \\
$((1,0),(1,3))$ & $8$ & $120$ & $90$ \\
$((1,1),(1,2))$ & $3$ & $15$ & $15$ \\
$((2,0),(0,3))$ & $18$ & $90$ & $45$ \\
$((2,1),(0,2))$ & $10$ & $100$ & $60$ \\
$((2,2),(0,1))$ & $4$ & $20$ & $15$ \\
\bottomrule
\end{tabular}
\end{center}
The inequality $m_\tau(\Sym^3\shU) \le \eta_\tau$ is strict for all HN types except $\tau = ((1,1),(1,2))$. The slopes of the subdimension vectors are $\mu(1,1)=1/2, \ \mu(1,2)=-1/3$ and the corresponding primitive 1-ps is $\lambda(t) = (\mathrm{diag}(t^3,t^{-2}),\mathrm{diag}(t^3,t^{-2},t^{-2}))$. By Remark \ref{gerbe}, we only need to calculate the cohomology of the critical weight subbundle over the quotient stack
$$[Z_\lambda^\circ/G_{(1,1)}\times G_{(1,2)}]\simeq N(3;1,1)\times N(3;1,2) \simeq \bp(W) \times \bp(W^*).$$
We write $(a,b)$ for the source and target lines of the $(1,1)$ block and $(c,D)$ for the source line and target plane of the $(1,2)$ block. Setting $h_1 = a\otimes b^{-1}$ and $\delta_1 = c\otimes(\det D)^{-1}$, the bundle $\shU$ is obtained by descending $(a\oplus c)\otimes h_1\otimes\delta_1$, and the highest-weight piece of $\Sym^3\shU$ is $(a\otimes h_1\otimes\delta_1)^{\otimes 3}$ (with weight $15$). By \eqref{eq:normal-two-block-det}, $\det\mathcal{N}_\lambda = a^{-5}\otimes b^2\otimes c^{-1}\otimes(\det D)^2$, whence
\[
\mathcal{C}_\lambda(\Sym^3\shU) = (a\otimes b^{-1})\boxtimes \bigl(c^2\otimes(\det D)^{-1}\bigr) \simeq \sho_{\bp(W)}(-1)\boxtimes\sho_{\bp(W^*)}(-1),
\]
which is acyclic on $\bp(W) \times \bp(W^*)$ (notice that for the second factor, $\sho_{\bp(W^*)}(1)$ is the determinant of the universal quotient bundle, thus equal to $\det(D\otimes c^{-1})=c^{-2}\otimes \det D$). Thus $H^{>0}(Y_{2,3},\Sym^3\shU) = 0$. Moreover, $H^0(Y_{2,3},\Sym^3\shU) = 0$ because $z\cdot\mathrm{Id}_{\bc^2}$ acts on degree-$n$ polynomials in $\bc[R_{(2,3)}]$ by $z^{n}$ and on the fiber of $\Sym^3\shU$ by $z^9$. This proves \eqref{eq:Kprime}.

Finally, tensoring \eqref{eq:small-exterior-square} with $\Lambda^2\shV \simeq \shV^*(-1)$ gives
\begin{equation*}\label{eq:K-decomposition}
\mathcal{J}_4 \simeq \shV^*(-2)^{\oplus 3} \oplus \bigl(\slf(\shU)\otimes\shV^*(-2)\bigr) \oplus \bigl(\shU\otimes\End(\shV)(-1)\bigr)^{\oplus 2} \oplus \bigl(\shV^*\otimes\shV^*(-2)\bigr).
\end{equation*}
Its first, second and fourth direct summands are equal to $\HHom(\sho_{Y_{2,3}}(2),\shV^*)^{\oplus 3}$, 
$\HHom(\shV(3),\slf(\shU)(1))$ and $\HHom(\shV(2),\shV^*)$ respectively, which are acyclic by the semiorthogonality relation in \eqref{eq:exceptional-collection}. For the remaining summand, \cite[Lemma~4.31]{Derquiver} establishes $\RG(Y_{2,3},\shU\otimes\End(\shV)(-1)) = \bc[-3]$, which proves \eqref{eq:K}.
\end{proof}

\section{\texorpdfstring{Unstable strata and weight analysis for $X = N(3;3,4)$}{Unstable strata and boundary weight analysis on X = N(3;3,4)}}\label{sec:unstable-strata}

\subsection{\texorpdfstring{Classification of unstable Harder--Narasimhan strata}{Classification of unstable Harder-Narasimhan strata on X}}

To compute the cohomology of the restricted bundles such as $E = \shE|_S,F = \shF|_S$ on a smooth 2-dimensional zero locus $S = Z(\sigma)$ via the Koszul complex of $\WW = \shF^{\oplus 2}\oplus\shE$, we must control the cohomology of the twisted exterior powers
\[
K_p(\mathcal{T}) := \mathcal{T}\otimes\Lambda^p\WW^*,\qquad 0 \le p \le 10,
\]
where $\WW^* = (\shF^*)^{\oplus 2}\oplus\shE^*$ and $\mathcal{T}$ ranges over $\sho_X$, $\shE$, $\shF$, $\shE^*$, $\shF^*$, $\End(\shE)$, $\End(\shF)$, $\shE^*\otimes\shF$, and $\shF^*\otimes\shE$. By Serre duality on the twelvefold $X$ (see \eqref{eq:serre-koszul} below), it suffices to analyze $K_p(\mathcal{T})$ in the range $0 \le p \le 5$ (and $0 \le p \le 4$ for $\shE^*,\shF^*$).

Table~\ref{tab:all-18-hn} lists all eighteen unstable Harder--Narasimhan types $\tau$ for the $3$-Kronecker quiver with dimension vector $(3,4)$ and stability parameter $\theta=(4,-3)$, together with their codimensions $c_\lambda = \mathrm{codim}_R(S_\lambda)$, the weights $\eta_\lambda = \wt_\lambda(\det\mathcal{N}_\lambda^*)$, the ten $\lambda$-weights of $\WW^*$ in nonincreasing order and the maximum weights $m_\lambda(K_p(\sho_X))$ for $p \le 5$. In the last column we record the maximum weights $m_\lambda(K_p(\mathcal{T}_{\Hom}))$ for $p \le 5$ when $\mathcal{T}_{\Hom}$ ranges over the homomorphism bundles $\End(\shE), \End(\shF), \shE^*\otimes\shF, \shF^*\otimes\shE$. In Appendix \ref{app:hn}, we include more basic information on the weights for unstable HN types for the reader's convenience.

\begin{table}[htbp]
\centering
\renewcommand{\arraystretch}{1.15}
\resizebox{\textwidth}{!}{%
\begin{tabular}{lrrlrr}
\toprule
\textbf{Harder--Narasimhan Type $\tau$} & $c_\lambda$ & $\eta_\lambda$ & \textbf{$\lambda$-Weights of $\WW^*$ (10-tuple)} & $\max_{p\le 5} m_\lambda(K_p(\sho_X))$ & $\max_{p\le 5} m_\lambda(K_p(\mathcal{T}_{\Hom}))$ \\
\midrule
$((1,0),(1,1),(1,2),(0,1))$ & $16$ & $385$ & $(70^2, 49^3, 42^4, 28)$ & $287$ & $329$ \\
$((1,0),(1,1),(1,3))$ & $15$ & $231$ & $(42^2, 28^3, 21^5)$ & $168$ & $189$ \\
$((1,0),(2,1),(0,3))$ & $25$ & $406$ & $(56^2, 49^5, 35^3)$ & $259$ & $280$ \\
$((1,0),(2,2),(0,2))$ & $18$ & $168$ & $(28^2, 21^6, 14^2)$ & $119$ & $133$ \\
$((1,0),(2,3),(0,1))$ & $13$ & $294$ & $(56^2, 35^7, 21)$ & $217$ & $252$ \\
$((1,0),(2,4))$ & $10$ & $70$ & $(14^2, 7^8)$ & $49$ & $56$ \\
$\mathbf{((1,1),(2,3))}$ & $\mathbf{4}$ & $\mathbf{28}$ & $\mathbf{(7^3, 0^7)}$ & $\mathbf{21}$ & $\mathbf{28}$ \\
$((2,0),(1,1),(0,3))$ & $28$ & $322$ & $(42^4, 35^3, 28^3)$ & $203$ & $217$ \\
$((2,0),(1,2),(0,2))$ & $24$ & $406$ & $(56^4, 42^4, 35^2)$ & $266$ & $287$ \\
$((2,0),(1,3),(0,1))$ & $22$ & $504$ & $(70^4, 49^5, 42)$ & $329$ & $357$ \\
$((2,1),(1,1),(0,2))$ & $17$ & $385$ & $(56^5, 49^3, 28^2)$ & $280$ & $308$ \\
$((2,1),(1,2),(0,1))$ & $14$ & $133$ & $(21^5, 14^4, 7)$ & $105$ & $119$ \\
$((2,1),(1,3))$ & $13$ & $91$ & $(14^5, 7^5)$ & $70$ & $77$ \\
$\mathbf{((2,2),(1,2))}$ & $\mathbf{6}$ & $\mathbf{42}$ & $\mathbf{(7^6, 0^4)}$ & $\mathbf{35}$ & $\mathbf{42}$ \\
$((3,0),(0,4))$ & $36$ & $252$ & $(28^6, 21^4)$ & $140$ & $147$ \\
$((3,1),(0,3))$ & $24$ & $168$ & $(21^7, 14^3)$ & $105$ & $112$ \\
$((3,2),(0,2))$ & $14$ & $98$ & $(14^8, 7^2)$ & $70$ & $77$ \\
$\mathbf{((3,3),(0,1))}$ & $\mathbf{6}$ & $\mathbf{42}$ & $\mathbf{(7^9, 0^1)}$ & $\mathbf{35}$ & $\mathbf{42}$ \\
\bottomrule
\end{tabular}%
}
\caption{Complete Harder--Narasimhan weight table for $X = N(3;3,4)$ in the Koszul range $p \le 5$. Equality $m_\lambda = \eta_\lambda$ occurs only on the three boldfaced strata.}\label{tab:all-18-hn}
\end{table}

Inspecting Table~\ref{tab:all-18-hn} immediately yields two key consequences:
\begin{enumerate}[(1)]
\item For $\mathcal{T} = \sho_X, \shE, \shF$ and every $0 \le p \le 5$, the strict inequality $m_\tau(K_p(\mathcal{T})) < \eta_\tau$ holds for all the eighteen unstable Harder--Narasimhan types. (Note that the weights of $K_p(\shE)$ and $K_p(\shF)$ are always smaller than or equal to those of $K_p(\sho_X)$.)
\item For $\mathcal{T} \in \{\End(\shE), \End(\shF), \shE^*\otimes\shF, \shF^*\otimes\shE\}$ and $1 \le p \le 5$ (as well as $\mathcal{T} \in \{\shE^*,\shF^*\}$ for $1 \le p \le 4$), the inequality $m_\tau(K_p(\mathcal{T})) \le \eta_\tau$ holds for all eighteen unstable HN types, and the equality $m_\tau(K_p(\mathcal{T})) = \eta_\tau$ can occur only for the three two-block HN types below:
\begin{equation*}\label{eq:three-types}
\begin{array}{c|cccc}
\toprule
\text{Critical HN type } \tau & c_\tau & \eta_\tau & \text{Weights of } \WW^* & \text{Koszul degrees } p \le 5 \text{ with } m_\lambda(K_p(\mathcal T)) = \eta_\lambda \\
\midrule
((1,1),(2,3)) & 4 & 28 & (7^3, 0^7) & p = 3, 4, 5 \\
((2,2),(1,2)) & 6 & 42 & (7^6, 0^4) & p = 5 \\
((3,3),(0,1)) & 6 & 42 & (7^9, 0^1) & p = 5 \\
\bottomrule
\end{array}
\end{equation*}
\end{enumerate}
We now analyze each of these three unstable Hesselink-Kempf-Ness strata in turn.

\subsection{\texorpdfstring{The unstable HN type $\tau = ((1,1),(2,3))$}{The boundary stratum tau = ((1,1),(2,3))}}\label{subsec:boundary-11-23}

For $\tau = ((1,1),(2,3))$, the two semistable blocks in the HN filtration have dimension vectors $(1,1)$ and $(2,3)$ respectively.
By Remark \ref{gerbe}, we only need to calculate the cohomology of critical weight subbundles over the quotient stack
$$[Z_\lambda^\circ/G_{(1,1)}\times G_{(2,3)}]\simeq N(3;1,1)\times N(3;2,3) \simeq \bp(W) \times Y_{2,3},$$
where $Y_{2,3} = N(3;2,3)$ is the auxiliary sixfold studied in Section~\ref{subsec:auxiliary-Y}. Let $(a,b)$ denote the source and target lines of the $(1,1)$ block, and let $(C_2,D_3)$ denote the source and target spaces of the $(2,3)$ block, so that $V_3 = a\oplus C_2$ and $V_4 = b\oplus D_3$. Writing $\chi^{-1} = \det(V_3)\otimes\det(V_4)^{-1} = a\otimes b^{-1}\otimes\det(C_2)\otimes\det(D_3)^{-1}$, the restrictions of the universal bundles $\shF^* = \mathcal U_3$ and $\shE^* = \mathcal U_4$ to $Z_\lambda^\circ$ decompose as follows:
\[
\shF^*|_{Z_\lambda^\circ} = (a\otimes\chi^{-1})\oplus(C_2\otimes\chi^{-1}),\qquad
\shE^*|_{Z_\lambda^\circ} = (b\otimes\chi^{-1})\oplus(D_3\otimes\chi^{-1}).
\]
On the first factor $\bp^2$, the line bundle $\ell := a^{-1}\otimes b$ descends to $\sho_{\bp^2}(1)$, while on the second factor $Y_{2,3}$, the bundles $C_2\otimes\chi^{-1}\otimes\ell$ and $D_3\otimes\chi^{-1}\otimes\ell$ descend to the universal bundles $\shU$ and $\shV$ respectively. Note also that $\sho_{Y_{2,3}}(-1) = \det\shV$ can be obtained by descending $(\det C_2)^3\otimes(\det D_3)^{-2}$.

The slopes are $\mu(1,1) = 1/2$ and $\mu(2,3) = -1/5$, so the corresponding primitive one-parameter subgroup is
\[
\lambda(t) = \bigl(\mathrm{diag}(t^5,t^{-2},t^{-2}),\,\mathrm{diag}(t^5,t^{-2},t^{-2},t^{-2})\bigr).
\]
Under the action of $\lambda$, the spaces $a, b$, and $\chi^{-1}$ have weights $5, 5$ and $2$ respectively, while $C_2\otimes\chi^{-1}$ and $D_3\otimes\chi^{-1}$ have weight $0$. Thus $\WW^*|_{Z_\lambda^\circ}$ has three weight-$7$ directions, namely $(a\otimes\chi^{-1})^{\oplus 2}\oplus(b\otimes\chi^{-1})$, and seven weight-$0$ directions, namely $\smallD\otimes\ell^{-1}$, where $\smallD = \shU^{\oplus 2}\oplus\shV$. For $p = 3+k$ ($k \ge 0$), the highest-weight summand of $\Lambda^{3+k}\WW^*$ has weight $21$ and is given by taking all the three weight-$7$ directions and $k$ weight-$0$ directions:
\begin{equation}\label{eq:first-top-exterior}
\bigl(\Lambda^{3+k}\WW^*\bigr)_{21} = a^5\otimes b^{-2}\otimes\bigl(\det C_2\otimes(\det D_3)^{-1}\bigr)^{\otimes 3}\otimes \Lambda^k(\smallD\otimes\ell^{-1}).
\end{equation}
Meanwhile, formula \eqref{eq:normal-two-block-det} gives
\begin{equation}\label{eq:first-normal-det}
\det\mathcal{N}_\lambda = a^{-7}\otimes b^3\otimes(\det D_3)^2\otimes(\det C_2)^{-1},
\end{equation}
whose weight is $-\eta_\lambda = -28$, as listed before.

Combining \eqref{eq:first-top-exterior} and \eqref{eq:first-normal-det}, we can compute the critical weight subbundles $$\mathcal{C}_\lambda(K_{3+k}(\mathcal{T})) = 
(\widetilde{K_{3+k}(\mathcal{T})}|_{Z_\lambda^\circ})_{28}\otimes\det\mathcal{N}_\lambda$$ on $\bp^2\times Y_{2,3}$. For example, for $\mathcal{T} = \End(\shE)$, the highest-weight piece of $\End(\shE|_{Z_\lambda^\circ})$ is $D_3^*\otimes b$ (of weight $7$), so
\begin{align}\label{eq:first-endE-calc}
\mathcal{C}_\lambda\bigl(K_{3+k}(\End(\shE))\bigr)
&= (D_3^*\otimes b)\otimes \bigl(\Lambda^{3+k}\WW^*\bigr)_{21}\otimes\det\mathcal{N}_\lambda \nonumber\\
&= (D_3^*\otimes b)\otimes a^{-2}\otimes b\otimes(\det C_2)^2\otimes(\det D_3)^{-1}\otimes\Lambda^k(\smallD\otimes\ell^{-1}) \\
&\simeq \ell^{2-k}\boxtimes\bigl(\Lambda^2\shV\otimes\Lambda^k\smallD\bigr),\nonumber
\end{align}
where the last isomorphism uses $\Lambda^2\shV \simeq \shV^*(-1) \leftrightarrow D_3^*\otimes(\det C_2)^2\otimes(\det D_3)^{-1}$. Performing the analogous calculations for other bundles $\mathcal{T}$ gives the following table on $\bp^2\times Y_{2,3}$:
\begin{equation}\label{eq:boundary-table}
\begin{array}{c|l}
\toprule
\mathcal{T} & \text{Critical weight subbundles } \mathcal{C}_\lambda\bigl(K_{3+k}(\mathcal{T})\bigr)  \\
\midrule
\End(\shF) & \ell^{1-k}\boxtimes(\shU\otimes\Lambda^k\smallD) \\
\End(\shE) & \ell^{2-k}\boxtimes(\Lambda^2\shV\otimes\Lambda^k\smallD) \\
\shE^*\otimes\shF & \ell^{2-k}\boxtimes(\shU\otimes\Lambda^k\smallD) \\
\shF^*\otimes\shE & \ell^{1-k}\boxtimes(\Lambda^2\shV\otimes\Lambda^k\smallD) \\
\shF^* & \ell^{-k}\boxtimes(\Lambda^k\smallD)(-1) \\
\shE^* & \ell^{1-k}\boxtimes(\Lambda^k\smallD)(-1) \\
\bottomrule
\end{array}
\end{equation}
In the following, we only care about the cases $k \in \{0,1,2\}$ ($p=3,4,5$) in the first three rows and $k \in \{0,1\}$ ($p=3,4$) in the last three rows.

By Lemma~\ref{lem:small-cohom} and the K\"unneth formula on $\bp^2\times Y_{2,3}$:
\begin{itemize}
\item For $k=0$ and $k=1$, every bundle in \eqref{eq:boundary-table} is acyclic because its $Y_{2,3}$-factor is acyclic by Lemma~\ref{lem:small-cohom}.
\item For $k=2$ ($p=5$), the critical weight subbundles for $\End(\shF)$ and $\shE^*\otimes\shF$ are $\ell^{-1}\boxtimes\mathcal{J}_3$ and $\sho_{\bp^2}\boxtimes\mathcal{J}_3$ respectively, both of which are acyclic by \eqref{eq:Kprime}.
\item For $k=2$ ($p=5$) and $\mathcal{T} = \End(\shE)$, the critical weight subbundle is $\sho_{\bp^2}\boxtimes\mathcal{J}_4$. By \eqref{eq:K}, its cohomology is concentrated entirely in degree $3$:
\[
H^j\bigl(\bp^2\times Y_{2,3},\,\sho_{\bp^2}\boxtimes\mathcal{J}_4\bigr) = \begin{cases} \bc^2, & j = 3, \\ 0, & j \neq 3. \end{cases}
\]
By equation \eqref{eq:boundary-degree-target} and the fact that $c_\lambda = 4$, to exclude the contribution to $H^5(X,K_5(\End(\shE)))$  from the stratum $S_\tau$, where $\tau=((1,1),(2,3))$, the vanishing $H^2(\bp^2\times Y_{2,3},\sho_{\bp^2}\boxtimes\mathcal{J}_4) = 0$ is exactly what we want.
\end{itemize}

\subsection{\texorpdfstring{The unstable HN type $\tau = ((2,2),(1,2))$}{The boundary stratum tau = ((2,2),(1,2))}}\label{subsec:boundary-22-12}

For $\tau = ((2,2),(1,2))$, the slopes are $\mu(2,2) = 1/2$ and $\mu(1,2) = -2/3$, and the corresponding primitive one-parameter subgroup is
\[
\lambda(t) = \bigl(\mathrm{diag}(t^3,t^3,t^{-4}),\,\mathrm{diag}(t^3,t^3,t^{-4},t^{-4})\bigr).
\]
The first subquotient in the HN filtration has dimension vector $(2,2)$ with source and target planes denoted by $(C_2,D_2)$; the corresponding quotient stack is
\[
\mathfrak{Z}_{2,2} = \left[\bigl(\Hom(C_2,D_2)\otimes W\bigr)^{\mathrm{ss}} \Big/ \bigl((GL(C_2)\times GL(D_2))/\bc^*\bigr)\right].
\]
The second subquotient in the HN filtration has dimension vector $(1,2)$ with source line $c$ and target plane $D$, and its fine moduli space is $\bp^2$. Thus by Remark \ref{gerbe}, we only need to calculate the cohomology of critical weight subbundles over the quotient stack
$$[Z_\lambda^\circ/G_{(2,2)}\times G_{(1,2)}]\simeq \mathfrak{Z}_{2,2}\times \bp(W^*).$$

Setting $\delta = \det(D_2)\otimes\det(C_2)^{-1}$ and $\rho = \delta^{-1}\otimes c\otimes(\det D)^{-1}$, the restrictions of the universal bundles decompose as follows:
\[
\shF^*|_{Z_\lambda^\circ} = \rho\otimes(C_2\oplus c),\qquad \shE^*|_{Z_\lambda^\circ} = \rho\otimes(D_2\oplus D).
\]
On the $\bp^2$ factor, let $\mathcal{Q}_2$ denote the universal rank-two quotient bundle fitting into $0 \to \sho_{\bp^2}(-1) \to \sho_{\bp^2}^{\oplus 3} \to \mathcal{Q}_2 \to 0$. It can be obtained by descending $D\otimes c^{-1}$. Then
\begin{equation}\label{eq:second-P2-identities}
\sho_{\bp^2}(1) = \det\mathcal{Q}_2 = c^{-2}\otimes\det D,\qquad \mathcal{Q}_2(-2) = c^3\otimes(\det D)^{-1}\otimes D^*.
\end{equation}

Here $\WW^*|_{Z_\lambda^\circ}$ has six weight-$7$ directions $\rho\otimes(C_2^{\oplus 2}\oplus D_2)$ and four weight-$0$ directions. In the range $p \le 5$, the equality $m_\lambda(K_p(\mathcal{T})) = \eta_\lambda = 42$ occurs only at $p=5$, where $\Lambda^5\WW^*$ contributes weight $35$ and $\mathcal{T} \in \{\End(\shE),\End(\shF),\shE^*\otimes\shF\}$ contributes weight $7$:
\begin{equation}\label{eq:second-top-exterior}
\begin{split}
 \bigl(\Lambda^5\WW^*\bigr)_{35}
 &=\rho^5\otimes\Lambda^5(C_2^{\oplus2}\oplus D_2)=\rho^5\otimes\det(C_2^{\oplus2}\oplus D_2)\otimes [C_2^{*\oplus2}\oplus D_2^*]\\
 &=\rho^5\otimes(\det C_2)^2\otimes \det(D_2)
              \otimes[C_2^{*\oplus2}\oplus D_2^*\bigr].
\end{split}
\end{equation}
Formula \eqref{eq:normal-two-block-det} gives
\begin{equation}\label{eq:second-normal-det}
\det\mathcal{N}_\lambda = (\det C_2)^{-5}\otimes(\det D_2)^2\otimes c^{-2}\otimes(\det D)^4,
\end{equation}
and the product of the two formulas above is $\delta^{-2}\otimes c^3\otimes(\det D)^{-1}\otimes((C_2^*)^{\oplus 2}\oplus D_2^*)$. Since the weight-$7$ direct summands of $\End(\shE)$, $\End(\shF)$, and $\shE^*\otimes\shF$ are $D_2\otimes D^*$, $C_2\otimes c^{-1}$ and $D_2\otimes c^{-1}$ respectively, we obtain the following table using \eqref{eq:second-P2-identities}:
\begin{equation}\label{eq:second-boundary}
\begin{array}{c|l}
\toprule
\mathcal{T} & \text{Critical weight subbundle } \mathcal{C}_\lambda\bigl(K_5(\mathcal{T})\bigr) \text{ on } \mathfrak{Z}_{2,2}\times\bp^2 \\
\midrule
\End(\shE) & \delta^{-2}\otimes\bigl((D_2\otimes C_2^*)^{\oplus 2}\oplus\End( D_2)\bigr)\boxtimes\mathcal{Q}_2(-2) \\
\End(\shF) & \delta^{-2}\otimes C_2\otimes\bigl((C_2^*)^{\oplus 2}\oplus D_2^*\bigr)\boxtimes\sho_{\bp^2}(-1) \\
\shE^*\otimes\shF & \delta^{-2}\otimes D_2\otimes\bigl((C_2^*)^{\oplus 2}\oplus D_2^*\bigr)\boxtimes\sho_{\bp^2}(-1) \\
\bottomrule
\end{array}
\end{equation}
By the K\"unneth formula, the last two rows of \eqref{eq:second-boundary} are acyclic because $\RG(\bp^2,\sho_{\bp^2}(-1)) = 0$. For the first row, twisting the Euler sequence on $\bp^2$ gives
\[
0 \longrightarrow \sho_{\bp^2}(-3) \longrightarrow \sho_{\bp^2}(-2)^{\oplus 3} \longrightarrow \mathcal{Q}_2(-2) \longrightarrow 0,
\]
whence
\begin{equation}\label{eq:Q2-minus2-cohom}
H^0(\bp^2,\mathcal{Q}_2(-2)) = 0,\qquad H^1(\bp^2,\mathcal{Q}_2(-2)) \simeq \bc,\qquad H^2(\bp^2,\mathcal{Q}_2(-2)) = 0.
\end{equation}
Although $\mathcal{Q}_2(-2)$ has a nonzero $H^1$, its $H^0$ vanishes; hence the K\"unneth formula gives $H^0(\mathfrak{Z}_{2,2}\times\bp^2,\mathcal{C}_\lambda(K_5(\End(\shE)))) = 0$. 
By equation \eqref{eq:boundary-degree-target} and the fact that $c_\lambda = 6$, to exclude the contribution to $H^5(X,K_5(\End(\shE)))$  from the stratum $S_\tau$, where $\tau=((2,2),(1,2))$, the vanishing $$H^0(\mathfrak{Z}_{2,2}\times\bp^2,\mathcal{C}_\lambda(K_5(\End(\shE)))) = 0$$ is exactly what we want.
\subsection{\texorpdfstring{The unstable HN type $\tau = ((3,3),(0,1))$}{The boundary stratum tau = ((3,3),(0,1))}}\label{subsec:boundary-33-01}

For $\tau = ((3,3),(0,1))$, the slopes are $\mu(3,3) = 1/2$ and $\mu(0,1) = -3$. The corresponding primitive 1-ps is \[
\lambda(t) = \bigl(\mathrm{diag}(t,t,t),\,\mathrm{diag}(t,t,t,t^{-6})\bigr).
\]
In the range $p \le 5$, for $\mathcal{T} = \End(\shE)$, the equality $m_\lambda(K_p(\mathcal{T})) = 42$ occurs only for $p=5$.
By Remark \ref{gerbe}, we only need to calculate the cohomology of the critical weight subbundle $\mathcal C_\lambda(K_5(\End(\mathcal E)))$ over the quotient stack
$$[Z_\lambda^\circ/G_{(3,3)}\times G_{(0,1)}]\simeq \mathfrak{Z}_{3,3},$$
where $
\mathfrak{Z}_{3,3} = \left[\bigl(\Hom(\bc^3,\bc^3)\otimes W\bigr)^{\mathrm{ss}} \Big/ \bigl((GL(\bc^3)\times GL(\bc^3))/\bc^*\bigr)\right].
$

Let $(C_3,D_3)$ be the source and target spaces of the $(3,3)$ block and $b$ be the target line of the $(0,1)$ block. Setting $\delta = \det(D_3)\otimes\det(C_3)^{-1}$, we have
\[
\shF^*|_{Z_\lambda^\circ} = (\delta\otimes b)^{-1}\otimes C_3,\qquad \shE^*|_{Z_\lambda^\circ} = \bigl((\delta\otimes b)^{-1}\otimes D_3\bigr)\oplus\delta^{-1},
\]
where the first summand in $\shE^*|_{Z_\lambda^\circ}$ has weight $7$ and the second summand $\delta^{-1}$ has weight $0$. The weight-$7$ piece of $\End(\shE)$ is $D_3\otimes b^{-1}$, and
\[
\bigl(\Lambda^5\WW^*\bigr)_{35} = (\delta\otimes b)^{-5}\otimes\Lambda^5(C_3^{\oplus 2}\oplus D_3),\qquad \det\mathcal{N}_\lambda = b^6\otimes(\det C_3)^{-3}\otimes\det D_3.
\]
Multiplying these together yields the critical weight subbundle on the quotient stack $\mathfrak{Z}_{3,3}$:
\begin{equation}\label{eq:third-boundary}
\mathcal{C}_\lambda\bigl(K_5(\End(\shE))\bigr) = \delta^{-4}\otimes(\det C_3)^{-2}\otimes D_3\otimes\Lambda^5(C_3^{\oplus 2}\oplus D_3).
\end{equation}
By equation \eqref{eq:boundary-degree-target} and the fact that $c_\lambda = 6$, to exclude the contribution to $H^5(X,K_5(\End(\shE)))$ (which is the only cohomology group that we care about) from the stratum $S_\tau$, where $\tau=((3,3),(0,1))$, we only need to show that $H^0(\mathfrak{Z}_{3,3},\mathcal{C}_\lambda(K_5(\End(\shE))))=0$.

For the dimension vector $(3,3)$ with stability parameter $(4,-3)$, any destabilizing subrepresentation (with respect to the slope-stability) of dimension $(a',b')$ should satisfy $a' > b'$, and the locus of such representations has codimension at least $\min_{a'>b'}(3a'(3-b') - \dim(\mathrm{Gr}(a',3)) - \dim(\mathrm{Gr}(b',3))) = 7 \ge 2$ (attained at $(a',b')=(1,0)$ or $(a',b')=(3,2)$). Thus global sections of \eqref{eq:third-boundary} over $\mathfrak{Z}_{3,3}$ extend to equivariant global sections on the representation space $\Hom(C_3,D_3)\otimes W$. Under the scaling action $z\cdot\mathrm{Id}_{C_3}$ (with $D_3$ fixed), the factor $\delta^{-4}\otimes(\det C_3)^{-2}$ has weight $z^6$, and in the decomposition
\[
\Lambda^5(C_3^{\oplus 2}\oplus D_3) = \bigoplus_{r=2}^5 \Lambda^r(C_3^{\oplus 2})\otimes\Lambda^{5-r}D_3,
\]
the $r$-th summand has weight $z^r$, giving total fiber weight $z^{6+r}$ with $6+r \ge 8 > 0$. Since $z\cdot\mathrm{Id}_{C_3}$ acts on degree-$n$ polynomials on $\Hom(C_3,D_3)\otimes W$ with weight $z^{n}$ ($n \ge 0$), no nonzero invariant section exists, so $H^0(\mathfrak{Z}_{3,3},\mathcal{C}_\lambda(K_5(\End(\shE))))=0$.

\section{Projective models}
\subsection{Cohomology of restricted bundles}\label{sec:restriction}

Having established all necessary vanishings on the unstable strata in Section~\ref{sec:unstable-strata}, we now compute the cohomology of $E = \shE|_S$ and $F = \shF|_S$ on any smooth 2-dimensional zero locus $S = Z(\sigma)$ of $\sigma\in H^0(X,\WW) = H^0(X,\shF^{\oplus 2}\oplus\shE)$.

For any vector bundle $\mathcal{T}$ on $X$, tensoring the Koszul resolution of $\sho_S$ by $\mathcal{T}$ gives the following spectral sequence
\begin{equation}\label{eq:koszul-ss}
\mathsf{E}_1^{-p,q} = H^q\bigl(X,K_p(\mathcal{T})\bigr) \Longrightarrow H^{q-p}(S,\mathcal{T}|_S),\qquad 0 \le p \le 10,
\end{equation}
where $K_p(\mathcal{T}) = \mathcal{T}\otimes\Lambda^p\WW^*$. Since $\dim X = 12$ and $\det\WW^* = K_X$, we have $K_p(\mathcal{T})^*\otimes K_X \simeq K_{10-p}(\mathcal{T}^*)$, and Serre duality on $X$ gives
\begin{equation}\label{eq:serre-koszul}
H^q\bigl(X,K_p(\mathcal{T})\bigr)^* \simeq H^{12-q}\bigl(X,K_{10-p}(\mathcal{T}^*)\bigr).
\end{equation}

\begin{lemma}\label{lem:structure-sheaf}
For every smooth 2-dimensional zero locus $S = Z(\sigma)$, one has $\RG(S,\sho_S) \simeq \bc[0]\oplus\bc[-2]$ and $K_S \simeq \sho_S$. In particular, $S$ is a K3 surface.
\end{lemma}

\begin{proof}
For $1 \le p \le 5$, Table~\ref{tab:all-18-hn} shows that $K_p(\sho_X) = \Lambda^p\WW^*$ satisfies the strict inequality $m_\lambda(K_p(\sho_X)) < \eta_\lambda$ on every unstable HN stratum, so $H^{>0}(X,K_p(\sho_X)) = 0$ by Lemma~\ref{Teleman}. Moreover, under the action of the one-parameter subgroup $\gamma(t) = (\mathrm{diag}(t,t,t),\mathrm{Id}_{V_4})$ in Lemma~\ref{lem:ambient-sections}, $\WW^*\leftrightarrow (V_3\otimes\chi^{-1})^{\oplus2}\oplus(V_4\otimes \chi^{-1})$ has strictly positive $\gamma$-weights ($4$ and $3$), whereas $\bc[R]=\Sym^*(V_3\otimes V_4^*\otimes W^*)$ has nonnegative $\gamma$-weights. Consequently, $H^0(X,K_p(\sho_X)) = 0$ for all $p \ge 1$ because any such global section should come from the $\gamma$-invariant part of $\bc[R]\otimes \Lambda^p\bigl((V_3\otimes\chi^{-1})^{\oplus2}\oplus(V_4\otimes \chi^{-1})\bigl)$ by Hartogs's lemma, which is zero. Thus $K_p(\sho_X)$ is acyclic for $1 \le p \le 5$, and Serre duality \eqref{eq:serre-koszul} implies that $K_p(\sho_X)$ is also acyclic for $6 \le p \le 9$.

At $p=0$, $K_0(\sho_X) = \sho_X$ satisfies $\RG(X,K_0(\sho_X)) = \bc[0]$ (since $X$ is Fano), and at $p=10$, $K_{10}(\sho_X) = K_X$ satisfies $\RG(X,K_{10}(\sho_X)) = \bc[-12]$. Thus the only nonzero terms on the $\mathsf{E}_1$-page of \eqref{eq:koszul-ss} are $\mathsf{E}_1^{0,0} = \bc$ (total degree $0$) and $\mathsf{E}_1^{-10,12} = \bc$ (total degree $2$). Since no differential can connect total degrees $0$ and $2$, the spectral sequence degenerates at $\mathsf{E}_1$, proving $\RG(S,\sho_S) \simeq \bc[0]\oplus\bc[-2]$. Finally, adjunction and \eqref{eq:canonical} give $K_S \simeq (K_X\otimes\det\WW)|_S \simeq \sho_S$.
\end{proof}

\begin{lemma}\label{lem:sections-restricted}
We have the following short exact sequences:
\begin{align}
0 \longrightarrow \bc \longrightarrow H^0(X,\shE) \longrightarrow H^0(S,E) \longrightarrow 0,\label{eq:sections-rank-four}\\
0 \longrightarrow \bc^5 \longrightarrow H^0(X,\shF) \longrightarrow H^0(S,F) \longrightarrow 0.\label{eq:sections-rank-three}
\end{align}
Moreover, $h^0(S,E) = 9$, $h^0(S,F) = 10$ and $H^{>0}(S,E) = H^{>0}(S,F) = 0$.
\end{lemma}

\begin{proof}
For $\mathcal{T} \in \{\shE,\shF\}$ and $0 \le p \le 5$, the strict inequality $m_\tau(K_p(\mathcal{T})) < \eta_\tau$ holds for all eighteen unstable HN types in Table~\ref{tab:all-18-hn}, so $H^{>0}(X,K_p(\mathcal{T})) = 0$. By the $\gamma$-weight argument in Lemma~\ref{lem:structure-sheaf}, $H^0(X,K_p(\mathcal{T})) = 0$ for all $p \ge 2$, while for $p=0,1$, Lemma~\ref{lem:ambient-sections} gives
\[
\begin{array}{c|ccc}
\toprule
\mathcal{T} & \RG(X,K_0(\mathcal{T})) & \RG(X,K_1(\mathcal{T})) & \RG(X,K_p(\mathcal{T}))\ (2 \le p \le 5) \\
\midrule
\shE & \bc^{10}[0] & \bc[0] & 0 \\
\shF & \bc^{15}[0] & \bc^5[0] & 0 \\
\bottomrule
\end{array}
\]
For $6 \le p \le 10$, Serre duality \eqref{eq:serre-koszul} reduces the calculation to $K_{p'}(\shE^*)$ and $K_{p'}(\shF^*)$ with $0 \le p' = 10 - p \le 4$. By the $\gamma$-weight argument, $H^0(X,K_{p'}(\shE^*)) = H^0(X,K_{p'}(\shF^*)) = 0$. For higher cohomology, the equality $m_\lambda(K_{p'}) = \eta_\lambda$ occurs only on the stratum $\tau = ((1,1),(2,3))$ at $p' = 3, 4$, corresponding to the last two rows of \eqref{eq:boundary-table} with $k = 0, 1$, which are acyclic by Lemma~\ref{lem:small-cohom}. Hence $K_{p'}(\shE^*)$ and $K_{p'}(\shF^*)$ are acyclic for $0 \le p' \le 4$.

It follows that 
\[
 \RG(S,E)=[\bc\longrightarrow\bc^{10}],\qquad
 \RG(S,F)=[\bc^5\longrightarrow\bc^{15}],
\]
with nonzero terms placed in degrees $-1,0$.
Since a coherent sheaf on $S$ has no cohomology in degree $-1$, both differentials are injective, yielding \eqref{eq:sections-rank-four}--\eqref{eq:sections-rank-three} and the stated dimensions of cohomology.
\end{proof}

\begin{lemma}\label{lem:simples}
The restricted bundles $E$ and $F$ on $S$ are simple: $\Hom_S(E,E) \simeq \Hom_S(F,F) \simeq \bc$.
\end{lemma}

\begin{proof}
Since $H^0(X,\End(\shE)) \simeq H^0(X,\End(\shF)) \simeq \bc$ by Lemma~\ref{lem:ambient-sections}, the Koszul spectral sequence \eqref{eq:koszul-ss} shows that it suffices to establish the following vanishings:
\begin{equation}\label{eq:diagonal-vanishing}
H^p\bigl(X,K_p(\End(\shE))\bigr) = H^p\bigl(X,K_p(\End(\shF))\bigr) = 0,\qquad 1 \le p \le 10.
\end{equation}
For $p=1,2$, the strict inequality $m_\lambda(K_p) < \eta_\lambda$ holds on all unstable strata. For $p=3,4$, the equality holds only for $\tau = ((1,1),(2,3))$ with $k=0,1$, in which cases the critical weight subbundles in \eqref{eq:boundary-table} are acyclic by Lemma~\ref{lem:small-cohom}. Thus $K_p(\End(\shE))$ and $K_p(\End(\shF))$ are acyclic for $1 \le p \le 4$.

At $p=5$, the equality $m_\tau(K_5) = \eta_\tau$ can occur on all the three HN types $\tau = ((1,1),(2,3))$, $((2,2),(1,2))$, and $((3,3),(0,1))$:
\begin{itemize}
\item For $\End(\shF)$, the equality occurs only for $((1,1),(2,3))$ and $((2,2),(1,2))$, and the corresponding critical weight subbundles in \eqref{eq:boundary-table} and \eqref{eq:second-boundary} are acyclic, so $K_5(\End(\shF))$ is acyclic.
\item For $\End(\shE)$, by \eqref{eq:boundary-degree-target}, the vanishing $H^5(X,K_5(\End(\shE))) = 0$ follows from the vanishings of\\
$H^{6-c_\lambda}(\mathcal{C}_\lambda(K_5(\End(\shE))))$ established in Sections~\ref{subsec:boundary-11-23}--\ref{subsec:boundary-33-01}, namely $H^2(\bp^2\times Y_{2,3},\sho_{\bp^2}\boxtimes\mathcal{J}_4) = 0$ for $((1,1),(2,3))$ and $H^0 = 0$ for $((2,2),(1,2))$ and $((3,3),(0,1))$.
\end{itemize}
Finally, for $6 \le p \le 10$, self-duality of $\End(\shE)$ and $\End(\shF)$ together with Serre duality \eqref{eq:serre-koszul} identifies $H^p(X,K_p(\End\mathcal{T}))^*$ with $H^{p'+2}(X,K_{p'}(\End\mathcal{T}))$, where $0 \le p' = 10 - p \le 4$. For $1 \le p' \le 4$, we have already shown that $K_{p'}$ is acyclic, while for $p'=0$, the rigidity $H^2(X,\End(\shE)) = H^2(X,\End(\shF)) = 0$ follows from Teleman quantization on $X$ (see \cite{belmans2023rigidity}).
\end{proof}

\begin{lemma}\label{lem:crosshom}
On every smooth 2-dimensional zero locus $S = Z(\sigma)$, one has
\[
\RHom_S(E,F) \simeq W^*[0] \simeq \bc^3[0]\qquad\text{and}\qquad \Hom_S(F,E) = 0.
\]
\end{lemma}

\begin{proof}
We apply \eqref{eq:koszul-ss} to $\mathcal{T} = \shE^*\otimes\shF$. For $0 \le p \le 5$, the equality $m_\tau(K_p(\shE^*\otimes\shF)) = \eta_\tau$ occurs only for $((1,1),(2,3))$ ($p=3,4,5$) and $((2,2),(1,2))$ ($p=5$), in which cases the critical weight subbundles in \eqref{eq:boundary-table} and \eqref{eq:second-boundary} are acyclic by Lemma~\ref{lem:small-cohom} and \eqref{eq:Kprime}. Thus $H^{>0}(X,K_p(\shE^*\otimes\shF)) = 0$ for $0 \le p \le 5$, while Lemma~\ref{lem:ambient-sections} and the $\gamma$-weight argument give $H^0(X,K_0(\shE^*\otimes\shF)) \simeq W^*$ and $H^0(X,K_p(\shE^*\otimes\shF)) = 0$ for $1 \le p \le 5$.

For $6 \le p \le 10$, Serre duality \eqref{eq:serre-koszul} reduces to considering $K_{p'}(\shF^*\otimes\shE)$ with $0 \le p' = 10 - p \le 4$. Here the equality $m_\tau = \eta_\tau$ holds only for $((1,1),(2,3))$ and $p'=3,4$ ($k=0,1$), where the critical weight subbundles $\ell^{1-k}\boxtimes(\Lambda^2\shV\otimes\Lambda^k\smallD)$ in \eqref{eq:boundary-table} are acyclic by Lemma~\ref{lem:small-cohom}. Since $H^0(X,K_{p'}(\shF^*\otimes\shE)) = 0$ by the $\gamma$-weight argument, $K_p(\shE^*\otimes\shF)$ is acyclic for all $1 \le p \le 10$. The Koszul spectral sequence therefore shows that $\RHom_S(E,F) \simeq W^*[0]$ at $(p,q)=(0,0)$, and Serre duality on $S$ gives $\Hom_S(F,E) \simeq H^2(S,E^*\otimes F)^* = 0$.
\end{proof}


\subsection{Degree and Mukai vectors}\label{degreePolarization}
For a general section $\sigma \in H^0(X,\WW)$, Proposition~\ref{cor:gg-on-S} ensures that $S = Z(\sigma)$ is smooth of dimension $2$ if it is nonempty. In fact, nonemptiness follows immediately from Lemma~\ref{lem:structure-sheaf} (if $S = \emptyset$, the Koszul complex would be exact on $X$, contradicting $\RG(S,\sho_S) \simeq \bc[0]\oplus\bc[-2]$). By Lemma~\ref{lem:structure-sheaf}, the smooth surface $S = Z(\sigma)$ is connected ($h^0(S,\sho_S)=1$), satisfies $h^1(S,\sho_S)=0$ and has trivial canonical bundle $K_S \simeq \sho_S$; thus $S$ is a K3 surface, and $h = c_1(L|_S)$ is an ample class. 
The self-intersection number $h^2$ is equal to $ \int_S c_1(L|_S)^2 = \int_X c_1(L)^2 c_3(\shF)^2 c_4(\shE)=40,$ which can be computed by the presentation of the Chow ring of $X$ in \cite{Franzen}; we will adopt an alternative method to calculate this self-intersection number $h^2$ below.

For any coherent sheaf $T$ on a K3 surface $S$, we use the Mukai vector convention
\[
v(T) = \ch(T)\sqrt{\td(S)} = \bigl(r,\,c_1,\,r + c_1^2/2 - c_2\bigr),
\]
equipped with the Mukai pairing $\langle (r,c,t),(r',c',t') \rangle = c\cdot c' - rt' - r't$, so that Riemann--Roch formula reads
\[
\chi(T,T') = -\langle v(T),v(T') \rangle,\qquad \chi(T) = r + t.
\]
Lemma~\ref{lem:sections-restricted} gives $\chi(S,E) = h^0(S,E) = 9$ and $\chi(S,F) = h^0(S,F) = 10$. Consequently,
\[
v(E) = (4,h,5),\qquad v(F) = (3,h,7).
\]
From lemma~\ref{lem:crosshom}, we have on $S$:
\[
 \chi(E,F) = \chi(S,E^*\otimes F) = 3.
\]
Now Riemann--Roch formula applies to $E^*\otimes F$ and gives
\[
3 = \chi(E,F) = -\langle (4,h,5),(3,h,7) \rangle = 4\cdot 7 + 3\cdot 5 - h^2 = 43 - h^2,
\]
whence $h^2 = 43 - 3 = 40$. 
In particular,
\[
-\chi(E,E)=\langle v(E),v(E) \rangle = 40 - 40 = 0,\qquad -\chi(F,F)= \langle v(F),v(F) \rangle = 40 - 42 = -2.
\]
By Lemma~\ref{lem:simples} and Serre duality on $S$, both $h^0(S,\End E) = h^2(S,\End E)$ and $h^0(S,\End F) = h^2(S,\End F)$ are equal to $1$. Since $\chi(E,E) = 0$ and $\chi(F,F) = 2$, it follows that
\[
\dim\Ext^\bullet_S(E,E) = (1,2,1),\qquad \dim\Ext^\bullet_S(F,F) = (1,0,1).
\]
Together with Lemmas~\ref{lem:sections-restricted} and \ref{lem:crosshom}, this completes the proof of Proposition~\ref{prop:cohom-summary}. Note that all assertions of Proposition~\ref{prop:cohom-summary} hold on every smooth 2-dimensional zero locus $S = Z(\sigma)$, prior to establishing slope stability or the primitivity of $h$.

\subsection{The Kodaira--Spencer map}\label{sec:deformations}

Let $\mathscr{U}_{\mathrm{sec}} \subset H^0(X,\WW)$ denote the nonempty Zariski open subset of sections whose zero locus $S = Z(\sigma)$ is a smooth K3 surface.

\begin{proposition}\label{prop:KS}
At every point $\sigma \in \mathscr{U}_{\mathrm{sec}}$, the Kodaira--Spencer map of the universal family over $\mathscr{U}_{\mathrm{sec}}$ has rank $19$. Its image is the nineteen-dimensional subspace of first-order deformations of $S$ preserving the polarization class $h$. 
\end{proposition}

\begin{proof}
Fix $S = Z(\sigma)$ and restrict the ambient resolution sequence \eqref{eq:tangent-X} for $T_X$ to $S$. Setting $Q_0 = \mathrm{coker}(\sho_S \to \End(E)\oplus\End(F))$, the restricted sequence breaks into two short exact sequences of vector bundles on $S$:
\begin{align*}
0 &\longrightarrow \sho_S \longrightarrow \End(E)\oplus\End(F) \longrightarrow Q_0 \longrightarrow 0,\\
0 &\longrightarrow Q_0 \longrightarrow (E^*\otimes F)^{\oplus 3} \longrightarrow T_X|_S \longrightarrow 0.
\end{align*}
In the first sequence, $H^0(\sho_S) \to H^0(\End(E)\oplus\End(F))$ is the diagonal scalar inclusion. The map $H^2(\sho_S)\rightarrow H^2(\End(E)\oplus\End(F))$ is Serre dual to the trace map, which is surjective; hence $H^2(\sho_S) \to H^2(\End(E)\oplus\End(F))$ is injective. Using Proposition~\ref{prop:cohom-summary}, the long exact sequence gives
\[
h^0(S,Q_0) = 1,\qquad h^1(S,Q_0) = 2,\qquad h^2(S,Q_0) = 1.
\]
By Lemma~\ref{lem:crosshom}, we have the equality $\RG(E^*\otimes F)^{\oplus 3}=\bc^9[0]$. The long exact sequence associated to the second short exact sequence therefore yields
\[
0 \longrightarrow \bc \longrightarrow \bc^9 \longrightarrow H^0(S,T_X|_S) \longrightarrow \bc^2 \longrightarrow 0,
\]
together with $H^1(S,T_X|_S) \simeq \bc$ and $H^2(S,T_X|_S) = 0$; in particular, $h^\bullet(S,T_X|_S) = (10,1,0)$.

The normal bundle  of $S$ inside $X$ is $N_{S/X} \simeq \WW|_S = F^{\oplus 2}\oplus E$. The Kodaira--Spencer map at $\sigma \in \mathscr{U}_{\mathrm{sec}}$ factors as
\[
T_\sigma\mathscr{U}_{\mathrm{sec}} = H^0(X,\WW) \longrightarrow H^0(S,N_{S/X}) \xrightarrow{\ \partial\ } H^1(S,T_S),
\]
where $\partial$ is the connecting morphism appearing in the long exact sequence associated to the short exact sequence $0 \to T_S \to T_X|_S \to N_{S/X} \to 0$.
By \eqref{eq:sections-rank-four}--\eqref{eq:sections-rank-three}, the restriction map $H^0(X,\WW) \to H^0(S,N_{S/X})$ is surjective, so the image of the Kodaira--Spencer map is equal to $\image(\partial)$. Using $h^1(T_X|_S)=1$ and lemma~\ref{lem:sections-restricted}, we can calculate $\dim(\mathrm{im}(\partial))=19$. Since every fiber of the family carries the restriction of $L$, this image is equal to the nineteen-dimensional hyperplane $\ker\bigl(H^1(S,T_S) \xrightarrow{\cup h}H^2(S,\End(T_S))=H^2(S,\sho_S)\bigr)$.
\end{proof}

Notice that Proposition~\ref{prop:KS} does not assume that $h$ is primitive. By the local Torelli theorem for K3 surfaces \cite[Chapters~6 and~7]{Huybrechts}, the period image of $\mathscr{U}_{\mathrm{sec}}$ contains an analytic open subset of the nineteen-dimensional period domain. Since the Noether--Lefschetz loci where $\rk\Pic(S) \ge 2$ form a countable union of effective divisors, the image of $\mathscr{U}_{\mathrm{sec}}$ contains members $S$ of Picard rank $1$.

\subsection{Primitivity of the polarization and slope-stability}\label{sec:primitive}

We now combine global generation properties and simplicity of vector bundles, and the existence of Picard-rank-one members to show that $h$ is a primitive polarization on $S$ and deduce slope-stability of $E$ and $F$. The following lemma is well-known:

\begin{lemma}\label{lem:ggquotient}
Let $T$ be a simple globally generated vector bundle of $\rk(T) > 1$ on a connected projective variety with $H^0(\sho) = \bc$. Then $\Hom(T,\sho) = 0$. If the variety is a smooth projective surface $S$ with $\Pic(S) = \bz. A$, where $A$ is ample, then every torsion-free quotient $\mathcal{Q}$ of $T$ with $0 < \rk(\mathcal{Q}) < \rk(T)$ satisfies $c_1(\mathcal{Q}) = aA$ for some integer $a \ge 1$.
\end{lemma}

\begin{proof}
If there exists a nonzero morphism $\varphi : T \to \sho_S$, global generation property of $T$ implies that some global section of $T$ maps to a nowhere-vanishing global section under $\varphi$. Then $\sho_S$ will be a direct summand of $T$ and this contradicts the simplicity of $T$.

For any torsion-free quotient $\mathcal{Q}$ of $T$ with $0 < r = \rk(\mathcal{Q}) < \rk(T)$, $\mathcal{Q}$ is globally generated, and $\det\mathcal{Q} = \sho_S(aA)$ is effective with $a \ge 0$. If $a = 0$, choosing $r$ generically independent global sections of $\mathcal{Q}$ induces a bundle morphism $\sho_S^{\oplus r} \to \mathcal{Q}^{**}$ whose determinant is a nonzero constant, hence we have an isomorphism $\sho_S^{\oplus r} \xrightarrow{\sim} \mathcal{Q} \xrightarrow{\sim} \mathcal{Q}^{**}$. Projecting onto a direct summand $\sho_S$ of $Q$ produces a nonzero homomorphism $T \to \sho_S$, which contradicts the first assertion. Thus $a \ge 1$.
\end{proof}

\begin{proposition}\label{prop:no-imprimitive}
Let $S$ be a K3 surface with $\Pic(S) = \bz. A$ and $A^2 = 10$. Then $S$ admits no simple globally generated rank-four vector bundle $T$ with Mukai vector $v(T) = (4,2A,5)$ and $h^0(S,T) = 9$.
\end{proposition}

\begin{proof}
Suppose such a bundle $T$ exists. We first show that $T$ is $A$-slope-semistable, of slope $\mu_A(T) = (2A\cdot A)/4 = 5$. By Lemma~\ref{lem:ggquotient}, any nontrivial torsion-free quotient $\mathcal{Q}$ of $T$ has rank $r \in \{1,2,3\}$ and $c_1(\mathcal{Q}) = aA$ with $a \ge 1$. The inequality $\mu_A(\mathcal{Q}) = 10a/r < 5$ can hold only when $(r,a) = (3,1)$. In that case, the saturated kernel of $T \to \mathcal{Q}$ is a rank-one reflexive sheaf on a smooth surface, thus it is equal to the line bundle $\sho_S(A)$ and fits into an exact sequence
\[
0 \longrightarrow \sho_S(A) \longrightarrow T \longrightarrow \mathcal{Q} \longrightarrow 0.
\]
By Kodaira vanishing and Riemann--Roch formula, $H^1(S,\sho_S(A))=H^2(S,\sho_S(A)) = 0$ and $h^0(S,\sho_S(A))$
$ = 2 + A^2/2 = 7$. This forces $h^0(S,\mathcal{Q}) = 9 - 7 = 2$, which is impossible for a globally generated sheaf $\mathcal{Q}$ of rank $3$. Thus $T$ is $A$-slope-semistable.

Next, by Kuleshov's theorem \cite[Theorem~2.1]{Kuleshov}, $S$ carries an $A$-slope-semistable spherical vector bundle $B$ with Mukai vector $v(B) = (2,A,3)$. Since the slope of any rank-one subsheaf must be an integer multiple of $10$, the bundle 
$B$ is actually slope-stable with $\mu_A(B) = 5 = \mu_A(T)$ and
\[
\chi(B,T) = -\langle (2,A,3),(4,2A,5) \rangle = -(20 - 10 - 12) = 2.
\]
If $T$ is slope-stable, any nonzero map between $B$ and $T$ would be an isomorphism between stable bundles of the same slope, which cannot happen. Then $\Hom(B,T) = \Hom(T,B) = 0$, which implies $\chi(B,T) = -\dim\Ext^1(B,T) \le 0$ and leads to a contradiction.

Thus $T$ must be strictly $A$-slope-semistable. Any saturated Jordan--H\"older filtration of $T$ has two stable factors with rank $2$ and first Chern class $A$, whose second Chern classes sum to $c_2(T) - A^2 = 19 - 10 = 9$. For a stable rank-two sheaf with $c_1 = A$, the Mukai bound $v^2 = 4c_2 - 18 \ge -2$ gives $c_2 \ge 4$, so the two factors have $c_2 = 4$ and $c_2 = 5$ respectively. Passing to the double dual shows that the factor with $c_2 = 4$ is locally free and has Mukai vector $(2,A,3)$, hence this factor is isomorphic to the spherical bundle $B$. Meanwhile, the other factor $C$ has Mukai vector $v(C) = (2,A,2)$ and $\chi(C,B) = 0$. Moreover, $\Hom(B,C) = 0$: any nonzero morphism $B \to C$ must be injective but the numerical invariants of $B$ are bigger.

There are now two possibilities for the Jordan--H\"older filtration:
\begin{itemize}
\item If the filtration is given by $0 \to B \to T \to C \to 0$, then $\Ext^1(C,B) \neq 0$ because the simple bundle $T$ cannot be equal to the direct sum of $B$ and $C$. Since $\chi(C,B) = 0$ and $\Hom(B,C) = 0$, Serre duality gives $\dim\Hom(C,B) = \dim\Ext^1(C,B) > 0$. Any nonzero map $C \to B$ yields a nontrivial endomorphism $T \twoheadrightarrow C \to B \hookrightarrow T$, contradicting simplicity of $T$.
\item If the filtration is given by $0 \to C \to T \to B \to 0$, then $C$ is also locally free (as the kernel of a surjection between vector bundles). Because $B$ and $C$ are nonisomorphic stable bundles of the same slope, we must have $\Hom(B,C) = \Hom(C,B) = 0$. Then $\chi(B,C) = 0$ implies $\Ext^1(B,C) = 0$. The extension therefore splits, again contradicting simplicity of $T$.\qedhere
\end{itemize}
\end{proof}

\begin{proof}[Proof of Theorem~\ref{thm:main}]
By Proposition~\ref{cor:gg-on-S} and Proposition~\ref{prop:KS}, we may choose a member $S = Z(\sigma)$ with $\sigma \in \mathscr{U}_{\mathrm{sec}}$ such that $\Pic(S) = \bz. A$ (with $A$ ample) and $E, F$ are globally generated on $S$. Write $h = mA$ with $m \in \bz_{\ge 1}$. Since $40 = h^2 = m^2 A^2$, $m$ can only be equal to $1$ or $2$. If $m=2$, then $A^2 = 10$ and $E$ is a simple globally generated rank-four bundle on $S$ with $v(E) = (4,2A,5)$ and $h^0(S,E) = 9$, contradicting Proposition~\ref{prop:no-imprimitive}. Hence $m=1$, and $h$ is primitive on $S$.

Because the divisibility in the integral cohomology groups of a deformation family of manifolds is locally constant, $h$ is primitive for \emph{every} member of $\mathscr{U}_{\mathrm{sec}}$. The family over $\mathscr{U}_{\mathrm{sec}}$ therefore defines a classifying morphism $\mathscr{U}_{\mathrm{sec}} \to \mathscr{M}_{21}$ to the irreducible nineteen-dimensional moduli space of primitively polarized K3 surfaces of genus $21$. By Proposition~\ref{prop:KS}, this classifying morphism is dominant.

Finally, both $E$ and $F$ are slope-stable when $\Pic(S)=\bz.h$ by Lemma~\ref{lem:ggquotient}. Since slope-stability is an open condition in flat families, $E$ and $F$ are $h$-slope-stable on a nonempty Zariski open subset of $\mathscr{U}_{\mathrm{sec}}$ which still dominates $\mathscr{M}_{21}$.
\end{proof}

\begin{Corollary}\label{cor:unirationality}
The coarse moduli space $\mathscr{M}_{21}$ of primitively polarized K3 surfaces of genus $21$ is unirational.
\end{Corollary}

\begin{proof}
The base $\mathscr{U}_{\mathrm{sec}}$ is a nonempty Zariski open subset of the forty-dimensional vector space $H^0(X,\WW)$, hence it is a rational variety. By Theorem~\ref{thm:main}(iii), the classifying morphism $\mathscr{U}_{\mathrm{sec}} \to \mathscr{M}_{21}$ is dominant, which exhibits a dominant rational map $\bp(H^0(X,\WW)) \dashrightarrow \mathscr{M}_{21}$.
\end{proof}

Recall that for any K3 surface $S$, Beauville and Voisin \cite{BV04} showed that all intersections of divisors and $c_2(T_S)$ in $\operatorname{CH}^2(S)$ lie in the cyclic subgroup $\mathbb{Z}. c_S$ generated by the canonical $0$-cycle class $c_S$ of a point on a rational curve in $S$. Let $\mathscr{M}_g^\circ \subset \mathscr{M}_g$ be the dense open subset parameterizing primitively polarized K3 surfaces of genus $g\ge3$ with trivial automorphism group, and let $\mathcal{S}_g^\circ \to \mathscr{M}_g^\circ$ be the universal family of such K3 surfaces. O'Grady \cite[Section 5]{Franchetta} proposed the generalized Franchetta conjecture, which predicted that the restriction of any codimension-two cycle class in $\operatorname{CH}^2(\mathcal{S}_g^\circ)$ to the fiber over $[S]\in \mathscr M_g^\circ$ is a multiple of $c_S$. 
\begin{Corollary}\label{franchetta}
    O'Grady's generalized Franchetta conjecture holds for $g = 21$.
\end{Corollary}
\begin{proof}
    By \cite[Theorem 1.2]{lu}, it suffices to prove that the universal family $\mathcal{S}_g^\circ$ is rationally connected. In our setting for $g=21$, it suffices to prove that the universal family over $\mathscr{U}_{\mathrm{sec}}$ is rational. The moduli space $X$ is rational by \cite[Theorem 6.4]{rational}. For points $x\in X-Z_{\mathrm{deg}}$, the spaces of elements $\sigma\in H^0(X,\mathcal W)$ vanishing at $x$ form a vector bundle over $ X-Z_{\mathrm{deg}}$, and this vector bundle is birational to the universal family over $\mathscr{U}_{\mathrm{sec}}$, which allows us to conclude. 
\end{proof}

\section{The Fourier--Mukai partner}\label{sec:fm}

The construction of projective models inside the quiver moduli space $X$ endows a general genus $21$ K3 surface $(S,h)$ with two slope-stable vector bundles: the rank-four bundle $E$ with $v(E) = (4,h,5)$, which deforms in a two-dimensional moduli space, and the rank-three spherical bundle $F$ with $ v(F) = (3,h,7)$, which is rigid. This indicates that the projective geometry around a fixed genus 21 K3 surface is closely related to its Fourier-Mukai partner $S'=M_S(4,h,5)$. We will investigate the relation between $S$ and $S'$ in the following sections.

Throughout Sections~\ref{sec:fm}--\ref{sec:hilb-schemes}, we assume that $(S,h)$ is a primitively polarized K3 surface of genus $21$ with $\Pic(S) = \bz. h$, $h^2 = 40$ and lies in the open locus where Proposition~\ref{fm:normalization} applies. Although the Picard-number-one assumption only holds for very general polarized K3 surfaces, main statements such as Theorem~\ref{thm:intro-hilb-square} and Theorem~\ref{thm:intro-hilb} will also hold for general K3 surfaces (of genus 21) using deformation argument as in \cite[Section 4]{Mukai2}.

In slope comparisons, we will normalize slopes by dividing by $h^2 = 40$, so a torsion-free sheaf with Mukai vector $(r,dh,t)$ has normalized slope $d/r$. Recall also that any stable sheaf $T$ on a K3 surface is simple, so Serre duality gives
\begin{equation}\label{eq:fm-square-bound}
v(T)^2 = \dim\Ext^1(T,T) - 2 \ge -2.
\end{equation}

\subsection{\texorpdfstring{The moduli surface $S' = M_{S}(4,h,5)$ and Mukai duality}{The moduli surface S' = M(4,h,5) and local freeness}}

Let
\[
S' = M_{S}(4,h,5)
\]
denote the moduli space of Gieseker-semistable sheaves on $S$ with the primitive Mukai vector $v = (4,h,5)$. Because $\Pic(S) = \bz. h$ and $\gcd(4,1)=1$, every Gieseker-semistable sheaf of Mukai vector $(4,h,5)$ is $h$-slope-stable, and $S'$ is a smooth projective K3 surface. Since $\gcd(4,40,5) = 1$, $S'$ is a fine moduli space admitting a universal sheaf.

\begin{lemma}\label{fm:local-free}
Every stable sheaf on $S$ with Mukai vector $(4,h,5)$ or $(3,h,7)$ is locally free. Moreover, there is a unique stable sheaf $F$ on $S$ with Mukai vector $(3,h,7)$.
\end{lemma}

\begin{proof}
If $T$ is a torsion-free slope-stable sheaf on $S$, its double dual $T^{**}$ is also slope-stable of the same rank and determinant. Putting $n = \operatorname{length}(T^{**}/T) \ge 0$, we have $v(T^{**}) = v(T) + (0,0,n)$ and
\[
v(T^{**})^2 = v(T)^2 - 2\rk(T)n.
\]
For $v(T) = (4,h,5)$ ($v^2 = 0$) or $v(T) = (3,h,7)$ ($v^2 = -2$), applying the inequality \eqref{eq:fm-square-bound} to $T^{**}$ forces $n = 0$, so $T$ is locally free.

The existence of some semistable sheaf $F$ with Mukai vector $(3,h,7)$ is guaranteed by \cite[Theorem 2.1]{Kuleshov}. Such $F$ is again automatically slope-stable and is locally free by the preceding paragraph.
If $F_1, F_2$ are nonisomorphic stable bundles with vector $(3,h,7)$, stability would force $\Hom(F_1,F_2) = \Hom(F_2,F_1) = 0$ because any such nonzero morphism must be an isomorphism, contradicting $\chi(F_1,F_2) = -\langle (3,h,7),(3,h,7) \rangle = 2 > 0$.
\end{proof}

Let $p : S\times S' \to S$ and $q : S\times S' \to S'$ denote the two projections. Mukai's duality theorem \cite[Proposition~1.1 and Theorem~1.2]{Mukai2} specializes to $(r,s)=(4,5)$ as follows.

\begin{proposition}\label{fm:normalization}
For a general polarized K3 surface $(S,h)$ of degree $40$, the moduli surface $S' = M_{S}(4,h,5)$ carries a primitive ample polarization $h'$ with $(h')^2 = 40$, and the rank-four universal bundle $\she$ on $S\times S'$ can be normalized so that
\begin{equation}\label{eq:fm-determinant}
\det\she \simeq p^*\sho_S(h)\otimes q^*\sho_{S'}(h').
\end{equation}
With this normalization, every horizontal slice $\she_s: = \she|_{\{s\}\times S'}$ ($s \in S$) is $h'$-slope-stable on $S'$ with Mukai vector $(4,h',5)$. This induces an isomorphism
\[
S \xrightarrow{\ \sim\ } M_{S'}(4,h',5),
\]
and the Fourier--Mukai functor $\Phi_{\she} = Rq_*(\she\otimes p^*(-)) : D^b(S) \xrightarrow{\sim} D^b(S')$ is an equivalence of derived categories with inverse $\Psi_{\she^*}[2] = Rp_*(\she^*\otimes q^*(-))[2]$.
\end{proposition}

Note that this Proposition doesn't require $\Pic(S)=\bz.h$. If we assume further that $\Pic(S)=\bz.h$, then we also have $\Pic(S')=\bz.h'$. For each $s'=[E] \in S'$, we have by definition $E \cong \she|_{S\times\{s'\}}$.

\begin{lemma}\label{fm:vanishing-lemma}
Let $S$ be a K3 surface with $\Pic(S) = \bz. h$ and $h^2 = 40$. If $A$ is a slope-stable vector bundle on $S$ of rank $r$ with first Chern class $c_1(A) = h$, then every nonsplit extension $B$ fitting into the short exact sequence
\begin{equation}\label{eq:extension-by-O}
0 \longrightarrow \sho_S \longrightarrow B \longrightarrow A \longrightarrow 0,
\end{equation}
is a slope-stable vector bundle.
\end{lemma}

\begin{proof}
Suppose a saturated subsheaf $C \subset B$ destabilizes $B$. Write $a = \rk(C) \le r$ and $c_1(C) = dh$. Since the normalized slope of $C$ satisfies $d/a \ge 1/(r+1) > 0$, we must have $d \ge 1$. Because $C \cap \sho_S$ has $c_1 \le 0$, the image of $C$ in $A$ must have positive first Chern class. If this image has rank strictly less than $r$, its normalized slope is at least $1/(r-1) > 1/r$, contradicting the slope-stability of $A$. The only remaining possibility is $a = r$, $C \cap \sho_S = 0$, and $d = 1$. Since $C$ is a saturated subsheaf of a vector bundle on a smooth surface, it is reflexive and hence locally free. The induced generically surjective morphism $C \to A$ between vector bundles of the same rank and determinant must be an isomorphism splitting \eqref{eq:extension-by-O}. This leads to a contradiction.
\end{proof}

\begin{lemma}\label{lem:vanishing-generation}
Let $A$ be a slope-stable vector bundle on $S$ with Mukai vector $v(A) = (r,h,t)$, where $\Pic(S) = \bz. h$ and $h^2 = 40$.
\begin{enumerate}[(i)]
\item If $(r+1)(t+1) > 21$, then $H^1(S,A) = H^2(S,A) = 0$ and $h^0(S,A) = r+t$.
\item If $(r+1)t > 21$, then $A$ is globally generated on $S$.
\end{enumerate}
\end{lemma}

\begin{proof}
(i) Slope-stability gives $\Hom(A,\sho_S) = 0$, so $H^2(S,A) = 0$ by Serre duality. If $H^1(S,A) \simeq \Ext^1(A,\sho_S)^* \neq 0$, then Lemma~\ref{fm:vanishing-lemma} produces a slope-stable extension $B$ with Mukai vector $v(B) = (r+1,h,t+1)$, whose square is
\[
v(B)^2 = 40 - 2(r+1)(t+1) < -2,
\]
contradicting \eqref{eq:fm-square-bound}. Riemann--Roch formula then gives $h^0(S,A) = \chi(S,A) = r+t$.

(ii) Suppose that global generation fails at some point $s \in S$, then there exists a surjection $A \twoheadrightarrow k(s)$ annihilating all global sections of $A$. Let $A_0 = \ker(A \twoheadrightarrow k(s))$. In the long exact sequence
\[
\Ext^1(A_0,\sho_S) \longrightarrow \Ext^2(k(s),\sho_S) \longrightarrow \Ext^2(A,\sho_S),
\]
the second arrow is Serre dual to the zero evaluation map $H^0(S,A) \to k(s)$, so the first arrow surjects onto $\Ext^2(k(s),\sho_S) \simeq \bc$. Notice that $\Ext^2(k(s),\sho_S)\cong \Ext^2_{\sho_{S,s}}(k(s),\sho_{S,s})$ is actually isomorphic to the space of local extensions $\Ext^1_{\sho_{S,s}}(A_{0,s},\sho_{S,s})$ by an argument using the long exact sequence. Choose an extension
\begin{equation}\label{eq:point-extension}
0 \longrightarrow \sho_S \longrightarrow B \longrightarrow A_0 \longrightarrow 0
\end{equation}
whose local extension class at $s$ is nonzero. Locally, $A_{0,s} \simeq \mathfrak{m}_s \oplus \sho_{S,s}^{\oplus(r-1)}$ and the nonzero local extension must resolve the direct summand $\mathfrak{m}_s$ via the Koszul complex $0 \to \sho_{S,s} \to \sho_{S,s}^{\oplus 2} \to \mathfrak{m}_s \to 0$; hence $B$ is locally free everywhere on $S$. A similar argument as in the proof of Lemma~\ref{fm:vanishing-lemma} shows that any saturated destabilizing subsheaf $C \subset B$ would be isomorphic to $A$ via a map factoring through the proper subsheaf $A_0 \subset A$, which is impossible. Thus $B$ is slope-stable of Mukai vector $v(B) = (r+1,h,t)$, satisfying $v(B)^2 = 40 - 2(r+1)t < -2$. Again \eqref{eq:fm-square-bound} leads to a contradiction.
\end{proof}


\begin{proposition}\label{fm:uniform}
For every $s'=[E]=[\she_{s'}] \in S'$ and every $s \in S$, the bundles $\she_{s'}$ on $S$ and $\she_s$ on $S'$ are globally generated, and satisfy
\[
h^0(S,\she_{s'}) = h^0(S',\she_s) = 9,\qquad H^{>0}(S,\she_{s'}) = H^{>0}(S',\she_s) = 0,\qquad \RHom_S(\she_{s'},F) \simeq \bc^3[0].
\]
Likewise, $F$ is globally generated on $S$ with $h^0(S,F) = 10$ and $H^{>0}(S,F) = 0$.
\end{proposition}

\begin{proof}
For $E=\she_{s'}$, $\she_s$ and $F$, Lemma~\ref{lem:vanishing-generation} establishes global generation, vanishing of higher cohomology, and the stated values of $h^0$.

Next, since $1/3 > 1/4$, slope-stability gives $\Hom(F,E) = 0$, hence $\Ext^2(E,F) = 0$. If $\Ext^1(E,F) \simeq \Ext^1(F,E)^* \neq 0$, we consider a nonsplit extension:
\[
0 \longrightarrow E \longrightarrow B \longrightarrow F \longrightarrow 0.
\]

By \cite[Lemma 2.1]{YoshiokaLemma}, $B$ must be slope-stable of Mukai vector $(7,2h,12)$, whose square is $160 - 168 = -8 < -2$, contradicting \eqref{eq:fm-square-bound}. Hence $\Ext^1(E,F) = 0$ and $\chi(E,F) = 3$ gives $\RHom_S(E,F) \simeq \bc^3[0]$.
\end{proof}

By Proposition~\ref{fm:uniform} and the cohomology-and-base-change theorem, the pushforwards
\begin{equation*}\label{eq:fm-pushforwards}
\mathscr{V} := p_*\she,\qquad \mathscr{V}' := q_*\she,\qquad (F')^* := q_*(\she^*\otimes p^*F)
\end{equation*}
are locally free of ranks $9$, $9$ and $3$ on $S$ and $S'$ respectively. All their higher direct images vanish, and the evaluation morphism $\mathscr{V}' \to \she_s$ obtained by pushing forward the evaluation map $\she\rightarrow\she|_{\{s\}\times S'}$ is surjective over $S'$ for every $s \in S$ (its fiber at $s' \in S'$ being the surjective evaluation map $H^0(S,\she_{s'}) \twoheadrightarrow \she_{s'}|_s=\she|_{(s,s')}$). Moreover, adjunction provides a canonical morphism on $S\times S'$:
\begin{equation}\label{fm:universal-map}
\epsilon : \she \longrightarrow p^*F\otimes q^*F' = F\boxtimes F'.
\end{equation}

\subsection{\texorpdfstring{Fourier-Mukai transform and spherical bundles}{The cohomological transform, stability of F' and V', and evaluation kernels}}
In order to investigate the bundles $\shv',F'$ and other Fourier-Mukai transforms, we first study the action of the Fourier-Mukai transform $\Phi_{\she}$ on the algebraic part of the Mukai lattice.
\begin{proposition}\label{fm:matrix}
The cohomological transform induced by $\Phi_{\she}$ on the algebraic part of the Mukai lattice sends $(r,dh,t)$ to $(r',d'h',t')$ according to the matrix below:
\begin{equation}\label{eq:fm-matrix}
\begin{pmatrix} r' \\ d' \\ t' \end{pmatrix}
=
\begin{pmatrix}
5 & 40 & 4 \\
1 & 9 & 1 \\
4 & 40 & 5
\end{pmatrix}
\begin{pmatrix} r \\ d \\ t \end{pmatrix}.
\end{equation}
In particular,
$v(\mathscr{V}') = (9,2h',9)$ and $v(F') = (3,h',7)$.
\end{proposition}

\begin{proof}
By Proposition~\ref{fm:normalization}, a skyscraper sheaf $k(s)$ ($v=(0,0,1)$) is mapped to $\she_s$ with $v(\she_s) = (4,h',5)$, giving the third column $(4,1,5)^{\mathsf{t}}$. On the other hand, by the symmetry between $S$ and $S'$ and the inverse of the Fourier--Mukai transform, we know that $\Phi_{\she}(4,-h,5)=(0,0,1)$. Because $\Phi_{\she}$ preserves the Mukai pairing, it must map $(0,-h,10)$ (which is orthogonal to both $(0,0,1)$ and $(4,-h,5)$) to $\pm(0,h',10)$. The image of $(0,-h,10)$ can only be $(0,h',10)$, because otherwise the image of $(1,0,0)$ is not an integer vector by calculation.

By definition $\Phi_{\she}(\sho_S) \simeq \mathscr{V}'[0]$, thus
applying \eqref{eq:fm-matrix} to $v(\sho_S) = (1,0,1)$ gives $v(\mathscr{V}') = (9,2h',9)$. By relative Serre duality and Proposition~\ref{fm:uniform}, $\Phi_{\she}(F^*) = Rq_*(\she\otimes p^*F^*) \simeq (q_*(\she^*\otimes p^*F))^*[-2] = F'[-2]$. Applying \eqref{eq:fm-matrix} to $v(F^*) = (3,-h,7)$ gives $v(F') = (3,h',7)$.
\end{proof}
Note that both $F'^*$ and $\mathscr V'$ are spherical by the fully-faithfulness of Fourier-Mukai transforms.
\begin{proposition}\label{fm:partner-bundle}
The rank-three bundle $F'$ on $S'$ is $h'$-slope-stable, spherical, and it is the unique stable bundle on $S'$ with Mukai vector $(3,h',7)$.
\end{proposition}

\begin{proof}
Suppose $(F')^*$ is not $h'$-slope-stable, and let $C \subset (F')^*$ be a semistable destabilizing subsheaf of rank $1$ or $2$. Its normalized slope is strictly greater than $\mu((F')^*) = -1/3$, hence $c_1(C)\ge 0$ and $\mu(C)\ge 0$. By definition of $(F')^* = q_*(\she^*\otimes p^*F)$, the inclusion $C \hookrightarrow (F')^*$ induces by adjunction a nonzero morphism $q^*C \to \she^*\otimes p^*F$ on $S\times S'$. Restricting to $\{s\}\times S'$ for a general $s \in S$ gives a nonzero map $C \to \she_s^*\otimes F_s \simeq (\she_s^*)^{\oplus 3}$. However, $\she_s^*$ is $h'$-slope-stable of normalized slope $-1/4 < 0 \le \mu(C)$, so nonzero map from $C$ to $(\she_s^*)^{\oplus 3}$ cannot exist and this leads to a contradiction. Thus $(F')^*$ and $F'$ are slope-stable, and the uniqueness of $F'$ follows from Lemma~\ref{fm:local-free}.
\end{proof}

We next prove that the rank-nine spherical bundle $\mathscr{V}' = q_*\she$ (as well as $\mathscr{V} = p_*\she$ on $S$) and the rank-five kernel $\Ker(\mathscr{V}' \to \she_s)$ are also slope-stable and globally generated.

\begin{lemma}\label{lem:finite-targets}
Let $C$ be a fixed slope-semistable torsion-free sheaf on a smooth projective surface. Among slope-stable vector bundles $B$ of the same slope as $C$, only finitely many such bundles up to isomorphism can satisfy $\Hom(C,B) \neq 0$.
\end{lemma}
\begin{proof}
Choose a Jordan--H\"older filtration of $C$ with torsion-free slope-stable subquotients. Any nonzero morphism $C \to B$ restricts to a nonzero morphism from one of these stable subquotients to $B$. Since a nonzero morphism between slope-stable torsion-free sheaves of equal slope must have full rank and is an isomorphism in codimension one, $B$ must be isomorphic to the reflexive hull (double dual) of that Jordan--H\"older factor, leaving only finitely many possibilities for $B$ up to isomorphism.
\end{proof}

\begin{proposition}\label{prop:V-stable}
The rank-nine bundle $\mathscr{V}' = q_*\she$ is $h'$-slope-stable on $S'$. By symmetry, the rank-nine bundle $\mathscr{V} = p_*\she$ is $h$-slope-stable on $S$ with Mukai vector $v(\mathscr{V}) = (9,2h,9)$.
\end{proposition}

\begin{proof}
Since $v(\mathscr{V}') = (9,2h',9)$ has normalized slope $2/9$ and $\gcd(2,9)=1$, no proper subsheaf of $\mathscr{V}'$ of rank $r \le 8$ can have normalized slope $2/9$. If $\mathscr{V}'$ were not slope-stable, the first factor in the Harder--Narasimhan filtration would be a slope-semistable subsheaf $C \subset \mathscr{V}'$ of rank $1 \le r \le 8$ and $c_1(C) = dh'$ with normalized slope $d/r > 2/9$.

By adjunction, the inclusion $C \hookrightarrow q_*\she$ corresponds to a nonzero morphism $q^*C \to \she$ on $S\times S'$, whose restriction to $\{s\}\times S'$ is some map from $C $ to $ \she_s$. Since $q^*C \to \she$ is nonzero, these restrictions cannot vanish for all $s \in S$. Stability of $\she_s$ (which has normalized slope $1/4$) therefore forces $d/r \le 1/4$. The only pairs of integers $(r,d)$ satisfying
\[
1 \le r \le 8,\qquad \frac{2}{9} < \frac{d}{r} \le \frac{1}{4}
\]
are $(r,d) = (4,1)$ and $(8,2)$, so $C$ is semistable of normalized slope $1/4$. However, the bundles $\she_s$ ($s \in S$) are pairwise nonisomorphic (as $\Hom(\she_{s_1},\she_{s_2}) = \Hom(k(s_1),k(s_2)) = 0$ for $s_1 \neq s_2$ by full faithfulness of $\Phi_{\she}$), so Lemma~\ref{lem:finite-targets} implies that $\Hom(C,\she_s) = 0$ for all but finitely many $s \in S$. Thus $q^*C \to \she$ vanishes on a dense open subset of $S\times S'$, hence this map vanishes everywhere since $\she$ is locally free, which leads to a contradiction.
\end{proof}

For each $s \in S$, Proposition~\ref{fm:uniform} shows that the evaluation map $\mathscr{V}' \to \she_s$ is surjective, defining a locally free rank-five kernel $K_s$ fitting into the exact sequence
\begin{equation}\label{eq:kernel-evaluation}
0 \longrightarrow K_s \longrightarrow \mathscr{V}' \longrightarrow \she_s \longrightarrow 0,\qquad v(K_s) = (9,2h',9) - (4,h',5) = (5,h',4).
\end{equation}
Note that $K_s = \Phi_{\she}(I_s)$ because $\Phi_\she(\sho_S)=\shv'[0]$ and $\Phi_\she(k(s))=\she_s[0]$.

\begin{proposition}\label{prop:K-stable}
For every $s \in S$, the rank-five bundle $K_s$ is $h'$-slope-stable, globally generated on $S'$, and satisfies $h^0(S',K_s) = 9$ and $H^{>0}(S',K_s) = 0$. Furthermore, $\mathscr{V}'$ is globally generated on $S'$ with $H^{>0}(S',\mathscr{V}') = 0$, $\mathscr{V}$ is globally generated on $S$ with $H^{>0}(S,\mathscr{V}) = 0$, and
\begin{equation}\label{eq:W-dimension}
U := H^0(S\times S',\she) = H^0(S',\mathscr{V}') = H^0(S,\mathscr{V})\qquad\text{has dimension } 18.
\end{equation}
\end{proposition}

\begin{proof}
A destabilizing subsheaf of $K_s$ of rank $1 \le a \le 4$ and $c_1 = dh'$ would have normalized slope $d/a \ge 1/5$, while its inclusion into the slope-stable bundle $\mathscr{V}'$ forces $d/a < 2/9$. This is impossible and shows that $K_s$ is $h'$-slope-semistable. The slope-stability of $K_s$ follows immediately.

An application of Lemma~\ref{lem:vanishing-generation} to $K_s$ with $(r,t) = (5,4)$ shows that $K_s$ is globally generated with $H^{>0}(S',K_s) = 0$ and $h^0(S',K_s) = 5+4 = 9$. Taking the associated long exact sequence of \eqref{eq:kernel-evaluation} yields
\[
0 \longrightarrow H^0(S',K_s) \longrightarrow H^0(S',\mathscr{V}') \longrightarrow H^0(S',\she_s) \longrightarrow 0,
\]
together with $H^{>0}(S',\mathscr{V}') = 0$ and $h^0(S',\mathscr{V}') = 9 + 9 = 18$. Since $K_s$ and $\she_s$ are globally generated and $H^0(S',\mathscr{V}') \twoheadrightarrow H^0(S',\she_s)$ is surjective, $\mathscr{V}'$ is globally generated as well. By the symmetry in Proposition \ref{fm:normalization}, the same results hold for $\mathscr{V} = p_*\she$ on $S$.
\end{proof}

\subsection{Remark on the projective model}
We have seen in the preceding section that the zero locus of a general global section of $\mathcal W$ is a genus 21 K3 surface $S$ together with two slope-stable vector bundles $E,F$ of ranks $4,3$ respectively. In view of the universal property of the quiver moduli $X=N(3;3,4)$, this embedding $S\hookrightarrow X$ is controlled by the evaluation map $F^*\otimes \Hom(F^*,E^*)\rightarrow E^*$, which can be viewed as a family of quiver representations of the 3-Kronecker quiver with dimension vector $(3,4)$. 

This slope-stable vector bundle $F$ with Mukai vector $(3,h,7)$ has been shown to be unique. For a general slope-stable bundle $[E]\in M_S(4,h,5)$, the evaluation map $F^*\otimes \Hom(F^*,E^*)\rightarrow E^*$ is a family of stable quiver representations by openness of stability for quiver representations. Thus the choice of such $[E]\in M_S(4,h,5)$ supplies a 2-dimensional family of embeddings of $S$ inside $X$ and explains geometrically the two excess parameters in the ambient parameter count:
\[
\dim H^0(X,\WW) - \dim\Aut(\WW) - \dim\operatorname{PGL}(W) = 40 - 11 - 8 = 21 = \dim\mathscr{M}_{21} + \dim M_S(4,h,5).
\]

However, it is not true that the evaluation map $F^*\otimes \Hom(F^*,E^*)\rightarrow E^*$ is a family of stable quiver representations for all $[E]\in M_S(4,h,5)$ in our genus 21 situation: the quiver-representation-stability requires that (pointwise on $S$) there is no subrepresentation with dimension vector $(1,1)$, which can be rephrased globally as follows: the family of one-dimensional subspaces of $F^*$ is parameterized by $\mathbb P_S(F^*)$ and the spaces $\Hom(F^*,E^*)$ are parametrized by $q_*(\she^*\otimes p^*F)=F'^*$ over $S'$. Thus (part of) the stability condition for quiver representation requires that the map 
$$\sho(-H)\boxtimes F'^*\subset F^*\boxtimes F'^*\rightarrow \she^*$$ should have rank at least 2 everywhere over $\bp_S(F^*)\times S'$, where $\sho(-H)$ is the tautological line bundle over the $\bp^2$-bundle $\bp_S(F^*)$. We can show that the rank 1 degeneracy locus of this bundle morphism under consideration is of codimension 6 and is a length 120 subscheme by Porteous's formula. 

This is in contrast with the genus 16 case. Recall that for a general genus 16 K3 surface $(S,h)$, the evaluation map $F_2^*\otimes\Hom(F_2^*,E_2^*)\rightarrow E_2^*$ is a family of stable quiver representations of dimension vector $(2,3)$ for the unique stable vector bundle $[F_2]\in M_S(2,h,8)$ and any $[E_2]\in M_S(3,h,5)$ (see \cite[Proposition 7.7]{Meng16}). In fact, the existence of such a subrepresentation with subdimension vector $(1,1)$ is again a codimension 6 condition, while the dimension of $\bp_S(F_2^*)\times M_S(3,h,5)$ is equal to 5 in that case.

\section{\texorpdfstring{Hyper-K\"ahler geometry of Hilbert squares and Hilbert cubes}{The geometric isomorphism of Hilbert squares}}\label{sec:hilb-schemes}

We now turn to the hyper-K\"ahler geometry of the Hilbert schemes of points on $(S,h)$ and its Fourier--Mukai partner $(S',h')$. In subsections 7.1-7.3, we establish an explicit isomorphism $S^{[2]} \xrightarrow{\sim} S'^{[2]}$ between their four-dimensional Hilbert squares and describe their common projective embedding into $\bp^{17}$. In subsections 7.4 and 7.5, we study the hyper-K\"ahler sixfolds $S^{[3]}\cong M_{S'}(3,h',6)$ and their common interesting four-dimensional subvariety $\bp_{S'}((F')^*)$.

As in Section~\ref{sec:fm}, we assume that $(S,h)$ is a primitively polarized K3 surface of genus $21$ with $\Pic(S) = \bz. h$, $h^2 = 40$ and lies in the open locus where Proposition~\ref{fm:normalization} applies. Although the Picard-number-one assumption only holds for very general polarized K3 surfaces,  main statements such as Theorem~\ref{thm:intro-hilb-square} and Theorem~\ref{thm:intro-hilb} will also hold for general K3 surfaces (of genus 21) using deformation argument as in \cite[Section 4]{Mukai2}.

\subsection{Length-two degeneracy loci}\label{sec:length-two}
Using the universal rank-four bundle $\she$ on $S\times S'$ and the rank-nine bundles $\mathscr{V} = p_*\she$ on $S$ and $\mathscr{V}' = q_*\she$ on $S'$, we will construct an explicit isomorphism $f : S^{[2]} \xrightarrow{\sim} S'^{[2]}$ by investigating the degeneracy loci of the restriction maps of $\she$, 
and construct an embedding of $S^{[2]}$ and $S'^{[2]}$ into $\bp(H^0(S\times S',\she)) \simeq \bp^{17}$.

For an arbitrary length 2 subscheme $Z\subset S$, we define $\she_Z := q_*(\she|_{Z\times S'})$, which is a rank 8 vector bundle over $S'$. The pushforward of the restriction map $\she\rightarrow\she|_{Z\times S'}$ gives us a bundle morphism between $\shv'$ and $\she_{Z}$, which is denoted by $e_Z$:
\begin{equation}\label{eq:main-evaluation}
e_Z : \mathscr{V}' \longrightarrow \she_Z.
\end{equation}
Before proceeding, note that for any $Z \in S^{[2]}$, the bundle $\she_Z = q_*(\she|_{Z\times S'})$ is an extension of two horizontal slices $\she_{s_1}$ and $\she_{s_2}$, so
\begin{equation}\label{eq:vector-PZ}
v(\she_Z) = (8,2h',10).
\end{equation}

The bundle morphism $\mathscr{V}' \xrightarrow{\ e_Z\ } \she_Z$ can be obtained by applying the Fourier--Mukai equivalence $\Phi_{\she}$ to the morphism $[\sho_S \to \sho_Z]$ over $S$. In order to study the kernel and cokernel of $e_Z$, we only need to investigate $\Phi_{\she}(I_Z)$.

Every length-two subscheme $Z \in S^{[2]}$ admits a filtration with two length-one factors: for a choice of $s_1 \in \Supp(Z)$, there is a short exact sequence of ideal sheaves on $S$:
\begin{equation}\label{eq:flag-ideals}
0 \longrightarrow I_Z \longrightarrow I_{s_1} \longrightarrow k(s_2) \longrightarrow 0,
\end{equation}
where $s_2\in S$. The equality $s_2= s_1$ holds when $Z$ is a double point supported at $s_1$ and the nonzero quotient $I_{s_1} \twoheadrightarrow k(s_1)$ records the tangent direction of $Z$ in that case.

Applying the Fourier--Mukai equivalence $\Phi_{\she}$ to \eqref{eq:flag-ideals}, and using $\Phi_{\she}(I_{s_1}) = K_{s_1}[0]$ (from \eqref{eq:kernel-evaluation}) and $\Phi_{\she}(k(s_2)) = \she_{s_2}[0]$, we obtain the following quasi-isomorphism:
\begin{equation}\label{eq:two-complexes}
\Phi_{\she}(I_Z) \simeq \bigl[K_{s_1} \xrightarrow{\ u_Z\ } \she_{s_2}\bigr] .
\end{equation}
Because $\Phi_{\she}$ is fully faithful and $I_{s_1} \twoheadrightarrow k(s_2)$ in \eqref{eq:flag-ideals} is nonzero, the morphism $u_Z : K_{s_1} \to \she_{s_2}$ is nonzero. Recall that $K_{s_1}$ and $\she_{s_2}$ are $h'$-slope-stable vector bundles with Mukai vectors $(5,h',4)$ and $(4,h',5)$ by Propositions~\ref{fm:normalization} and \ref{prop:K-stable}.

\begin{lemma}\label{lem:rank-five-four}
Let $A$ and $B$ be $h'$-slope-stable vector bundles on $S'$ with Mukai vectors $v(A) = (5,h',4)$ and $v(B) = (4,h',5)$ respectively. Then every nonzero homomorphism $u : A \to B$ is generically surjective and fits into an exact sequence on $S'$:
\begin{equation}\label{eq:rank-five-four}
0 \longrightarrow \sho_{S'} \longrightarrow A \xrightarrow{\ u\ } B \longrightarrow Q \longrightarrow 0,\qquad \length(Q) = 2.
\end{equation}
\end{lemma}

\begin{proof}
Let $J = \image(u) \subset B$ have rank $r \in \{1,2,3,4\}$ and first Chern class $c_1(J) = dh'$. If $r \le 3$, then $J$ is simultaneously a proper torsion-free quotient of $A$ and a low-rank proper subsheaf of $B$, so slope-stability of $A$ and $B$ yield
\begin{equation}\label{eq:impossible-interval}
\frac{1}{5} < \frac{d}{r} < \frac{1}{4},
\end{equation}
which cannot be achieved. Thus we must have $r = 4$.

The kernel $L = \ker(u)$ is a saturated rank-one subsheaf of the vector bundle $A$ on the smooth surface $S'$, hence it is a line bundle. Let $D \ge 0$ be the first Chern class of the torsion sheaf $\coker(u)$. Comparing determinants in $0 \to L \to A \to B \to \coker(u) \to 0$ gives $L \simeq \sho_{S'}(D)$. Since $D$ is effective and $\Pic(S') = \bz .h'$, we have $D \sim dh'$ with $d \ge 0$. Slope-stability of $A$ requires $d < 1/5$, which forces $d = 0$. Therefore $D = 0$, $L \simeq \sho_{S'}$, and $Q = \coker(u)$ is zero-dimensional. Simple calculations show that $\length(Q) = 2$.
\end{proof}

\begin{proposition}
    For every $Z \in S^{[2]}$, the degeneracy locus $D_7(e_Z) \subset S'$ defined by the $8\times 8$ maximal minors of \eqref{eq:main-evaluation} has length $2$.
\end{proposition}
\begin{proof}
Since $K_{s_1}$ and $\she_{s_2}$ are $h'$-slope-stable with Mukai vectors $(5,h',4)$ and $(4,h',5)$, the above lemma applies to $u_Z:K_{s_1} \rightarrow \she_{s_2}$ and establishes $\mathcal{H}^0(\Phi_{\she}(I_Z)) \simeq \sho_{S'}$ and $\mathcal{H}^1(\Phi_{\she}(I_Z)) \simeq Q_Z$ since $\Phi_\she(I_Z)$ is quasi-isomorphic to the complex $[K_{s_1} \rightarrow \she_{s_2}]$. Now applying $\Phi_\she$ to the short exact sequence $0\rightarrow I_Z\rightarrow \sho_S\rightarrow \sho_Z\rightarrow0$ gives
the exact sequence
\begin{equation}\label{eq:evaluation-long-exact}
0 \longrightarrow \sho_{S'} \xrightarrow{\ s_Z\ } \mathscr{V}' \xrightarrow{\ e_Z\ } \she_Z \longrightarrow Q_Z \longrightarrow 0,\qquad \length(Q_Z) = 2.
\end{equation}

 For the zero-dimensional sheaf $Q_Z$ on the smooth surface $S'$, $\HHom(Q_Z,\sho_{S'}) = \EExt^1_{\sho_{S'}}(Q_Z,\sho_{S'}) = 0$ and the only nonzero sheaf Ext into $\sho_{S'}$ is $\EExt^2_{\sho_{S'}}(Q_Z,\sho_{S'})$ with length $2$. Dualizing \eqref{eq:evaluation-long-exact} therefore yields the exact sequence
\begin{equation}\label{eq:dual-long-exact}
0 \longrightarrow \she_Z^* \xrightarrow{\ e_Z^*\ } (\mathscr{V}')^* \xrightarrow{\ s_Z^*\ } \sho_{S'} \longrightarrow \EExt^2_{\sho_{S'}}(Q_Z,\sho_{S'}) \longrightarrow 0.
\end{equation}
The last term is a length-two quotient of $\sho_{S'}$, hence it is equal to $\sho_{Z'}$ for a length-two subscheme $Z' \in S'^{[2]}$:
\begin{equation}\label{eq:ideal-resolution}
0 \longrightarrow \she_Z^* \xrightarrow{\ e_Z^*\ } (\mathscr{V}')^* \xrightarrow{\ s_Z^*\ } I_{Z'} \longrightarrow 0.
\end{equation}
In particular, the morphism $e_Z$ only degenerates at finitely many points and $Z'$ is the zero locus of the section $s_Z\in H^0(S',\mathscr{V}')$.

The kernel of the map $e_Z: \mathscr{V}' \rightarrow \she_Z$ is reflexive and thus locally free of rank $1$. In our situation where $\det\mathscr{V}' \simeq \det\she_Z$ and the degeneracy locus has codimension $2$, the kernel of the map $e_Z$ must be equal to $\sho_{S'}$.

Our section $s_Z$ obviously gives a morphism $\sho_{S'}\rightarrow \mathscr{V}'$ serving as the kernel map for $e_Z$. There is another map which can serve as the kernel map: taking the eighth exterior power of the map $e_Z: \mathscr{V}' \rightarrow \she_Z$, dualizing and tensoring with $\det\mathscr{V}' \simeq \Lambda^8\she_Z$ gives us a morphism $\sho_{S'}\rightarrow \mathscr{V}'$, whose image lies in the kernel of $e_Z$ by the cofactor identity. By the above paragraph these two maps must coincide with the kernel map up to a nonzero scalar, and in particular they have the same degeneracy locus when viewed as global sections of $\mathscr{V}'$.
Consequently,
\begin{equation}\label{eq:minors-equal-zero}
Z' = D_7(e_Z) = Z(s_Z) \in S'^{[2]},\qquad \length(Z') = 2.
\end{equation}
This proves the statement in the proposition and shows that \eqref{eq:dual-long-exact} coincides with the Eagon--Northcott complex \cite[Appendix A2.6]{Eisenbud} for the map $e^*_Z$.
\end{proof}

\subsection{Isomorphism between Hilbert squares}\label{sec:families}

Let $T = S^{[2]}$, let $p_{S'} : S'\times T \to S'$ and $p_T : S'\times T \to T$ denote the two projections, and let $\mathcal{Z} \subset S\times T$ be the universal length-two subscheme. Pulling back $\she$ from $S\times S'$ to $\mathcal{Z}\times S' \subset S\times T\times S'$ and pushing forward along the finite flat degree-two projection $\mathcal{Z}\times S' \to T\times S'$ produces the rank-eight universal vector bundle $\mathcal{B}$ on $T\times S'$. As in the last subsection, we have a universal evaluation morphism
\begin{equation}\label{eq:universal-evaluation}
e : p_{S'}^*\mathscr{V}' \longrightarrow \mathcal{B}.
\end{equation}
Let $D\mathcal{Z} := D_7(e) \subset S'\times T$ be the degeneracy locus defined by the $8\times 8$ maximal minors. Since the formation of maximal minors commutes with base change, the fiber of $D\mathcal{Z} \to T$ over every point $[Z] \in T=S^{[2]}$ is precisely the length-two scheme $Z' = D_7(e_Z)$ of \eqref{eq:minors-equal-zero}.



The map $D\mathcal{Z} \to T$ is proper with fibers being subschemes of $S'$ of constant length $2$, hence it is also flat. The universal property of the Hilbert scheme $S'^{[2]}$ gives us a global morphism
\begin{equation}\label{eq:f-morphism}
f : T = S^{[2]} \longrightarrow S'^{[2]},\qquad Z \longmapsto D_7(e_Z).
\end{equation}
In view of the symmetry between $S$ and $S'$ displayed in Proposition \ref{fm:normalization}, we have another map 
\begin{equation}
g : T' = S'^{[2]} \longrightarrow S^{[2]},\qquad Z' \longmapsto D_7(e_{Z'}),
\end{equation}
sending a length 2 subscheme $Z'\subset S'$ to the degeneracy locus of the corresponding morphism $e_{Z'}: \shv\rightarrow \she_{Z'}:=p_*(\she|_{S\times Z'})$.
\begin{proposition}\label{isoHilbertSquare}
    The assignment $f(Z) = D_7(e_Z)$ defines an isomorphism between smooth projective hyper-K\"ahler fourfolds
\[
f : S^{[2]} \xrightarrow{\ \sim\ } S'^{[2]},
\]
whose inverse $g : S'^{[2]} \xrightarrow{\sim} S^{[2]}$ is given by the exact same construction with $S$ and $S'$ exchanged.
\end{proposition}
\begin{proof}
To show that $f$ and $g$ are inverse to each other, we make the following observation. We adopt the notation $\DD_{S'}(A) = R\mathscr Hom_{\sho_{S'}}(A,\sho_{S'})$ for the derived dualization functor. Dualizing the two-term complex $\Phi_{\she}(I_Z) \simeq [\mathscr{V}' \xrightarrow{e_Z} \she_Z]$ and using the resolution \eqref{eq:ideal-resolution} gives a canonical isomorphism in $D^b(S')$:
\begin{equation}\label{eq:ideal-transform}
\DD_{S'}\Phi_{\she}(I_Z) \simeq I_{f(Z)}[0].
\end{equation}
Let $\Psi_{\she} = Rp_*(\she\otimes q^*(-))$ be the Fourier--Mukai transform from $ D^b(S')$ to $D^b(S)$, and let $\Psi_{\she^*}$ be $Rp_*(\she^*\otimes q^*(-))$. Since $K_S \simeq \sho_S$ and $K_{S'} \simeq \sho_{S'}$, the adjoint formula for Fourier--Mukai transforms and Grothendieck--Verdier duality give
\begin{equation}\label{eq:inverse-phi}
\Phi_{\she}^{-1} \simeq \Psi_{\she^{*}}[2],\qquad \DD_S\Psi_{\she}(A) \simeq \Psi_{\she^{*}}(\DD_{S'}A)[2].
\end{equation}
Combining the two identities in \eqref{eq:inverse-phi} yields the inverse of the contravariant equivalence $\DD_{S'}\circ\Phi_{\she}$:
\begin{equation}\label{eq:anti-inverse}
\bigl(\DD_{S'}\circ\Phi_{\she}\bigr)^{-1} \simeq \DD_S\circ\Psi_{\she}.
\end{equation}
By symmetry, we also have an isomorphism $\DD_S\Psi_{\she}(I_{Z'}) \simeq I_{g(Z')}[0]$ for all $Z' \in S'^{[2]}$. Equations \eqref{eq:ideal-transform} and \eqref{eq:anti-inverse} then give $I_{g(f(Z))} \simeq I_Z$ and $I_{f(g(Z'))} \simeq I_{Z'}$, which shows that $f $ and $g$ are inverse to each other.
\end{proof}


\subsection{\texorpdfstring{Embeddings of Hilbert squares}{Projective interpretations: 3-plane incidence, jet matrices, and closed immersions in P17}}\label{sec:projective}
We begin from fixing some notations. On the Hilbert square $T=S^{[2]}$, let $H_2$ denote the divisor class induced by $h$ and let $2\delta_2$ denote the exceptional divisor of the Hilbert--Chow morphism $S^{[2]} \to S^{(2)}$ (with $H'_2$ and $2\delta'_2$ defined analogously on $S'^{[2]}$). The Beauville--Bogomolov quadratic form \cite{Beauville} satisfies
\[
q(H_2) = 40,\qquad q(\delta_2) = -2,\qquad (H_2,\delta_2) = 0.
\]

In order to construct a projective embedding of $S^{[2]}\cong S'^{[2]}$, we come to exploit further the mechanism beneath the fact that $f$ and $g$ are inverse to each other in terms of global sections of $\she$.
Taking global sections in \eqref{eq:evaluation-long-exact} yields
\begin{equation}\label{eq:unique-product-section}
H^0\bigl(S\times S',\,\she\otimes p^*I_Z\bigr) = \ker\bigl(H^0(S',\mathscr{V}') \xrightarrow{H^0(e_Z)} H^0(S',\she_Z)\bigr) = \bc\cdot s_Z \simeq \bc.
\end{equation}
Let $\sigma_Z \in U = H^0(S\times S',\she)$ be a generator of the above space. It is a global section of $\she$ vanishing on $Z\times S'$; viewed as the section $s_Z$ of $\mathscr{V}' = q_*\she$ on $S'$, its zero scheme is $Z' = f(Z)$ by \eqref{eq:minors-equal-zero}. Therefore
\begin{equation}\label{eq:two-vanishings}
\sigma_Z|_{Z\times S'} = 0\qquad\text{and}\qquad \sigma_Z|_{S\times f(Z)} = 0.
\end{equation}
By the symmetry in Proposition \ref{isoHilbertSquare}, $\sigma_Z$ also spans the one-dimensional $\ker(H^0(e_{Z'}))$ associated to $Z' = f(Z)$ on the other side, and viewed as a section of $\mathscr{V} = p_*\she$ on $S$, its zero scheme is $g(Z') = Z$. 

Over $S'$, since $\she_Z$ is an extension of two horizontal bundles $\she_{s_1}, \she_{s_2}$, Proposition~\ref{fm:uniform} gives $h^0(S',\she_Z) = 18$ and $H^{>0}(S',\she_Z) = 0$ for every $Z \in T=S^{[2]}$. Thus $\mathcal{R} := (p_T)_*\mathcal{B}$ is a rank-eighteen vector bundle on $T$, and pushing forward \eqref{eq:universal-evaluation} along $p_T : S'\times T \to T$ induces a bundle morphism
\begin{equation}\label{eq:global-evaluation-family}
U\otimes\sho_T \longrightarrow \mathcal{R},
\end{equation}
whose fiberwise kernel is one-dimensional by \eqref{eq:unique-product-section}. Hence \eqref{eq:global-evaluation-family} has constant rank $17$, and its kernel $\mathcal{L} := \ker(U\otimes\sho_T \to \mathcal{R})$ is a line subbundle defining a morphism
\begin{equation*}\label{eq:j-definition}
j : T = S^{[2]} \longrightarrow \bp(U) \simeq \bp^{17},\qquad j^*\sho_{\bp(U)}(-1) = \mathcal{L}.
\end{equation*}
Exchanging $S$ and $S'$ defines $j' : S'^{[2]} \to \bp(U)$, and \eqref{eq:two-vanishings} implies $j = j'\circ f$.
\begin{thm}\label{closed-embedding}
The morphisms $j,j'$ are closed embeddings 
\[
j : S^{[2]} \hookrightarrow \bp(U) \simeq \bp^{17},\qquad j' : S'^{[2]} \hookrightarrow \bp(U) \simeq \bp^{17}
\]
satisfying $j = j'\circ f$. For each $Z \in S^{[2]}$, the point $j(Z) = [\sigma_Z]$ is the unique global section of $\she$, up to scale, vanishing identically on $Z\times S'$ (and simultaneously on $S\times f(Z)$). The pullbacks of $\sho_{\bp(U)}(1)$ are equal to $\sho_{S^{[2]}}(H_2 - 4\delta_2), \ \sho_{S'^{[2]}}(H'_2 - 4\delta'_2)$ respectively.
\end{thm}
\begin{proof}
To prove that $j$ and $j'$ are closed embeddings, note that the tautological line $\mathcal{L} \subset U\otimes\sho_T = H^0(S,\mathscr{V})\otimes\sho_T$ defines a section of $\mathscr{V}\boxtimes\mathcal{L}^*$ on $S\times T$. Let $\mathcal{Y} \subset S\times T$ be the zero scheme of this section. By construction, this section vanishes on the universal length-two family $\mathcal{Z} \subset S\times T$, so $\mathcal{Z} \subset \mathcal{Y}$ and $\mathcal{I}_{\mathcal{Y}} \subset \mathcal{I}_{\mathcal{Z}}$. Since $\sho_{\mathcal{Z}}$ is flat over $T$, the map $\mathcal{I}_{\mathcal{Z}}\otimes k(t) \to \mathcal{I}_{\mathcal{Z}_t}\subset\sho_{S\times\{t\}}$ is injective for every closed point $t \in T$.
Over any slice $\{t\}\times S$, we have the scheme-theoretical equality $\mathcal{Y}_t = \mathcal{Z}_t$ by the discussion above. Consequently, the map $\mathcal{I}_{\mathcal{Y}}\otimes k(t)\rightarrow \mathcal{I}_{\mathcal{Z}}\otimes k(t)$ is surjective, and $(\mathcal{I}_{\mathcal{Z}}/\mathcal{I}_{\mathcal{Y}})\otimes k(t) = 0$ for all $t \in T$. By Nakayama's lemma, $\mathcal{I}_{\mathcal{Y}} = \mathcal{I}_{\mathcal{Z}}$, which shows that the degeneracy locus of the section of $\mathscr{V}\boxtimes\mathcal{L}^*$ is exactly the universal length-two subscheme $\mathcal{Z} \subset S\times T$.

Now for any two morphisms from an arbitrary test scheme to $T = S^{[2]}$ whose compositions with $j : T \to \bp(U)$ are the same, they pull back the same section of $\mathscr{V}\boxtimes\mathcal{L}^*$. Hence they also pull back the same family of length-two subschemes by the above reasoning, and must coincide by the universal property of $S^{[2]}$. This shows that $j$ is a proper monomorphism and thus a closed embedding. The same holds for $j'$ and proves the statement of the theorem.

On $S'\times T$, by definition, the line bundle $\mathcal{L}$ fits into the following sequence of sheaves
\begin{equation}\label{eq:family-kernel}
0 \longrightarrow p_T^*\mathcal{L} \longrightarrow p_{S'}^*\mathscr{V}' \xrightarrow{\ e\ } \mathcal{B} \longrightarrow \mathcal{Q} \longrightarrow 0,
\end{equation}
whose pushforward via $p_T$ recovers \eqref{eq:global-evaluation-family} and whose restriction to every slice $S'\times\{Z\}$ is the exact sequence \eqref{eq:evaluation-long-exact}. In this sequence, $Q$ is defined as the cokernel and it is flat over $T$ since the fiberwise Hilbert polynomial is constant. 
An argument using the fiberwise exact sequence \eqref{eq:evaluation-long-exact} and Nakayama's lemma as above allows us to conclude that the above sequence is exact. Since $\Supp(\mathcal{Q}) = \mathcal{Z}'$ is finite over the fourfold $T$, it has codimension $2$ in $S'\times T$ and leads to $\det(\mathcal{Q}) \simeq \sho_{S'\times T}$. Consequently, $p_T^*\mathcal{L} \simeq \det(p_{S'}^*\mathscr{V}')\otimes(\det\mathcal{B})^{-1}$. Restricting to $\{s'\}\times T$, we obtain $\mathcal{L}=(\det\mathcal B|_{\{s'\}\times T})^{-1}$. The vector bundle $\mathcal{B}|_{\{s'\}\times T} $ over $T=S^{[2]}$ is actually equal to the tautological transform of the bundle $\she_{s'}$ over $S$. The determinant of this tautological transform is equal to $\sho_T(H_2-4\delta_2)$ by \cite[Corollary 2.4]{Wandel}, which finishes the proof.
\end{proof}

\begin{remark}
   Since $q(H_2 - 4\delta_2) = 40 - 32 = 8$ and $H_2 - 4\delta_2$ is ample, Riemann--Roch formula for $S^{[2]}$ and Kodaira vanishing give
\begin{equation*}\label{eq:twenty-one}
h^0\bigl(S^{[2]},\,\sho_{S^{[2]}}(H_2 - 4\delta_2)\bigr) = \chi\bigl(\sho_{S^{[2]}}(H_2 - 4\delta_2)\bigr) = \binom{q(H_2-4\delta_2)/2 + 3}{2} = \binom{7}{2} = 21,
\end{equation*}
 and the Fujiki formula shows that the degree of the embedded fourfold is $3\cdot 8^2 = 192$. 
\end{remark}

We conclude this subsection with an identification of the Picard group of $S^{[2]}\cong S'^{[2]}$.
From the equality $j = j'\circ f$ and Theorem \ref{closed-embedding}, the induced map $f^* : \Pic(S'^{[2]}) \xrightarrow{\sim} \Pic(S^{[2]})$ satisfies
\begin{equation*}\label{eq:invariant-polarization}
f^*(H'_2 - 4\delta'_2) = H_2 - 4\delta_2.
\end{equation*}
We restrict ourselves to the case where $\Pic(S)=\bz.h$ and $\Pic(S')=\bz.h'$. In this situation, the action of $f^*$ on the Picard lattice, which preserves the Beauville--Bogomolov form and fixes the positive vector $H_2 - 4\delta_2$, must be either the identity ($f^*H'_2 = H_2$, $f^*\delta'_2 = \delta_2$) or the reflection across $H_2 - 4\delta_2$. Because $H_2$ is the pullback of the ample generator on the symmetric square $S^{(2)}$, the associated ring $\bigoplus_{m \ge 0} H^0\big(S^{[2]}, \mathcal{O}_{S^{[2]}}(mH_2)\big)$ recovers the Hilbert--Chow contraction $S^{[2]} \to S^{(2)}$. If $f^*H'_2 = H_2$, the induced isomorphism of the associated rings would yield $S^{(2)} \simeq S'^{(2)}$, and comparing their singular diagonals would give $S \simeq S'$. However, under the Picard rank one assumption, $S$ and $S' = M_{S}(4,h,5)$ are nonisomorphic Fourier--Mukai partners using \cite[Proposition~2.1]{MMY}, because we know $S$ admits an isomorphism with $ M_{S}(1,h,20)$ by sending $s\in S$ to $I_s\otimes\sho_S(h)$. Consequently, $f^*$ cannot map $H_2'$ to $H_2$ and should be equal to the reflection across $H_2 - 4\delta_2$: 
\begin{equation*}\label{eq:picard-action}
f^*H'_2 = 9H_2 - 40\delta_2,\qquad f^*\delta'_2 = 2H_2 - 9\delta_2.
\end{equation*}

\subsection{\texorpdfstring{The $\bp^2$-bundle $Y'$ inside the Hilbert cube $S^{[3]}$}{The Brill-Noether P2-bundle Y' inside the Hilbert cube}}\label{sec:p2bundle}

In the previous subsection, the isomorphism between Hilbert squares $S^{[2]}$ and $ S'^{[2]}$ was constructed by considering the restriction maps of the rank-four universal bundle $\she$ alone. In the following, we will incorporate the rank-three spherical bundle $F$ and the three-dimensional arrow spaces $\Hom(\she_{s'},F)$ to study the six-dimensional Hilbert cube $S^{[3]}$.


\begin{lemma}\label{fm:rank-three}
For each $s'=[E] \in S'$, every nonzero homomorphism $u : E \to F$ is generically surjective and fits into an exact sequence on $S$:
\begin{equation}\label{eq:fm-arrow-length-three}
0 \longrightarrow \sho_S \longrightarrow E \xrightarrow{\ u\ } F \longrightarrow Q_u \longrightarrow 0,\qquad \operatorname{length}(Q_u) = 3.
\end{equation}
\end{lemma}

\begin{proof}
Assume $\operatorname{im}(u) \subset F$ has rank $r \in \{1,2,3\}$ and $c_1(\operatorname{im}(u)) = dh$. Since $E$ and $F$ are slope-stable of normalized slopes $1/4$ and $1/3$, we have the inequality $1/4 < d/r \le 1/3$, which admits no integer solution $(r,d)$ with $r \in \{1,2\}$. Thus $r = 3$ and $d = 1$. The kernel of $u$ is a rank-one reflexive sheaf with $c_1 = 0$, hence it is equal to $\sho_S$, and the cokernel $Q_u$ has zero-dimensional support with Mukai vector
\[
v(Q_u) = v(F) - v(E) + v(\sho_S) = (3,h,7) - (4,h,5) + (1,0,1) = (0,0,3),
\]
so $\operatorname{length}(Q_u) = 3$.
\end{proof}

Let
\[
\pi : Y' := \bp_{S'}\bigl((F')^*\bigr) \longrightarrow S'
\]
be the projective bundle over $S'$ whose fiber over $s'=[E] \in S'$ is $\bp\Hom_S(E,F)$. For each point $[u : E \to F] \in Y'$, dualizing the exact sequence in Lemma~\ref{fm:rank-three} gives
\begin{equation}\label{hilb:dual-sequence}
0 \longrightarrow F^* \xrightarrow{\ u^*\ } E^* \longrightarrow I_Z \longrightarrow 0,\qquad \operatorname{length}(Z) = 3,
\end{equation}
where $I_Z = \image(E^* \to \sho_S)$ is the ideal sheaf of a length-three subscheme $Z = \operatorname{Supp}(Q_u) \in S^{[3]}$. This construction works in families and defines a morphism $ Y' \to S^{[3]}$.

On $S^{[3]}$, let $F^{[3]}$ be the rank-nine tautological transform associated to $F$, and consider the evaluation morphism
\[
\ev_3 : H^0(S,F)\otimes\sho_{S^{[3]}} \simeq \sho_{S^{[3]}}^{\oplus 10} \longrightarrow F^{[3]},
\]
which is fiberwise given by $H^0(S,F)\rightarrow H^0(Z,F|_Z)$.
Let $\mathfrak{B} \subset S^{[3]}$ be the degeneracy locus defined by the $9\times 9$ maximal minors of $\ev_3$, whose closed points are
\[
\mathfrak{B} = \{Z \in S^{[3]} : h^0(S,F\otimes I_Z) \ge 2\}.
\]

\begin{thm}\label{hilb:bn}
The morphism $Y' = \bp_{S'}((F')^*) \to S^{[3]}$ is a closed embedding that identifies $Y'$ scheme-theoretically with $\mathfrak{B}$. In particular, $\mathfrak{B}$ is smooth of dimension $4$, and $h^0(S,F\otimes I_Z) = 2$ for every $Z \in \mathfrak{B}$.
\end{thm}

\begin{proof}
For any $Z \in S^{[3]}$, Riemann--Roch formula, stability of $F^*$ and Serre duality show that $\chi(I_Z,F^*) $
$= 1$, 
$\Hom(I_Z,F^*) = 0$ and $\Ext^2(I_Z,F^*) \simeq H^0(S,F\otimes I_Z)^*$, so
\begin{equation}\label{hilb:ext-dimension}
\dim\Ext^1_S(I_Z,F^*) = h^0(S,F\otimes I_Z) - 1.
\end{equation}
Consequently, the image of $Y'$ must lie inside $\mathfrak{B}\subset S^{[3]}$ because the nontrivial extension \eqref{hilb:dual-sequence} shows that $\dim\Ext^1_S(I_Z,F^*)\ge1$.

Conversely, for $Z \in \mathfrak{B}$, choose a nontrivial extension
\begin{equation}\label{hilb:reverse-extension}
0 \longrightarrow F^* \longrightarrow B \longrightarrow I_Z \longrightarrow 0.
\end{equation}
We claim that $B$ is $h$-slope-stable of normalized slope $-1/4$. In fact, any saturated destabilizing subsheaf $C \subset B$ has rank $\le 3$ and $c_1(C) = dh$ with $d \ge 0$. Since $F^*$ is stable of slope $-1/3$, the image of $C$ in $I_Z$ must have rank $1$ with nonpositive first Chern class. If $\rk(C) \ge 2$, then $C \cap F^*$ has rank $1$ or $2$ and $c_1(C \cap F^*) \le -h$ by the stability of $F^*$, contradicting $d \ge 0$. Thus $\rk(C) = 1$, $c_1(C) = 0$ and $C \xrightarrow{\sim} I_W \subset I_Z$ splits the pullback of \eqref{hilb:reverse-extension} to $I_W$ for some 0-dimensional subscheme $W$ containing $Z$. However, $\Ext^1(I_Z,F^*) \to \Ext^1(I_W,F^*)$ is injective because $\Ext^1(I_Z/I_W,F^*) = 0$ for the finite-length sheaf $I_Z/I_W$, which contradicts the choice of a nonsplit extension.

Since $v(B) = (4,-h,5)$ is isotropic with respect to the Mukai pairing, a similar argument as in Lemma~\ref{fm:local-free} shows that $B$ is locally free. Thus $B \simeq E^*$ for a unique bundle $[E] \in S'$. Applying $\Hom(F^*,-)$ to \eqref{hilb:reverse-extension} and using $\Hom(F^*,F^*) = \bc$, $\Hom(F^*,E^*) = \Hom(E,F) \simeq \bc^3$ and $\Ext^1(F^*,F^*) = 0$ yields
\[
0 \longrightarrow \bc \longrightarrow \bc^3 \longrightarrow H^0(S,F\otimes I_Z) \longrightarrow 0,
\]
whence $h^0(S,F\otimes I_Z) = 2$ and $\dim\Ext^1(I_Z,F^*) = 1$. Thus the extension \eqref{hilb:reverse-extension} is unique up to isomorphism and $\ev_3$ has constant rank $8$ along $\mathfrak{B}$. This shows that the set-theoretical map from $Y'$ to $\mathfrak{B}$ is bijective. 

In the characteristic zero setting, to show that $Y'\cong \mathfrak{B}$, it remains to show that $\mathfrak{B}$ is smooth. For a point $[Z]\in \mathfrak{B}$, the kernel and cokernel of the fiber map $H^0(S,F)\rightarrow H^0(Z,F|_Z)$ are equal to $H^0(F\otimes I_Z)$ and $H^1(F\otimes I_Z)\cong \Ext^1(I_Z,F^*)^*=\bc$ respectively. To show that the differentials of the maximal minors are nonzero at $[Z]\in \mathfrak{B}$, we only need to show that the following map is surjective:
$$\Ext^1(I_Z,I_Z)\rightarrow \Hom\bigl(H^0(F\otimes I_Z),H^1(F\otimes I_Z)\bigr),\qquad \xi\rightarrow [f\rightarrow\xi\cup f],$$
where $\Ext^1(I_Z,I_Z)$ is the tangent space at $[Z]$ in $S^{[3]}$.
Let $\alpha\in \Ext^1(I_Z,F^*)=\Hom(I_Z,F^*[1])$ denote the extension class corresponding to the extension sequence \eqref{hilb:dual-sequence}, which is actually equal to the generator of $\Ext^1(I_Z,F^*)\cong H^1(F\otimes I_Z)^*=\bc.$ Under this notation, the above map is Serre dual to the map 
$$\Hom(F^*,I_Z)\rightarrow \Hom(I_Z,I_Z[1]), \qquad f\rightarrow f[1]\circ\alpha.$$
 By applying $\Hom(-,I_Z)$ to the sequence \eqref{hilb:dual-sequence}, we only need to prove that $\Hom(E^*,I_Z)=\bc$ in order to conclude that the last map is injective. This last claim $\Hom(E^*,I_Z)=\bc$ can be obtained by applying $\Hom(E^*,-)$ to \eqref{hilb:dual-sequence} and invoking Proposition \ref{fm:uniform}.
\end{proof}

\subsection{\texorpdfstring{The Hilbert cube $S^{[3]}$ as a moduli space on $S'$}{The entire Hilbert cube as a moduli space on S'}}
In this subsection, we construct an isomorphism between $S^{[3]}$ and $M_{S'}(3,h',6)$ which is compatible with the embeddings of $Y'=\bp_{S'}((F')^*)$.\label{sec:isomorphicCube}
\begin{thm}\label{hilb:global-fm}
For every $Z \in S^{[3]}$, the shifted Fourier--Mukai transform
\[
K_Z := \Phi_{\she}(I_Z)[1]
\]
is an $h'$-slope-stable torsion-free sheaf on $S'$ of Mukai vector $(3,h',6)$, and $\Phi_{\she}(-)[1]$ induces an isomorphism of smooth projective hyper-K\"ahler sixfolds
\[
S^{[3]} \xrightarrow{\ \sim\ } M_{S'}(3,h',6).
\]
Under this isomorphism, the subvariety $Y' \simeq \mathfrak{B} \subset S^{[3]}$ is precisely the locus of sheaves $K_Z$ that are not locally free, and its points correspond to the elementary modifications
\[
0 \longrightarrow K_Z \longrightarrow F' \longrightarrow k(s') \longrightarrow 0.
\]
\end{thm}

\begin{proof}
First, for any $s'=[E] \in S'$ and any nonzero global section $\sigma_0 \in H^0(S,E)$, stability of $E$ and $\Pic(S) = \bz. h$ imply that $\sigma_0$ has no divisorial zero locus. Thus $\sigma_0$ realizes $\sho_S$ as a saturated subsheaf of $E$ and the cokernel $C = \mathrm{coker}(\sho_S \xrightarrow{\sigma_0} E)$ is torsion-free. The Mukai vector of $C$ is equal to $(3,h,4)$ and one can prove that $C$ is slope-stable by using the slope-stability of $E$.

Writing $n = \operatorname{length}(C^{**}/C)$ for the length of the zero scheme of $\sigma_0$, the bound \eqref{eq:fm-square-bound} for $C^{**}$ gives $16 - 6n = v(C^{**})^2 \ge -2$, so $n \le 3$. If the equality $n = 3$ holds, then $C^{**}$ must be equal to $ F$ and the zero locus of $\sigma_0$ is equal to the degeneracy locus of the induced map $E\rightarrow F$. Combining these with the discussion in the last subsection, for any length 3 subscheme $Z\subset S$, it follows that $H^0(S,E\otimes I_Z) \neq 0$ if and only if $Z \in Y'\subset S^{[3]}$ and $s' =[E]= \pi(Z)$, in which case $h^0(S,E\otimes I_Z) = 1$.

Now for a fixed length 3 subscheme $[Z] \in S^{[3]}$, we will define a stable sheaf $[K_Z]\in M_{S'}(3,h',6)$, which gives us the desired isomorphism. We apply $\Phi_{\she}$ to $0 \to I_Z \to \sho_S \to \sho_Z \to 0$. Using $\Phi_{\she}(\sho_S) = \mathscr{V}'[0]$ (rank $9$) and $\Phi_{\she}(\sho_Z) = \she_Z[0]$ where $\she_Z = q_*(\she|_{Z\times S'})$ (rank $12$), the preceding fiberwise calculation shows that the evaluation map $\mathscr{V}' \to \she_Z$ has injective fiber maps outside at most one single point $s'$, which is equal to $[E] = \pi(Z)\in S'$ when $[Z]\in Y'\subset S^{[3]}$. Thus $\mathscr{V}' \to \she_Z$ is injective as a sheaf morphism, and its cokernel $K_Z := \Phi_{\she}(I_Z)[1]$ sits in the short exact sequence
\begin{equation}\label{hilb:resolution}
0 \longrightarrow \mathscr{V}' \longrightarrow \she_Z \longrightarrow K_Z \longrightarrow 0.
\end{equation}
The sheaf $K_Z$ has depth $\ge 1$ at every closed points by Auslander--Buchsbaum formula, and has Mukai vector $v(K_Z) = (3,h',6)$ by \eqref{eq:fm-matrix}. Since  $K_Z$ is locally free outside at most one point, the bound on its local depth implies that it is torsion-free.

The torsion-free sheaf $K_Z$ is slope-stable: because $\she_Z$ is an extension of three slices $\she_{s_i}$, each of which is stable of normalized slope $1/4$, $\she_Z$ is semistable of normalized slope $1/4$. Any saturated destabilizing subsheaf of the rank-three sheaf $K_Z$ would have $c_1 \ge h'$, so its quotient would have rank $\le 2$ and normalized slope $\le 0$,  which is impossible when the quotient is viewed as a torsion-free quotient of $\she_Z$. 

Since the first nonzero map in \eqref{hilb:resolution} remains injective on every fiber over $S^{[3]}$, $K_Z$ forms a flat family of stable sheaves by the local criterion for flatness, and defines a morphism $S^{[3]} \to M_{S'}(3,h',6)$ between smooth irreducible sixfolds. This map is injective by the fully-faithfulness of the Fourier-Mukai transform, hence an isomorphism (in characteristic 0).

Finally, the locally freeness of $K_Z$ is equivalent to requiring the injective map in \eqref{hilb:resolution} to be a bundle inclusion. Thus the above analysis shows that $K_Z$ is locally free exactly when $[Z]\notin Y'$. 
In the case $[Z]\in Y'$, the sheaf $K_Z=\Phi_{\she}(I_Z)[1]$ can be calculated using the short exact sequence \eqref{hilb:dual-sequence}. Recalling that $\Phi_\she(F^*)=F'[-2]$ and $\Phi_\she(E^*)=k(s')[-2]$, we see that $K_Z$ is equal to the kernel of the map $F'\rightarrow k(s')$ and its non-locally-free locus is exactly $s'=[E]\in S'$.
\end{proof}

\appendix

\section{The full Harder--Narasimhan weight certificate}\label{app:hn}

The following table lists the weights of $\mathcal U_3,\mathcal U_4$, the codimension $c_\tau = \mathrm{codim}_R(S_\tau)$ of the strata and the weight $\eta_\tau = \wt_\tau(\det\mathcal{N}_\tau^*)$ for all the eighteen unstable Harder-Narasimhan types. All these calculations can be done using the software QuiverTools \cite{BFP:QuiverTools}.

\begin{center}\small
\begin{longtable}{@{}lrrrr@{}}
\toprule
HN type&$c_\tau$&$\mathcal U_3$ weights&$\mathcal U_4$ weights&$\eta_\tau$\\\midrule
\endfirsthead
\toprule HN type&$c_\lambda$&$\mathcal U_3$ weights&$\mathcal U_4$ weights&$\eta_\lambda$\\\midrule
\endhead
$(10,11,12,01)$&16&$(70,49,42)$&$(49,42,42,28)$&385\\
$(10,11,13)$&15&$(42,28,21)$&$(28,21,21,21)$&231\\
$(10,21,03)$&25&$(56,49,49)$&$(49,35,35,35)$&406\\
$(10,22,02)$&18&$(28,21,21)$&$(21,21,14,14)$&168\\
$(10,23,01)$&13&$(56,35,35)$&$(35,35,35,21)$&294\\
$(10,24)$&10&$(14,7,7)$&$(7,7,7,7)$&70\\
$(11,23)$&4&$(7,0,0)$&$(7,0,0,0)$&28\\
$(20,11,03)$&28&$(42,42,35)$&$(35,28,28,28)$&322\\
$(20,12,02)$&24&$(56,56,42)$&$(42,42,35,35)$&406\\
$(20,13,01)$&22&$(70,70,49)$&$(49,49,49,42)$&504\\
$(21,11,02)$&17&$(56,56,49)$&$(56,49,28,28)$&385\\
$(21,12,01)$&14&$(21,21,14)$&$(21,14,14,7)$&133\\
$(21,13)$&13&$(14,14,7)$&$(14,7,7,7)$&91\\
$(22,12)$&6&$(7,7,0)$&$(7,7,0,0)$&42\\
$(30,04)$&36&$(28,28,28)$&$(21,21,21,21)$&252\\
$(31,03)$&24&$(21,21,21)$&$(21,14,14,14)$&168\\
$(32,02)$&14&$(14,14,14)$&$(14,14,7,7)$&98\\
$(33,01)$&6&$(7,7,7)$&$(7,7,7,0)$&42\\
\bottomrule
\end{longtable}
\end{center}
\bibliographystyle{alpha}
\bibliography{bibliography}
\texttt{jmeng.maths@gmail.com} \\
Institut de Math\'ematiques de Toulouse; UMR 5219, Universit\'e de Toulouse; CNRS, UPS, F-31062 Toulouse Cedex 9, France

\end{document}